\documentclass[10pt,preprint]{elsarticle}
\usepackage[width=16cm,height=20cm]{geometry}
\usepackage{graphicx}
\usepackage{amssymb}
\usepackage{amsmath}
\usepackage{amsthm}
\usepackage{bm}
\usepackage{color}
\usepackage[bookmarks]{hyperref}
\usepackage{lineno}
\usepackage{aliases}

\newcommand{\vel}{\bm{v}}
\newcommand{\p}{\mathcal{P}}
\newcommand{\R}{\mathcal{R}}

\newcommand{\con}{\bm{U}}
\newcommand{\cons}{\con^{(1)}}
\newcommand{\conp}{\con^{(2)}}

\newcommand{\prim}{\bm{Q}}

\newcommand{\f}{\bm{F}}
\newcommand{\fx}{\bm{F}_1}
\newcommand{\fy}{\bm{F}_2}
\newcommand{\fs}{\f^{(1)}}
\newcommand{\fp}{\f^{(2)}}
\newcommand{\m}{\bm{m}}
\newcommand{\s}{\bm{S}}
\newcommand{\D}{\bm{D}}
\newcommand{\tD}{\tilde{\D}}
\newcommand{\B}{\bm{B}}
\newcommand{\Bp}{\B^{(2)}}
\newcommand{\A}{\bm{A}}

\newcommand{\dt}{\Delta t}
\newcommand{\dx}{\Delta x}
\newcommand{\dy}{\Delta y}

\newcommand{\uad}{\mathcal{U}}

\newcommand{\reva}[1]{#1}
\newcommand{\revb}[1]{#1}
\newcommand{\corr}[1]{#1}

\begin{document}

\begin{frontmatter}
\title{A path conservative finite volume method for a shear shallow water model}
\author[1]{Praveen Chandrashekar}
\ead{praveen@tifrbng.res.in}
\author[2]{Boniface Nkonga}
\ead{boniface.nkonga@unice.fr}
\author[1]{Asha Kumari Meena}
\ead{asha@tifrbng.res.in}
\author[2]{Ashish Bhole}
\ead{ashish.bhole@unice.fr}
\address[1]{TIFR Centre for Applicable Mathematics, Bangalore-560065, India.}
\address[2]{Universit\'e C\^ote d'Azur, Inria, CNRS, LJAD, 06108 Nice Cedex 2, France.}
\begin{abstract}
\revb{The shear shallow water model provides an approximation for shallow water flows by including the effect of vertical shear in the model.} This model can be derived from the depth averaging process by including the second order velocity fluctuations which are neglected in the classical shallow water approximation. The resulting model has a non-conservative structure which resembles the 10-moment equations from gas dynamics. This structure facilitates the development of path conservative schemes and we construct HLL, 3-wave and 5-wave HLLC-type solvers. An explicit and semi-implicit MUSCL-Hancock type second order scheme is proposed for the time integration. Several test  cases including roll waves show the performance of the proposed modeling and numerical strategy.
\end{abstract}
\begin{keyword}
Shear shallow water model, non-conservative system, path conservative scheme,  approximate Riemann solver, finite volume method.
\end{keyword}
\end{frontmatter}

\section{Introduction}
\begin{figure}
\begin{center}
\includegraphics[width=0.7\textwidth]{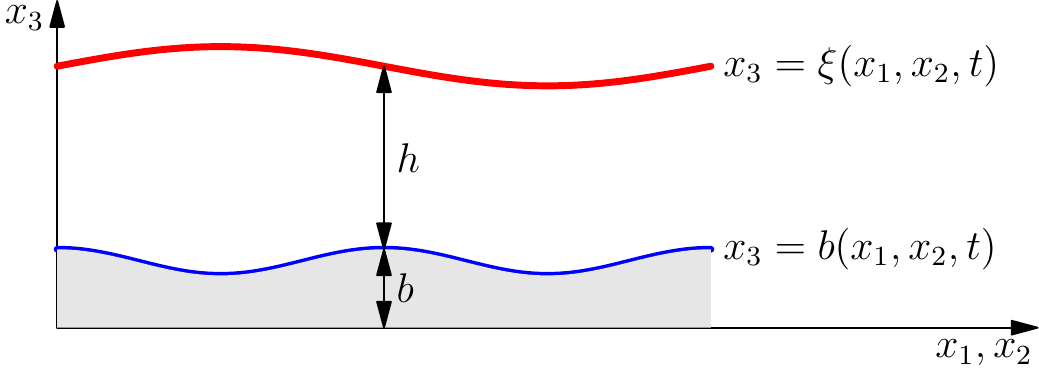}
\end{center}
\caption{Shallow water approximation: The free surface is given by $x_3 = \xi(x_1,x_2,t)$ and the bottom surface is given by $x_3 = b(x_1,x_2,t)$.}
\label{fig:shallow}
\end{figure}

The shallow water equations, sometimes also called the Saint-Venant equations, are used to model the flow of fluids in situations where the depth of fluid is small relative to the horizontal scale of the flow field variations~\cite{Whitham1999}. They have been used to model the flow on the scale of the  atmosphere and ocean,  and have been applied for tsunami prediction, storm surges, flow around structures, etc. Because the model has only two independent spatial variables and does not require tracking the free surface, it provides a simpler approximation than the full three dimensional Euler equations with a free surface. The shallow water equations are derived from the incompressible Euler or Navier-Stokes equations by averaging them over the depth coordinate. The horizontal velocity is assumed to be weakly varying in the vertical coordinate which implies that the vertical shear is negligible. This allows us to ignore the second order velocity fluctuations and leads to a closed set of equations which we may be called the {\em classical shallow water equations} and are given by
\begin{subequations}
\label{eq:csw}
\begin{align}
\df{h}{t} + \nabla\cdot(h\vel) &= 0 \\
\df{(h\vel)}{t} + \nabla\cdot \left( h \vel \otimes \vel + \frac{gh^2}{2} I \right) &= -g h \nabla b -C_f |\vel| \vel
\end{align}
\end{subequations}
where $h,\vel=(v_1,v_2)$ are the water depth and velocity, $b$ is the bottom topography, see Figure~(\ref{fig:shallow}), $g >0$ is the acceleration due to gravity, and we have included a frictional term to model the bottom friction with $C_f$ being the Chezy coefficient. The horizontal velocity $\vel$ is a depth average of the three dimensional velocity field. Since the vertical shear is neglected, the classical shallow water equations cannot model large scale eddies ('roller') that appear near the surface and behind the hydraulic jump. Under the assumption of smallness of horizontal vorticity, a more general model called the {\em shear shallow water} (SSW) model can be derived~\cite{Teshukov2007,Richard2012,Richard2013,Gavrilyuk2018} which includes the second order velocity fluctuation terms, and can be written as
\begin{subequations}
\label{eq:ssw}
\begin{align}
\df{h}{t} + \nabla\cdot(h\vel) &= 0 \\
\df{(h\vel)}{t} + \nabla\cdot \left( h \vel \otimes \vel + \frac{gh^2}{2} I + h \p \right) &= -g h \nabla b -C_f |\vel| \vel \\
\df{\p}{t} + \vel\cdot\nabla\p + (\nabla\vel) \p + \p (\nabla\vel)^\top &= \mathcal{D}
\end{align}
\end{subequations}
where the symmetric stress tensor $\p$ comes from the second order velocity fluctuations related to the perturbations from the depth average, and comprises of three independent components $\p_{11}, \p_{12}, \p_{22}$, and $\mathcal{D}$ is the dissipation tensor. Using Stokes type hypothesis, we can relate the dissipation tensor $\mathcal{D}$  linearly with the stress tensor $\p$, and furthermore if we want the model to reduce to the classical shallow water model in the limit $\p = 0$, we get~\cite{Gavrilyuk2018}
\[
\mathcal{D} = -\frac{2 \alpha |\vel|^3}{h} \p
\]
where the coefficient $\alpha$ is a function of the invariants of $\p$. The equations~\eqref{eq:ssw} lead to an equation for ``total energy"
\begin{equation}
\label{eq:tote}
\df{}{t}(he) + \nabla\cdot \left[ h e \vel + (g h^2/2 + h \p)\vel \right] = - C_f |\vel|^3 - Q
\end{equation}
where
\[
e = \half |\vel|^2 + \half \trace(\p) + \half g h^2, \qquad Q = - \half h \trace(\mathcal{D})
\]
Following~\cite{Richard2013,Gavrilyuk2018}, the coefficient $\alpha$ in the dissipation term is given by
\[
\alpha = \max\left(0, C_r \frac{T - \phi h^2}{T^2} \right), \qquad T = \trace(\p) = \p_{11} + \p_{22}
\]
Moreover, the quantities $C_f, C_r, \phi$ are model constants that must be calibrated using experiments. The equations for $\p$ would contain third order velocity fluctuations, see \ref{sec:avg}, which can be ignored if the horizontal shear is weak~\cite{Teshukov2007} or modeled in some way so as to close the set of equations. The depth averaging process which is described in \ref{sec:avg} reminds us of the well known Reynolds averaging of Navier-Stokes equations for turbulent flows which leads to a hierarchy of equations due to the non-linear nature and have to be closed with some turbulence model\footnote{However, unlike the Reynolds average, the depth average does not commute with differentiation.}. A similar model as above has been studied for the Favre-averaged compressible Navier-Stokes equations used to model turblent flows~\cite{Berthon2002}. The SSW model~\eqref{eq:ssw} is non-conservative since the equations for $\p$ cannot be put in conservation form. The solution of non-conservative equations is a tricky issue when discontinuities arise since we need a proper notion of weak solution. The jump conditions will depend on the particular form of the equations and not all forms may yield the same jump conditions.

The numerical solution of the set of equations~\eqref{eq:ssw} has been addressed in~\cite{Gavrilyuk2018,Bhole2019}, both of which use a splitting approach. The approach in~\cite{Gavrilyuk2018} splits the equations into two sub-systems called $a$-waves (acoustic) and $b$-waves (shear), and develops an approximate Riemann solver for each one independently. Each of these sub-systems is also augmented with the energy conservation equation~\eqref{eq:tote} which is used to derive some jump conditions required to develop the Riemann solvers. The approach in~\cite{Bhole2019} also uses the same acoustic and shear sub-systems and develops fluctuation splitting schemes for each sub-system on unstructured grids, but does not make use of the total energy equation~\eqref{eq:tote}.

In the present work, we cast the SSW equations in a particular non-conservative form which  is similar to the 10-moment equations~\cite{Levermore1998,Berthon2006} from gas dynamics. In this model, instead of equations for the stress $\p$, we have equations for an {\em energy tensor} $E$, while the mass and momentum equations remain unchanged. This form of the equations naturally arises when we perform the depth averaging of the 3-D Euler equations and the derivation is given in \ref{sec:avg}. In fact, the equation for the energy tensor appears in~\cite{Teshukov2007} but it has not been used by any of the researchers to develop a numerical approximation. We suggest that the form of the equations is important and hence we retain the equation structure arising from depth averaging to build a numerical approximation. The non-conservative terms in this form contain only derivatives of the water depth $h$ unlike model~\eqref{eq:ssw} which has derivatives of $\vel,\p$ in the non-conservative terms. The presence of only the derivatives of $h$ in the non-conservative terms facilitates the construction of path conservative schemes~\cite{DalMaso1995}. By using the generalized {\em Rankine-Hugoniot (RH) jump conditions} arising from taking a linear path in the state space, we build HLL-type Riemann solvers for the new system. We construct the HLL, a 3-wave HLLC and a 5-wave HLLC solver, with the last one including all the waves in the Riemann problem. Unlike previous works, we do not split the model in several sub-systems but instead we construct a unified Riemann solver for the full system. A higher order version of the scheme is constructed following the MUSCL-Hancock approach~\cite{VanLeer1997} where we make the source terms implicit. The resulting semi-implicit scheme is solved exactly. \reva{While such path conservative schemes provide a framework to construct stable numerical approximations, we should mention that the theoretical analysis of such schemes is not well developed. The knowledge of the correct path may not be known and even when it is known, it is not guaranteed that the numerical scheme will converge to the correct solution, as shown in~\cite{Abgrall2010} in case of Euler equations.}

The rest of the paper is organized as follows. The model form used in the current work is introduced in Section~(\ref{sec:10mom}). The notions of path conservative scheme are discussed in Section~(\ref{sec:1d}) and the first order scheme is presented. The structure of the states across discontinuities using the generalized RH conditions are presented in Section~(\ref{sec:jump}). Sections~(\ref{sec:hll}),~(\ref{sec:hllc3}),~(\ref{sec:hllc5}) present the HLL, 3-wave HLLC and the 5-wave HLLC approximate Riemann solvers, respectively. The higher order versions of the scheme using MUSCL-Hancock-type approach is given in Section~(\ref{sec:1dsecond}) in 1-D and Section~(\ref{sec:2d}) in 2-D. Then we present a set of test cases in Section~(\ref{sec:num}). Section~(\ref{sec:sum}) makes a summary of the work and draws some conclusions. The depth averaging of the shallow water equations that leads to the model form~\eqref{eq:tenmom} is shown in \ref{sec:avg} and the solution of the semi-implicit scheme is presented in \ref{sec:imp}.

\section{Re-formulation of the SSW model}
\label{sec:10mom}
The SSW model~\eqref{eq:ssw} can be written in an almost conservative form. To do this, we define the symmetric tensors
\[
\R_{ij} := h \p_{ij}, \qquad E_{ij} := \half \R_{ij} + \half h v_i v_j, \qquad 1 \le i,j \le 2
\]
The quantity $\R_{ij}$ has units of stress while $E_{ij}$ has units of energy per unit volume. Then an elementary computation shows that the set of equations~\eqref{eq:ssw} can be written as the following set of non-conservative equations
\begin{equation}
\label{eq:tenmom}
\df{\con}{t} + \df{\fx}{x_1} + \df{\fy}{x_2} + \B_1 \df{h}{x_1} + \B_2 \df{h}{x_2} =  \s
\end{equation}
where
\[
\con = \begin{bmatrix}
h \\
h v_1 \\
h v_2 \\
E_{11} \\
E_{12} \\
E_{22} \end{bmatrix}, \quad
\fx = \begin{bmatrix}
h v_1 \\
\R_{11} + h v_1^2 + \half g h^2 \\
\R_{12} + h v_1 v_2 \\
(E_{11} + \R_{11}) v_1 \\
E_{12} v_1 + \half (\R_{11} v_2 + \R_{12} v_1) \\
E_{22} v_1 + \R_{12} v_2 \end{bmatrix}, \quad
\fy = \begin{bmatrix}
h v_2 \\
\R_{12} + h v_1 v_2 \\
\R_{22} + h v_2^2 + \half g h^2 \\
E_{11} v_2 + \R_{12} v_1 \\
E_{12} v_2 + \half (\R_{12} v_2 + \R_{22} v_1) \\
(E_{22} + \R_{22}) v_2 \end{bmatrix}
\]
\[
\B_1 = \begin{bmatrix}
0 \\
0 \\
0 \\
g h v_1 \\
\half g h v_2 \\
0 \end{bmatrix}, \qquad
\B_2 = \begin{bmatrix}
0 \\
0 \\
0 \\
0 \\
\half g h v_1 \\
g h v_2 \end{bmatrix}, \qquad
\s = \begin{bmatrix}
0 \\
-g h \df{b}{x_1} - C_f |\vel| v_1 \\
-g h \df{b}{x_2} - C_f |\vel| v_2 \\
- g h v_1 \df{b}{x_1} - \alpha |\vel|^3 \p_{11} - C_f |\vel| v_1^2 \\
- \half g h v_2 \df{b}{x_1} - \half g h v_1 \df{b}{x_2} - \alpha |\vel|^3 \p_{12}  - C_f |\vel| v_1 v_2 \\
- g h v_2 \df{b}{x_2} -\alpha |\vel|^3 \p_{22}  - C_f |\vel| v_2^2
\end{bmatrix}
\]
\revb{The fluxes can be written in terms of $\con$ using the following transformation
\[
h=U_1, \quad v_1 = U_2/U_1, \quad v_2 = U_3/U_1, \quad \R_{11} = 2 U_4 - U_2^2/U_1, \quad \R_{12} = 2 U_5 - U_2 U_3/U_1, \quad \R_{22} = 2 U_6 - U_3^2/U_1
\]
}
Note that the vectors $\B_1(\con), \B_2(\con)$ are linear in $\con$ and in fact depend only on the momentum density 
\[
\m = h\vel
\]
It is usual to write the non-conservative terms in terms of a matrix-vector product, but since only the derivatives of $h$ appear in the non-conservative terms of this model, it is more convenient to write it in terms of the vectors $\B_1,\B_2$ as above. Coincidentally, this model is identical to the 10-moment model~\cite{Levermore1998,Berthon2006} from gas dynamics except for the presence of gravity terms and the non-conservative terms. While the set of equations~\eqref{eq:tenmom} can be obtained by manipulating equations~~\eqref{eq:ssw}, the more fundamental way to obtain these equations is by the depth averaging process which is performed in \ref{sec:avg}.  \corr{It is also known that the quantity
\[
\eta = \eta(\con) :=  -h \log \left( \frac{\det{\R}}{h^4}\right)
\]
is a convex entropy function~\cite{Berthon2006}. Using the SSW equations and ignoring the source terms, we can derive the entropy equation
\begin{equation}
\df{\eta}{t} + \nabla\cdot(\eta \vel) = 0
\label{eq:enteqn}
\end{equation}
In fact~\cite{Gavrilyuk2018}, shows that the ``total energy" equation~\eqref{eq:tote} and the entropy equation~\eqref{eq:enteqn} are the only two additional conservation laws we can derive from the SSW model, but the ``total energy" is not a convex function of $\con$.} The fact that the entropy $\eta$ is a convex function of $\con$ and this variable set arises naturally when we perform the depth averaging, indicates that the set of variables $\con$ \revb{can be useful for numerical modeling also. While any set of independent variables is fine for smooth solutions, the computation of correct weak solutions will depend on jump conditions which can be different for different set of independent variables.} Hence we propose to construct numerical schemes starting from the SSW equations as given in~\eqref{eq:tenmom}. The numerical strategy we use is based on the concept of path conservative schemes~\cite{Pares2006} applied to system~\eqref{eq:tenmom}, and we recall the basic notions of path conservative schemes in a 1-D version of the above model in the next section. We will assume throughout the paper that the bottom topography $b(x_1,x_2)$ is a continuous function and independent of time, in which case the terms containing this quantity can be treated as source terms and need not be included in the Riemann solver.
\section{Model in 1-D and notion of path conservative scheme}
\label{sec:1d}
We will consider the 1-D SSW model which can be written as
\begin{equation}
\label{eq:ssw1d}
\df{\con}{t} + \df{\f(\con)}{x} + \B(\m) \df{h}{x} = \s(\con)
\end{equation}
where $\f = \fx$, $\B = \B_1$ and the source term is given by
\[
\s = \begin{bmatrix}
0 \\
-g h \df{b}{x} - C_f |\vel| v_1 \\
 - C_f |\vel| v_2 \\
-\alpha |\vel|^3 \p_{11} - g h v_1 \df{b}{x} - C_f |\vel| v_1^2 \\
-\alpha |\vel|^3 \p_{12} - \half g h v_2 \df{b}{x}  - C_f |\vel| v_1 v_2 \\
-\alpha |\vel|^3 \p_{22}  - C_f |\vel| v_2^2
\end{bmatrix}
\]
For simplicity of notation, we will sometimes write the velocity components as $(u,v)= (v_1,v_2)$. The system of equations~\eqref{eq:ssw1d} is a hyperbolic system with eigenvalues~\cite{Gavrilyuk2018,Berthon2002}
\[
\lambda_1 = u - \sqrt{gh + 3 \p_{11}}, \quad \lambda_2 = u - \sqrt{\p_{11}}, \quad \lambda_3 = \lambda_4 = u, \quad \lambda_5 = u+ \sqrt{\p_{11}}, \quad \lambda_6 = u + \sqrt{g h + 3 \p_{11}}
\]
The first and last eigenvalues correspond to genuinely non-linear characteristic fields in the sense of Lax~\cite{Godlewski1996}, while the remaining eigenvalues correspond to linearly degenerate characteristic fields~\cite{Gavrilyuk2018}. Hence $\lambda_1, \lambda_6$ are associated with shock/rarefaction waves while the remaining eigenvalues give rise to shear/contact waves.

Ignoring the source term for the moment as they do not contain derivatives of $\con$, let us write the non-conservative system as
\[
\df{\con}{t} + \A(\con) \df{\con}{x} = 0
\]
If we have discontinuous solutions, then we have to give a meaning to the derivative term which can be done by integration by parts if $\A$ is the gradient of a flux function as in case of conservation laws. If $\A$ is not the gradient of a flux, then the non-conservative product is interpreted as a Borel measure~\cite{DalMaso1995}. This definition requires the choice of a smooth path $\Psi : [0,1] \times \uad \times \uad \to \uad$ connecting the two states $\con_L, \con_R$ across the jump discontinuity at $x=x_0$ such that
\[
\Psi(0; \con_L, \con_R) = \con_L, \qquad \Psi(1; \con_L, \con_R) = \con_R
\]
where $\uad$ is the set of admissible states. Then the non-conservative product is defined as the Borel measure~\cite{DalMaso1995,Gosse2001}
\[
\mu(x_0) = \left[ \int_0^1 \A(\Psi(\xi; \con_L, \con_R)) \dd{\Psi}{\xi} \ud \xi \right] \delta(x_0)
\]
where $\delta$ is the Dirac delta function. This viewpoint is equivalent to the definition of non-conservative product proposed by Volpert~\cite{volpert1967}. Using this notion, a theory of weak solutions can be developed based on which the Riemann problem has usual structure as for conservative systems, leading to shocks or rarefaction waves corresponding to genuinely non-linear characteristic fields and contact waves corresponding to linearly degenerate fields. Across a point of discontinuity moving with speed $S$, a weak solution has to  satisfy the generalized Rankinge-Hugoniot jump condition
\[
\int_0^1 \left[ \A(\Psi(\xi; \con_L, \con_R)) - S I \right] \dd{\Psi}{\xi} \ud \xi = 0
\]
The choice of the correct path is a difficult question and has to be derived from a regularized model motivated from the physical background of the problem. In practice, it is usual to consider the linear path
\[
\Psi(\xi; \con_L, \con_R) = \con_L + \xi (\con_R - \con_L)
\]
Then the jump condition for our model~\eqref{eq:ssw1d} becomes
\begin{equation}
\label{eq:rh}
\int_0^1 \A(\Psi(\xi; \con_L, \con_R))  \dd{\Psi}{\xi} \ud \xi = \f_R - \f_L + \B(\m_L, \m_R) (h_R - h_L) = S(\con_R - \con_L)
\end{equation}
where
\[
\B(\m_L, \m_R) = \B\left( \frac{\m_L + \m_R}{2} \right)
\]
The source term does not make any contribution to the jump conditions since it does not contain derivative of $\con$. The Riemann problem is the building block of a finite volume method and this approach can be used for non-conservative systems also~\cite{Gosse2001,Pares2006}. The main idea is to split the {\em fluctuation} into two parts corresponding to left moving and right moving waves arising in the Riemann solution, where the fluctuation is defined as
\[
\D(\con_L, \con_R) = \int_0^1 \A(\Psi(\xi; \con_L, \con_R))  \dd{\Psi}{\xi} \ud \xi = \D^-(\con_L,\con_R) + \D^+(\con_L, \con_R)
\]
The splitting of the fluctuation can be performed using a Roe-type Riemann solver or HLL-type Riemann solver, the latter being the approach taken in the present work. Assume that there are $m$ waves in the Riemann solution with $m-1$ intermediate states. Let us denote the wave speeds as $S_j$, $j=1,\ldots,m$ and the intermediate states as $\con_j^*$, $j=1,\ldots,m-1$ with $\con_0^* = \con_L$ and $\con_m^* = \con_R$. The fluctuation splitting is given by
\[
\D^\pm(\con_L,\con_R) = \sum_{j=1}^m S_j^\pm (\con^*_{j+1} - \con^*_j)
\]
where
\[
S^- = \min(0, S), \qquad S^+ = \max(0,S)
\]
Let us consider a partition of the domain into disjoint cells of size $\dx$. Let $\con_j^n$ denote the approximation of the cell average value in the $j$'th cell at time $t=t_n$. The first order scheme is given by
\[
\con_j^{n+1} = \con_j^n - \frac{\dt}{\dx} ( \D^{+,n}_\jmh + \D^{-,n}_\jph) + \dt \s(\con_j^{n+\theta}), \qquad
\D_\jph^{\pm,n} = \D^{\pm}(\con_j^n, \con_{j+1}^n)
\]
For $\theta=0$ we obtain an explicit scheme and for $\theta=1$ we obtain a semi-implicit scheme; however the coupling in the semi-implicit scheme is only local to the cell. An exact solution process for the semi-implicit scheme is explained in the \ref{sec:imp}. If the system is conservative, i.e., $\A = \f'(\con)$ for some $\f$, then the above scheme can be written in conservation form with some numerical flux function~\cite{Pares2006}.
\section{Linear waves and jump conditions}
\label{sec:jump}
In this section, we study the structure of the states forming the linearly degenerate waves. Let us define the average and jump operators by
\[
\avg{\cdot} = \frac{ (\cdot)_L + (\cdot)_R}{2}, \qquad \jump{ \cdot } = (\cdot)_R - (\cdot)_L
\]
Then the jump conditions across a discontinuity moving with speed $S$ lead to the following set of equations.
\begin{align*}
\jump{h u} &= S \jump{h} \\
\jump{\R_{11} + h u^2 + \half g h^2} &= S \jump{h u} \\
\jump{\R_{12} + h u v} &= S \jump{h v} \\
\jump{E_{11} u + \R_{11} u} +g \avg{h u} \jump{h} &= S \jump{E_{11}} \\
\jump{E_{12} u + \half (\R_{11} v + \R_{12} u)} + \half g \avg{h v} \jump{h} &= S \jump{E_{12}} \\
\jump{E_{22} u + \R_{12} v} &= S \jump{E_{22}}
\end{align*}
When we construct the HLLC solvers in later sections, we will use the information deduced in the following two  sub-sections to decide the structure of the intermediate states.
\subsection{Contact wave}
\label{sec:jumpc}
Since the contact wave is linearly degenerate, the normal velocity $u$ is continuous across the contact wave, $u_L = u_R = u$, and equal to the speed of the contact wave $S = u$. The jump conditions across the contact wave then yield the following four relations
\begin{equation}
\jump{\R_{11} + \half g h^2} = 0, \qquad \jump{\R_{12}} = 0, \qquad \jump{\R_{11} v} + g \avg{hv} \jump{h} = 0, \qquad \jump{\R_{12} v} = 0
\label{eq:jumpmid}
\end{equation}
From the first and third conditions, we obtain
\[
\left( \avg{\R_{11}} + \frac{1}{4} g \jump{h}^2 \right) \jump{v} = 0
\]
Since we expect $\R_{11}$ to be positive, the first factor cannot be zero and hence we require that $\jump{v} = 0$, so that both velocity components are continuous across the middle wave. The second condition of~\eqref{eq:jumpmid} shows that $\R_{12}$ is also continuous across the middle wave. 
\paragraph{Remark}
In the above derivation, we deduced that $v$ is continuous across the contact wave but the jump conditions admit another solution in some special cases. If the two states are such that
\[
\R_{11}^L = \R_{11}^R = \R_{12}^L = \R_{12}^R =0, \qquad h_L = h_R
\]
then $v_L \ne v_R$ is an admissible set of states that satisfies all the jump conditions. The assumption $\R_{11}^L = \R_{11}^R = 0$ implies that there is no vertical shear which does not hold in practical situations that we are interested in. Hence  it is reasonable to ignore this solution in the construction of the Riemann solver. 
\subsection{Shear waves}
\label{sec:jumps}
Let us now consider the simple waves corresponding to the eigenvalues $\lambda_2, \lambda_5$, \revb{both of which correspond to linearly degenerate eigenvectors.} If the two states $\con_L, \con_R$ correspond to a $\lambda_2$ wave, then they lie on the same integral curve of the eigenvector corresponding to the eigenvalue $\lambda_2 = u - \sqrt{\p_{11}}$. In terms of the variables $(h, u, v, \p_{11}, \p_{12}, \p_{22})$, the corresponding eigenvector is~\cite{Gavrilyuk2018}
\[
\bm{r}_2 = [ 0, \ 0, \ -c, \ 0, \ c^2, \ 2 \p_{12} ]^\top
\]
where $c = \sqrt{\p_{11}}$, and the integral curve is given by
\[
\frac{\ud h}{0} = \frac{\ud u}{0} = \frac{\ud v}{-c} = \frac{\ud \p_{11}}{0} = \frac{\ud \p_{12}}{c^2} = \frac{\ud \p_{22}}{2 \p_{12}}
\]
We immediately see that $h, u, \p_{11}$ are constant along the integral curve and such constants are also called {\em Riemann invariants}. This implies that $\lambda_2$ has the same value in the two states which is consistent with the fact that we have a linearly degenerate field. From the remaining equations, we can deduce that $v \sqrt{\p_{11}} +  \p_{12}$ and $\det \p$ are also invariant along the integral curve. 

Similarly, if we consider a $\lambda_5$ wave, the corresponding eigenvector is 
\[
\bm{r}_5 = [ 0, \ 0, \ c, \ 0, \ c^2, \ 2 \p_{12} ]^\top
\]
and the integral curve is given by
\[
\frac{\ud h}{0} = \frac{\ud u}{0} = \frac{\ud v}{c} = \frac{\ud \p_{11}}{0} = \frac{\ud \p_{12}}{c^2} = \frac{\ud \p_{22}}{2 \p_{12}}
\]
We deduce that the quantities $h, u, \p_{11}, v \sqrt{\p_{11}} - \p_{12}, \det \p$ are invariant along the integral curve.
\paragraph{Remark}
In contrast to the Euler equations, the contact and shear waves are here associated to different eigenvalues when $\R_{11} > 0$. In the limit when $\R_{11}$ goes to zero, the contact and shear waves merge and the system loses its strict hyperbolicity since the eigenvectors $\bm{r}_2, \bm{r}_5$ become parallel. Then the tangential velocity $v$ is free to take any value along the integral curve. In order to capture the interactions between contact and shear wave, as well as the transition when $\R_{11}$ becomes small, accurate numerical scheme should also considered in these waves. Indeed, interactions of contact and shear waves seems to play an important role for the considered model, as is observed in the two dimensional roll wave tests, where only the five wave solver which includes all waves in the Riemann solver is able to produce realistic looking solutions.
\section{HLL Riemann solver}
\label{sec:hll}
The HLL Riemann solver~\cite{Harten1983a,Einfeldt1988} includes only the slowest wave $S_L$ and fastest wave $S_R$ in the Riemann problem. There is an intermediate state $\con_*$ between the two waves. The jump conditions~\eqref{eq:rh} across the two waves are given by
\[
\f_* - \f_L + \B(\m_L, \m_*) (h_* - h_L) = S_L (\con_* - \con_L)
\]
\[
\f_R - \f_* + \B(\m_*, \m_R) (h_R - h_*) = S_R (\con_R - \con_*)
\]
Eliminating $\f_*$ we obtain the intermediate state
\begin{equation}
\con_* = \frac{1}{S_R - S_L} \left[ (\con_R S_R - \con_L S_L) - (\f_R - \f_L) - \B(\m_L, \m_*)(h_* - h_L) - \B(\m_*,\m_R)(h_R - h_*) \right]
\label{eq:int1d}
\end{equation}
This looks like an implicit equation for the intermediate state due to the non-conservative terms, but we can split this equation into two parts; define
\[
\con = \begin{bmatrix}
\cons \\
\conp \end{bmatrix}, \quad \cons = \begin{bmatrix}
h \\
\m \end{bmatrix} = \begin{bmatrix}
h \\
h v_1 \\
h v_2 \end{bmatrix}, \qquad \conp = \begin{bmatrix}
E_{11} \\
E_{12} \\
E_{22} \end{bmatrix}
\]
with similar splitting of $\f,\B$. The first three equations of~\eqref{eq:int1d} do not contain the non-conservative terms and we obtain the corresponding intermediate state
\[
\cons_* = \begin{bmatrix}
h_* \\
\m_* \end{bmatrix} =
\frac{1}{S_R - S_L} \left[ \cons_R S_R - \cons_L S_L - (\fs_R - \fs_L) \right]
\]
The last three equations of~\eqref{eq:int1d} then yield
\[
\conp_* = \frac{1}{S_R - S_L} \left[ \conp_R S_R - \conp_L S_L - (\fp_R - \fp_L) - \Bp(\m_L, \m_*)(h_* - h_L) - \Bp(\m_*,\m_R)(h_R - h_*) \right]
\]
Hence we have explicit solution for the full intermediate state. The split fluctuations are then obtained from 
\[
\D^\pm(\con_L,\con_R)  = S_L^\pm (\con_* - \con_L) + S_R^\pm (\con_R - \con_*)
\]
We estimate the minimum and maximum speeds in the Riemann problem as follows
\[
S_L = \min\{ \lambda_1(\con_L), \lambda_1(\avg{\con}) \}, \qquad S_R = \max\{ \lambda_6(\con_R), \lambda_6(\avg{\con}) \}
\]
which is similar to estimates used for Euler equations, see e.g.,~\cite{Einfeldt1988}.
\section{HLLC3 Riemann solver: 3 waves}
\label{sec:hllc3}
{\renewcommand{\arraystretch}{1.5}
\begin{table}
\begin{center}
\begin{tabular}{|c|c|c|c|c|c|c|}
$h_L$ & & $h_{*L}$ &  & $h_{*R}$  & & $h_R$ \\
$u_L$ & & $u_*$ &  & $u_*$  & & $u_R$ \\
$v_L$ & &  $v_{*}$ &  & $v_{*}$  & & $v_R$ \\
$\R_{11}^L$ & $S_L$ &  $\R_{11}^{*L}$ & $u_*$ & $\R_{11}^{*R}$  & $S_R$ & $\R_{11}^R$ \\
$E_{11}^L$ & &  $E_{11}^{*L}$ & & $E_{11}^{*R}$  & & $E_{11}^R$ \\
$\R_{12}^L$ & &  $\R_{12}^{*}$ &  & $\R_{12}^{*}$  & & $\R_{12}^R$ \\
$E_{22}^L$ & & $E_{22}^{*L}$ & & $E_{22}^{*R}$  & & $E_{22}^R$
\end{tabular}
\end{center}
\caption{Intermediate states for the 3-wave solver. The wave speeds are shown in between the states.}
\label{tab:hllc3}
\end{table}
The HLL solver does not include the linearly degenerate waves like the contact wave, which get excessively diffused in the numerical results. In the HLLC solver~\cite{Toro1994}, the contact wave of speed $S_M$ is also included in the wave model so that we have three waves $S_L < S_M < S_R$ and two intermediate states $\con_{*L}$, $\con_{*R}$. We will determine the two intermediate states by satisfying the jump conditions across the three waves. Since the contact wave is linearly degenerate, the normal velocity in the two intermediate states is same and equal to the speed of the contact wave
\[
u_{*L} = u_{*R} = u_* = S_M
\]
From the discussion on jump conditions in section~(\ref{sec:jumpc}) we know that $v$ and $\R_{12}$ are continuous across this wave so that $v_{*L} = v_{*R} = v_*$ and $\R_{12}^{*L} = \R_{12}^{*R} = \R_{12}^*$. These states are shown in Table~(\ref{tab:hllc3}).

We first consider the jump conditions across the $S_L, S_R$ waves. The jump conditions for the continuity equation yield
\begin{equation}
h_{*\alpha} = h_\alpha \frac{S_\alpha - u_\alpha}{S_\alpha - u_*}, \qquad \alpha = L, R
\label{eq:inth}
\end{equation}
The jump conditions for the $x$-momentum equation yields
\[
\R_{11}^{*\alpha} = \R_{11}^\alpha + (h_\alpha u_\alpha^2 - h_{*\alpha} u_*^2) + \half g (h_\alpha^2 - h_{*\alpha}^2) + S_\alpha (h_{*\alpha} u_* - h_\alpha u_\alpha), \qquad \alpha = L, R
\]
Using~\eqref{eq:inth}, the above equation can be written as
\begin{equation}
\R_{11}^{*\alpha} = \R_{11}^\alpha + h_\alpha (S_\alpha - u_\alpha) (u_* - u_\alpha) + \half g (h_\alpha^2 - h_{*\alpha}^2)
\label{eq:intrxx}
\end{equation}
Substituting~\eqref{eq:intrxx} into the first jump condition in~\eqref{eq:jumpmid}, we can find an expression for the contact wave speed $S_M$
\begin{equation}
S_M = u_* = \frac{ h_R u_R (S_R - u_R) - h_L u_L (S_L - u_L) - (\R_{11}^R - \R_{11}^L) - \half g (h_R^2 - h_L^2)}{h_R (S_R - u_R) - h_L (S_L - u_L)}
\label{eq:ustar}
\end{equation}
The jump conditions from the $y$ momentum equation yield
\begin{equation}
\R_{12}^{*L} = \R_{12}^L + h_L (S_L - u_L) (v_* - v_L), \qquad \R_{12}^{*R} = \R_{12}^R + h_R (S_R - u_R) (v_* - v_L)
\label{eq:intrxy}
\end{equation}
and since $R_{12}$ must be continuous across the contact wave, we obtain from $\R_{12}^{*L} = \R_{12}^{*R}$ an expression for $v_*$ as
\[
v_* = \frac{ h_R v_R (S_R - u_R) - h_L v_L (S_L - u_L) - (\R_{12}^R - \R_{12}^L)}{h_R (S_R - u_R) - h_L (S_L - u_L)}
\]
We can observe that $u_*, v_*$ are equal to the velocity of the intermediate state of the HLL solver. From the last three jump conditions corresponding to the energy equations, we obtain
\begin{align}
\label{eq:exx}
E_{11}^{*\alpha} &= \frac{1}{S_\alpha - u_*} \left[ (S_\alpha - u_\alpha) E_{11}^\alpha + R_{11}^{*\alpha} u_* - R_{11}^\alpha u_\alpha + \half g (h_\alpha u_\alpha + h_{*\alpha} u_*) (h_{*\alpha} - h_\alpha) \right] \\
E_{12}^{*\alpha} &= \frac{1}{S_\alpha - u_*} \left[ (S_\alpha - u_\alpha) E_{12}^\alpha + \half (\R_{11}^{*\alpha} v_* + \R_{12}^* u_*) - \half (\R_{11}^{\alpha} v_\alpha + \R_{12}^\alpha u_\alpha) + \frac{1}{4} g (h_\alpha v_\alpha + h_{*\alpha} v_{*\alpha}) (h_{*\alpha} - h_\alpha) \right] \\
E_{22}^{*\alpha} &= \frac{1}{S_\alpha - u_*} \left[ (S_\alpha - u_\alpha) E_{22}^\alpha + R_{12}^* v_* - R_{12}^\alpha v_\alpha \right]
\end{align}
The intermediate states are finally given by
\[
\con_{*\alpha} = \left[ h_{*\alpha}, \ h_{*\alpha} u_*, \ h_{*\alpha} v_*, \ E_{11}^{*\alpha}, \ E_{12}^{*\alpha}, \ E_{22}^{*\alpha} \right]^\top
\]
and the fluctuations are  given by
\[
\D^\pm(\con_L, \con_R) = S_L^\pm (\con_{*L} - \con_L) + S_M^\pm (\con_{*R} - \con_{*L}) + S_R^\pm (\con_R - \con_{*R})
\]

\section{HLLC5 Riemann solver: 5-waves}
\label{sec:hllc5}
The 3-wave HLLC solver improves upon the HLL solver by including the linearly degenerate contact wave. But the SSW model contains two more linearly degenerate waves and we can try to construct a multi-state HLL solver by including all linearly degenerate waves in the Riemann solution. When we include all the five waves, there are four intermediate states which are shown in Table~(\ref{tab:hllc5}). Note that the quantities $h, u, \R_{11}$ are continuous across the two shear waves $S_{*L}, S_{*R}$ as deduced in Section~(\ref{sec:jumps}).
{\renewcommand{\arraystretch}{1.5}
\begin{table}
\begin{center}
\begin{tabular}{|c|c|c|c|c|c|c|c|c|c|c|}
$h_L$ & & $h_{*L}$ & & $h_{*L}$ &  & $h_{*R}$ & & $h_{*R}$ & & $h_R$ \\
$u_L$ & & $u_*$ & & $u_*$ &  & $u_*$ & & $u_*$ & & $u_R$ \\
$v_L$ & & $v_{*L}$ & & $v_{**}$ &  & $v_{**}$ & & $v_{*R}$ & & $v_R$ \\
$\R_{11}^L$ & $S_L$ & $\R_{11}^{*L}$ & $S_{*L}$ & $\R_{11}^{*L}$ & $u_*$ & $\R_{11}^{*R}$ & $S_{*R}$ & $\R_{11}^{*R}$ & $S_R$ & $\R_{11}^R$ \\
$E_{11}^L$ & & $E_{11}^{*L}$ & & $E_{11}^{*L}$ & & $E_{11}^{*R}$ & & $E_{11}^{*R}$ & & $E_{11}^R$ \\
$\R_{12}^L$ & & $\R_{12}^{*L}$ & & $\R_{12}^{**}$ &  & $\R_{12}^{**}$ & & $\R_{12}^{*R}$ & & $\R_{12}^R$ \\
$E_{22}^L$ & & $E_{22}^{*L}$ & & $E_{22}^{**L}$ & & $E_{22}^{**R}$ & & $E_{22}^{*R}$ & & $E_{22}^R$
\end{tabular}
\end{center}
\caption{Intermediate states for the 5-wave solver. The wave speeds are shown in between the states.}
\label{tab:hllc5}
\end{table}
}
\paragraph{Jump across $S_L$, $S_R$ waves}
The jump condition across $S_L,S_R$ waves is
\[
\f_{*\alpha} - \f_\alpha + \B(\m_\alpha, \m_{*\alpha})(h_{*\alpha} - h_\alpha) = S_\alpha (\con_{*\alpha} - \con_\alpha), \qquad \alpha = L, R
\]
which can be written as
\[
\underbrace{
\begin{bmatrix}
h_{*\alpha} u_* \\
\R_{11}^{*\alpha} + h_{*\alpha} u_*^2 + \half g h_{*\alpha}^2 \\
\R_{12}^{*\alpha} + h_{*\alpha} u_* v_{*\alpha} \\
E_{11}^{*\alpha} u_* + \R_{11}^{*\alpha} u_* \\
E_{12}^{*\alpha} u_* + \half(\R_{11}^{*\alpha} v_{*\alpha} + \R_{12}^{*\alpha} u_*) \\
E_{22}^{*\alpha} u_* + \R_{12}^{*\alpha} v_{*\alpha}
\end{bmatrix}}_{\f_{*\alpha}} - S_\alpha \underbrace{\begin{bmatrix}
h_{*\alpha} \\
h_{*\alpha} u_* \\
h_{*\alpha} v_{*\alpha} \\
E_{11}^{*\alpha} \\
E_{12}^{*\alpha} \\
E_{22}^{*\alpha} \end{bmatrix}}_{\con_{*\alpha}} +  \begin{bmatrix}
0 \\
0 \\
0 \\
\half g (h_\alpha u_\alpha + h_{*\alpha} u_*) \\
\frac{1}{4}g (h_\alpha v_\alpha + h_{*\alpha} v_{*\alpha}) \\
0
\end{bmatrix} (h_{*\alpha} - h_\alpha)
= \f_\alpha - S_\alpha \con_\alpha
\]
The first jump condition gives $h_{*\alpha}$ which is identical to equation~\eqref{eq:inth}. The second jump condition gives $\R_{11}^{*\alpha}$ and it is identical to equation~\eqref{eq:intrxx}. The third and fifth conditions are coupled and their solution yields $v_{*\alpha}$, $\R_{12}^{*\alpha}$. Define $m_\alpha = h_\alpha (u_\alpha - S_\alpha)$ and $p = \R_{11} + \half g h^2$. There is a common value of $p$ in all the four intermediate states which is given by
\[
p_* = \R_{11}^{*L} + \half g h_{*L}^2 = \R_{11}^{*R} + \half g h_{*R}^2 = \frac{m_R p_L - m_L p_R - m_L m_R (u_R - u_L)}{m_R - m_L}
\]
where the last equality can be derived from the $h$ and $x$-momentum jump conditions across the $S_{L}, S_{R}$ waves. Then the third and fifth jump conditions can be written as
\begin{align*}
m_\alpha (v_{*\alpha} - v_\alpha) + h_{*\alpha} \p_{12}^{*\alpha} &= h_\alpha \p_{12}^\alpha \\
\left[ p_* - \half g h_\alpha h_{*\alpha} \right] (v_{*\alpha} - v_\alpha) + m_\alpha \p_{12}^{*\alpha} &= [m_\alpha + h_\alpha (u_\alpha - u_*)] \p_{12}^\alpha
\end{align*}
whose solution is given by
\[
v_{*\alpha} = v_\alpha + \left[ \frac{m_\alpha (h_\alpha - h_{*\alpha}) - h_\alpha h_{*\alpha}(u_\alpha - u_*)}{m_\alpha^2 - h_{*\alpha} p_* + \half g h_\alpha h_{*\alpha}^2} \right] \p_{12}^\alpha, \qquad \p_{12}^{*\alpha} = \left[ \frac{m_\alpha^2 - h_\alpha p_* + \half g h_\alpha^2 h_{*\alpha} + m_\alpha h_\alpha (u_\alpha - u_*) }{ m_\alpha^2 - h_{*\alpha} p_* + \half g h_\alpha h_{*\alpha}^2} \right] \p_{12}^\alpha
\]
From the above solution, we can compute $\R_{12}^{*\alpha} = h_{*\alpha} \p_{12}^{*\alpha}$ and $E_{12}^{*\alpha} = \half \R_{12}^{*\alpha} + \half h_{*\alpha} u_* v_{*\alpha}$. The fourth jump condition gives $E_{11}^{*\alpha}$ which is identical to equation~\eqref{eq:exx}. The sixth jump condition gives $E_{22}^{*\alpha}$
\[
E_{22}^{*\alpha} = \frac{1}{S_\alpha - u_*} \left[ (S_\alpha - u_\alpha) E_{22}^\alpha + R_{12}^{*\alpha} v_{*\alpha} - R_{12}^\alpha v_\alpha \right]
\]
The second and fourth wave speeds can now be estimated as
\[
\p_{11}^{*\alpha} = \frac{\R_{11}^{*\alpha}}{h_{*\alpha}}, \qquad S_{*L} = u_* - \sqrt{ \p_{11}^{*L}}, \qquad S_{*R} = u_* + \sqrt{ \p_{11}^{*R}}
\]
\paragraph{Jump across $u_*$ wave}
All the jump conditions are satisfied provided we choose $u_*$ as given in equation~\eqref{eq:ustar}.
\paragraph{Jump across $S_{*L}$ wave}
The non-conservative terms are absent for this wave since $h$ is continuous. The first, second and fourth jump conditions are automatically satisfied. The third jump condition yields
\[
(\R_{12}^{**} + h_{*L} u_* v_{**}) - (\R_{12}^{*L} + h_{*L} u_* v_{*L}) = (u_* - \sqrt{\p_{11}^{*L}})(h_{*L} v_{**} - h_{*L} v_{*L})
\]
\[
\R_{12}^{**} = \R_{12}^{*L} - h_{*L} \sqrt{\p_{11}^{*L}} (v_{**} - v_{*L})
\]
The sixth jump condition yields
\[
E_{22}^{**L} = E_{22}^{*L} - \frac{1}{\sqrt{\p_{11}^{*L}}} ( \R_{12}^{**} v_{**} - \R_{12}^{*L} v_{*L})
\]
We can now check that the fifth jump condition is also satisfied.
\paragraph{Jump across $S_{*R}$ wave}
The non-conservative terms are absent for this wave since $h$ is continuous. The first, second and fourth jump conditions are automatically satisfied. The third jump condition yields
\[
(\R_{12}^{**} + h_{*R} u_* v_{**}) - (\R_{12}^{*R} + h_{*R} u_* v_{*R}) = (u_* + \sqrt{\p_{11}^{*R}})(h_{*R} v_{**} - h_{*R} v_{*R})
\]
\[
\R_{12}^{**} = \R_{12}^{*R} + h_{*R} \sqrt{\p_{11}^{*R}} (v_{**} - v_{*R})
\]
The sixth jump condition yields
\[
E_{22}^{**R} = E_{22}^{*R} + \frac{1}{\sqrt{\p_{11}^{*R}}} ( \R_{12}^{**} v_{**} - \R_{12}^{*R} v_{*R})
\]
We can now check that the fifth jump condition is also satisfied.

All the intermediate quantities have been determined except $v_{**}$; we have two estimates of $\R_{12}^{**}$ from the jumps across $S_{*L}$, $S_{*R}$ waves. Setting these equal to one another, we get an equation for $v_{**}$
\[
v_{**} = \frac{ h_{*L} v_{*L} \sqrt{\p_{11}^{*L}} + h_{*R} v_{*R} \sqrt{\p_{11}^{*R}} - (\R_{12}^{*R} - \R_{12}^{*L})}{ h_{*L} \sqrt{\p_{11}^{*L}} + h_{*R} \sqrt{\p_{11}^{*R}} }
\]
Finally, the fluctuations can be computed as
\[
\D^\pm(\con_L,\con_R) = S_L^\pm (\con_{*L} - \con_L) + S_{*L}^\pm (\con_{**L} - \con_{*L}) + u_*^\pm (\con_{**R} - \con_{**L}) + S_{*R}^\pm (\con_{*R} - \con_{**R}) + S_R^\pm (\con_R - \con_{*R})
\]
\section{Second order scheme in 1-D}
\label{sec:1dsecond}
We will follow a predictor-corrector approach like MUSCL-Hancock scheme which can also be thought of as an ADER scheme. We have to reconstruct the solution variables in each cell in order to get a better representation of the solution, and it is found to be beneficial in terms of getting non-oscillatory solutions to reconstruct primitive variables, which are taken to be
\[
\prim = [h, \ v_1, \ v_2, \ \R_{11}, \ \R_{12}, \ \R_{22} ]^\top
\]
\paragraph{Step 1} The first step is to predict the solution at half time level using a local formulation, i.e., without any information from neighbouring cells. We estimate the spatial variation of the solution by a limiter function, e.g.,
\[
\Delta \prim_j^n = \minmod\left( \beta(\prim_j^n - \prim_{j-1}^n), \half(\prim_{j+1}^n - \prim_{j-1}^n), \beta( \prim_{j+1}^n - \prim_j^n) \right), \qquad \beta \in [1,2]
\]
and transform to conserved variables
\[
\Delta\con_j^n = \df{\con}{\prim}(\prim_j^n) \Delta\prim_j^n
\]
Now we can estimate the solution at the faces by linear reconstruction,
\[
\con_{\jmh,R}^{n} = \con_j^n - \half \Delta \con_j^n, \qquad \con_{\jph,L}^{n} = \con_j^n + \half \Delta \con_j^n
\]
Then the solution in the cell is evolved by half a time step using the PDE
\[
\con_j^\nph = \con_j^n + \frac{\dt}{2} \partial_t \con_j^n
\]
where
\[
\partial_t \con_j^n = - \frac{ \f(\con_{\jph,L}^{n}) - \f(\con_{\jmh,R}^{n})}{\dx} - \B(\con_j^n) \frac{\Delta h_j^n}{\dx} +  \s(\con_j^{n+\theta})
\]
The update equation is explicit if $\theta=0$ and is an implicit equation for $\con_j^\nph$ if $\theta=\half$; once we obtain $\con_j^\nph$, then we estimate the solution at the cell faces at the half time level as
\[
\con_{\jmh,R}^{\nph} = \con_j^n - \half \Delta \con_j^n + \frac{\dt}{2} \partial_t \con_j^n, \qquad
\con_{\jph,L}^{\nph} = \con_j^n + \half \Delta \con_j^n + \frac{\dt}{2} \partial_t \con_j^n
\]
\paragraph{Step 2} The second step uses the predicted solution at half time level to update the solution to next time level
\[
\con_j^{n+1} = \con_j^n + \dt \left[ -\frac{1}{\dx}( \tD_\jmh^+ + \tD_\jph^- ) - \frac{ \f(\con_{\jph,L}^{\nph}) - \f(\con_{\jmh,R}^{\nph}) }{\dx} - \B(\con_j^\nph) \frac{\Delta h_j^n}{\dx} + \s(\con_j^\nph) \right]
\]
where
\[
\tD_\jmh^+ = \D^+(\con_{\jmh,L}^{\nph}, \con_{\jmh,R}^{\nph}), \qquad \tD_\jph^- = \D^-(\con_{\jph,L}^{\nph}, \con_{\jph,R}^{\nph})
\]
Note that the coupling between the cells occurs in this stage via the computation of fluctuations that involve the solution from the neighbouring cells.
\section{Two dimensional scheme}
\label{sec:2d}
We consider a Cartesian mesh with cells of size $(\Delta x, \Delta y)$ where the cell centers are indexed as $(j,k)$ and we use half indices $(\jph,k)$ and $(j,\kph)$ to denote the cell faces. The MUSCL-Hancock type scheme can now be naturally extended to two dimensions by the following two step process.
\paragraph{Step 1}
We estimate the the spatial variation of the primitive variables along the two directions
\[
\Delta_x \prim_{j,k}^n = \minmod\left( \beta(\prim_{j,k}^n - \prim_{j-1,k}^n), \half(\prim_{j+1,k}^n - \prim_{j-1,k}^n), \beta( \prim_{j+1,k}^n - \prim_{j,k}^n) \right), \qquad \beta \in [1,2]
\]
\[
\Delta_y \prim_{j,k}^n = \minmod\left( \beta(\prim_{j,k}^n - \prim_{j,k-1}^n), \half(\prim_{j,k+1}^n - \prim_{j,k-1}^n), \beta( \prim_{j,k+1}^n - \prim_{j,k}^n) \right), \qquad \beta \in [1,2]
\]
and transform to conserved variables
\[
\Delta_x \con_{j,k}^n = \df{\con}{\prim}(\prim_{j,k}^n) \Delta_x \prim_{j,k}^n, \qquad \Delta_y \con_{j,k}^n = \df{\con}{\prim}(\prim_{j,k}^n) \Delta_y \prim_{j,k}^n
\]
Now we can estimate the solution at the face centers by linear reconstruction,
\[
\con_{\jmh,k,R}^{n} = \con_{j,k}^n - \half \Delta_x \con_{j,k}^n, \qquad \con_{\jph,k,L}^{n} = \con_{j,k}^n + \half \Delta_x \con_{j,k}^n
\]
\[
\con_{j,\kmh,R}^{n} = \con_{j,k}^n - \half \Delta_y \con_{j,k}^n, \qquad \con_{j,\kph,L}^{n} = \con_{j,k}^n + \half \Delta_y \con_{j,k}^n
\]
Then the solution in the cell is evolved by half a time step using the PDE~\eqref{eq:tenmom}
\[
\con_{j,k}^\nph = \con_{j,k}^n + \frac{\dt}{2} \partial_t \con_{j,k}^n
\]
where
\begin{align*}
\partial_t \con_{j,k}^n = & - \frac{ \fx(\con_{\jph,k,L}^{n}) - \fx(\con_{\jmh,k,R}^{n})}{\dx} - \frac{ \fy(\con_{j,\kph,L}^{n}) - \fy(\con_{j,\kmh,R}^{n})}{\dy} \\
& - \B_1(\con_{j,k}^n) \frac{\Delta h_{j,k}^n}{\dx} - \B_2(\con_{j,k}^n) \frac{\Delta h_{j,k}^n}{\dy} +  \s(\con_{j,k}^{n+\theta})
\end{align*}
The above update is explicit if $\theta=0$ and is an implicit equation for $\con_{j,k}^\nph$ if $\theta=\half$; once we obtain $\con_{j,k}^\nph$, then we estimate the solution at the cell faces by a linear extrapolation from the cell centers
\[
\con_{\jmh,k,R}^{\nph} = \con_{j,k}^\nph - \half \Delta_x \con_{j,k}^n, \qquad
\con_{\jph,k,L}^{\nph} = \con_{j,k}^\nph + \half \Delta_x \con_{j,k}^n
\]
\[
\con_{j,\kmh,k,R}^{\nph} = \con_{j,k}^\nph - \half \Delta_y \con_{j,k}^n, \qquad
\con_{j,\kph,L}^{\nph} = \con_{j,k}^\nph + \half \Delta_y \con_{j,k}^n
\]
\paragraph{Step 2}
The second step uses the predicted solution at half time level to update the solution to next time level
\begin{align*}
\con_{j,k}^{n+1} = \con_{j,k}^n + \dt \Bigg[ & - \frac{ \fx(\con_{\jph,k,L}^{\nph}) + \tD^-_{\jph,k} - \fx(\con_{\jmh,R}^{\nph}) + \tD^+_{\jmh,k}}{\dx} \\
& - \frac{ \fy(\con_{j,\kph,k,L}^{\nph}) + \tD^-_{j,\kph} - \fy(\con_{j,\kmh,R}^{\nph}) + \tD^+_{j,\kmh}}{\dy} \\
& - \B_1(\con_{j,k}^\nph) \frac{\Delta_x h_{j,k}^n}{\dx} - \B_2(\con_{j,k}^\nph) \frac{\Delta_y h_{j,k}^n}{\dy} + \s(\con_{j,k}^\nph) \Bigg]
\end{align*}
This completes the description of the higher order scheme in 2-D. The first order scheme in 2-D is easily obtained by eliminating the reconstruction step and the predictor step.
\section{Numerical results}
\label{sec:num}
The methods developed in this paper are applied to some 1-D and 2-D test cases. The 1-D tests are performed using a purely 1-D code. In all the test cases, we take the acceleration due to gravity as $g = 9.81 \ m/s^2$. The time step is chosen from a CFL condition of the form
\[
\Delta t = \textrm{CFL} \frac{1}{\max_{j,k} \left[ \frac{\lambda_x(\con_{j,k})}{\dx} + \frac{\lambda_y(\con_{j,k})}{\dy} \right]}
\]
where $\lambda_x, \lambda_y$ are the maximum wave speeds along the $x,y$ directrions respectively, and we use CFL = 0.5 in all test cases. During each computation, we monitor positivity of the values of $h, \p_{11}, \p_{22}$ and we do not perform any artificial clipping of negative values. \reva{In the numerical tests, we compare the results from the present method with those from~\cite{Bhole2019} on a fine mesh. Both methods are based on Riemann solver idea but the computational cost can be different. The current method is found to be about 8\% more expensive compared to the method in~\cite{Bhole2019}, since the new Riemann solver has more arithmetic operations and due to the need to transform to conserved variables to calculate the fluctuations.}
\subsection{1-D shear test problem}
\label{sec:shear}
This is a Riemann problem without source terms which gives rise to two shear waves.  The initial depth, normal velocity and the stress tensor are constant in space and are given by 
\[
h = 0.01 \ m, \qquad v_1 = 0, \qquad \p_{11}  = \p_{22} = 10^{-4} \ m^2/s^2, \qquad \p_{12}=0 \ m^2/s^2
\]
Only the transverse velocity has an initial discontinuity located at the middle of the domain and is given by
\[
v_2 = \begin{cases}
\ \ 0.2 \ m/s, & x < 0.5 \ m \\
-0.2 \ m/s, & x > 0.5 \ m
\end{cases}
\]
The computations are performed at first and second order using different Riemann solvers and mesh resolutions on the domain $[0,1]$. The first order results on a coarse mesh (500 cells) and a fine mesh (2000 cells) are shown in Figure~(\ref{fig:shearO1b}). The second order results are shown in Figure~(\ref{fig:shearO2b}) where a coarse mesh of 200 cells and a fine mesh of 2000 cells are used. The reference results are obtained using the method in~\cite{Bhole2019} on a grid of 10000 cells. The solution consists of two shear waves with discontinuities only in the transverse velocity, the $\p_{12}$ and the $\p_{22}$ components of the
stress tensor. It is not surprising that the numerical
error is larger for solvers that do not have these waves in their
intrinsic structure (HLL and HLLC3). However, the differences are
reduced when second order accurate schemes are used and globally the
mesh convergence is achieved for all schemes and are in accordance
with the solution obtained in~\cite{Gavrilyuk2018},~\cite{Bhole2019}. \reva{The plots of $\p_{22}$ in Figures~(\ref{fig:shearO1b}),~(\ref{fig:shearO2b}) show spurious spike in the middle of the computational domain, which have also been observed with other methods in the literature. The HLLC5 solver gives the solution with smallest amplitude of these spikes.}

\begin{figure}
	\begin{center}
		\begin{tabular}{ccc}
			\includegraphics[width=0.40\textwidth]{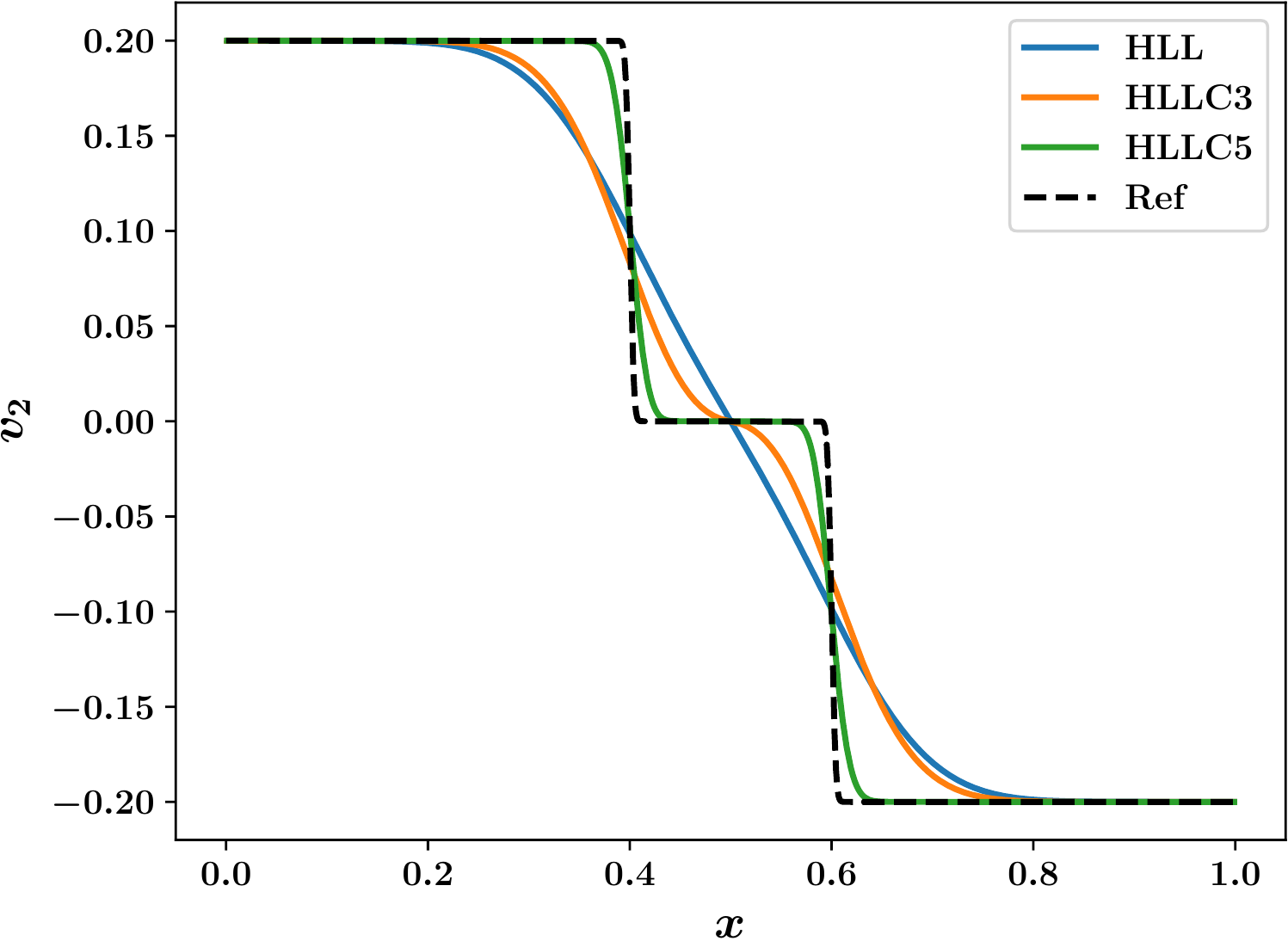} &
         \includegraphics[width=0.40\textwidth]{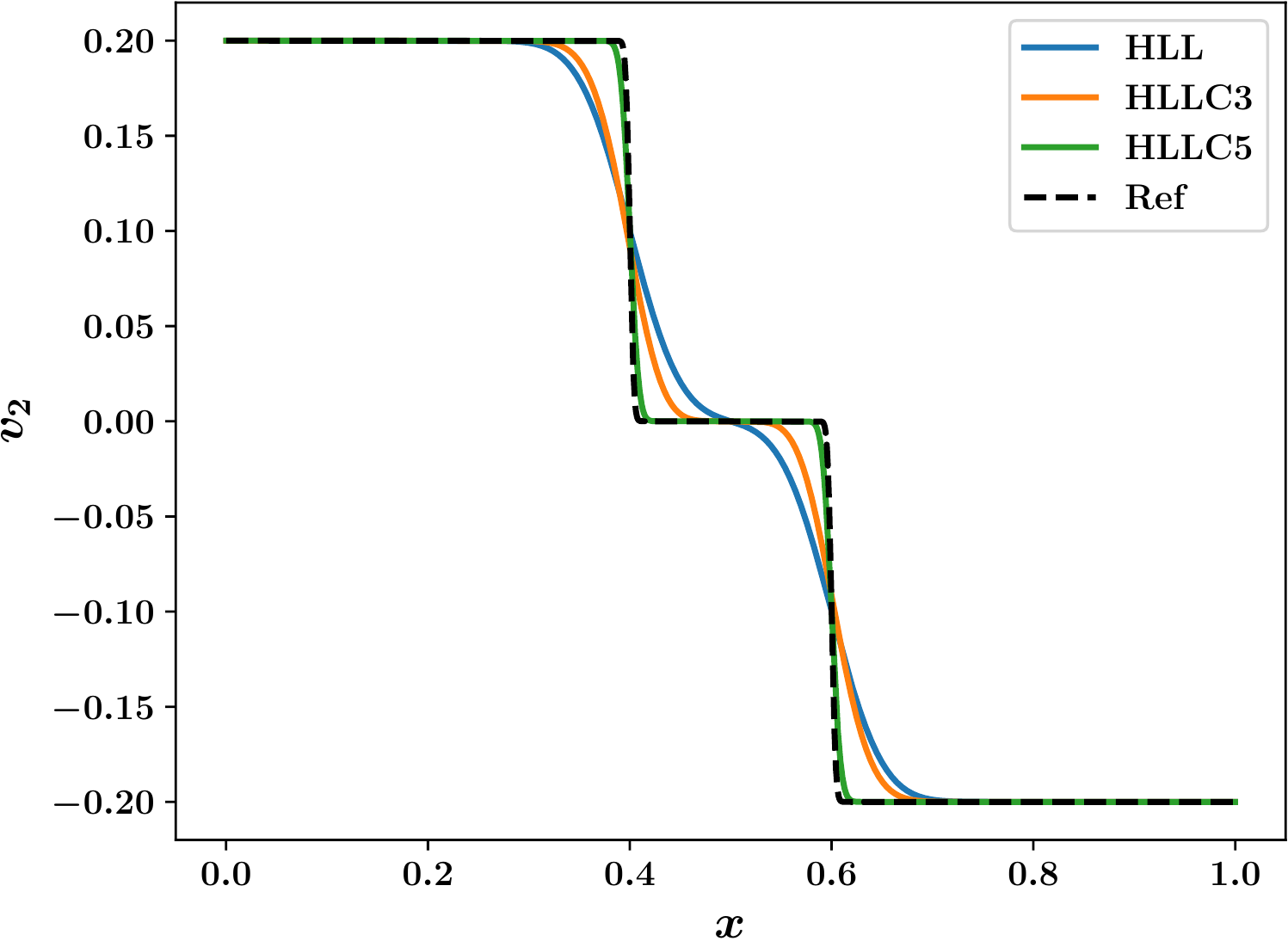} \\
			(a) $v_{2}$, 500 cells & (b) $v_{2}$, 2000 cells \\
			\includegraphics[width=0.40\textwidth]{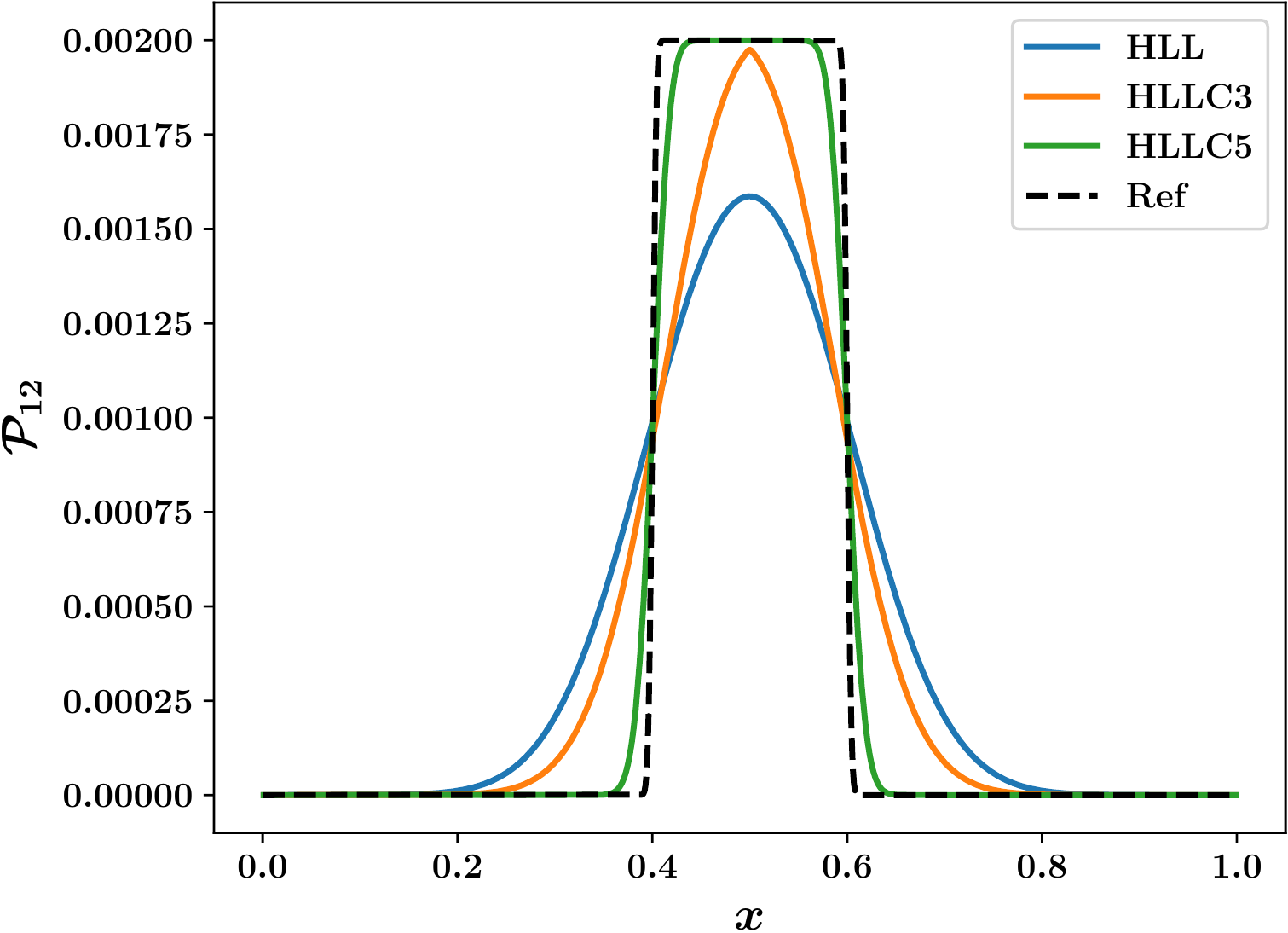} &
         \includegraphics[width=0.40\textwidth]{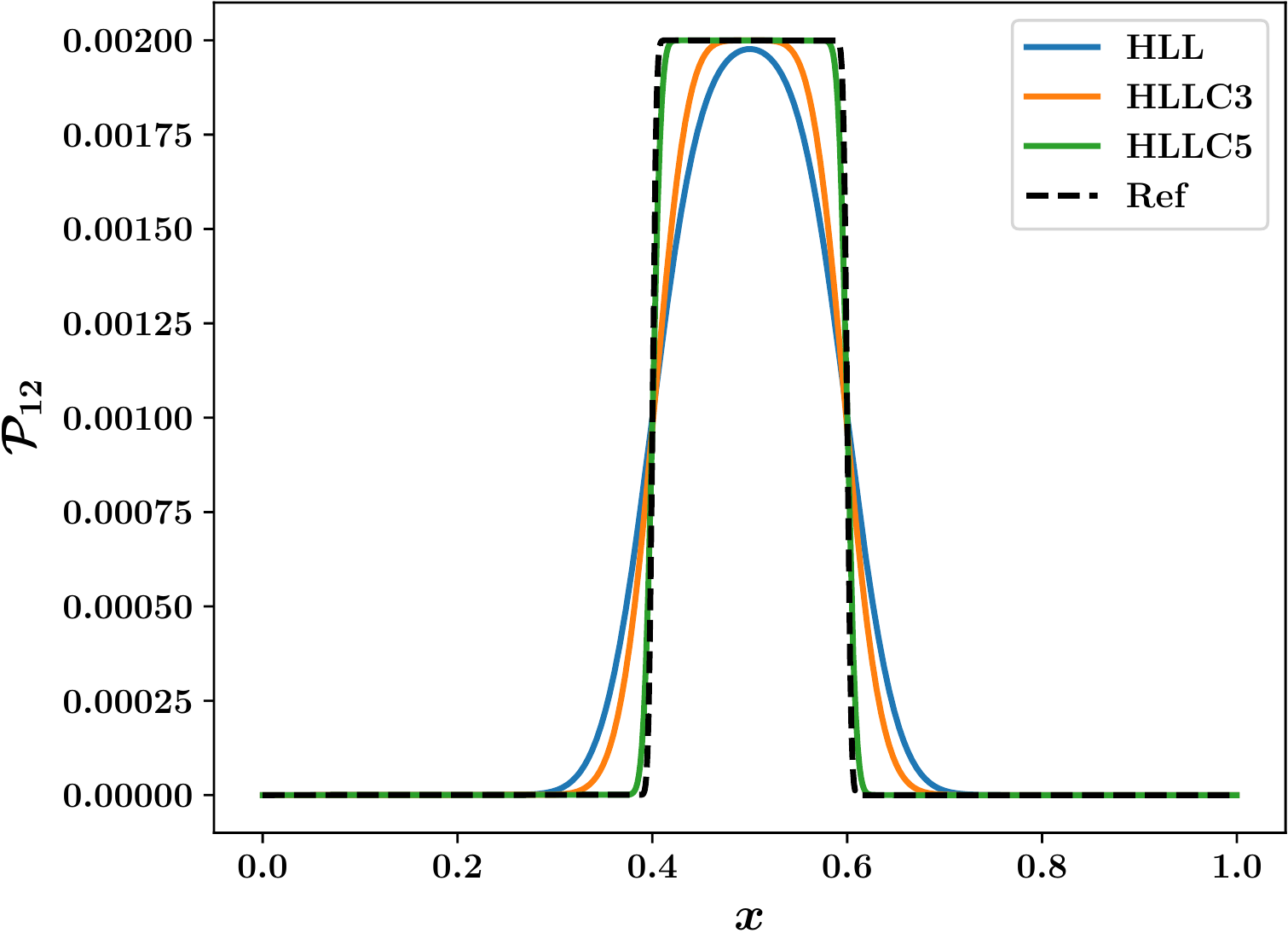} \\
			(a) $\p_{12}$, 500 cells & (b) $\p_{12}$, 2000 cells \\
			\includegraphics[width=0.40\textwidth]{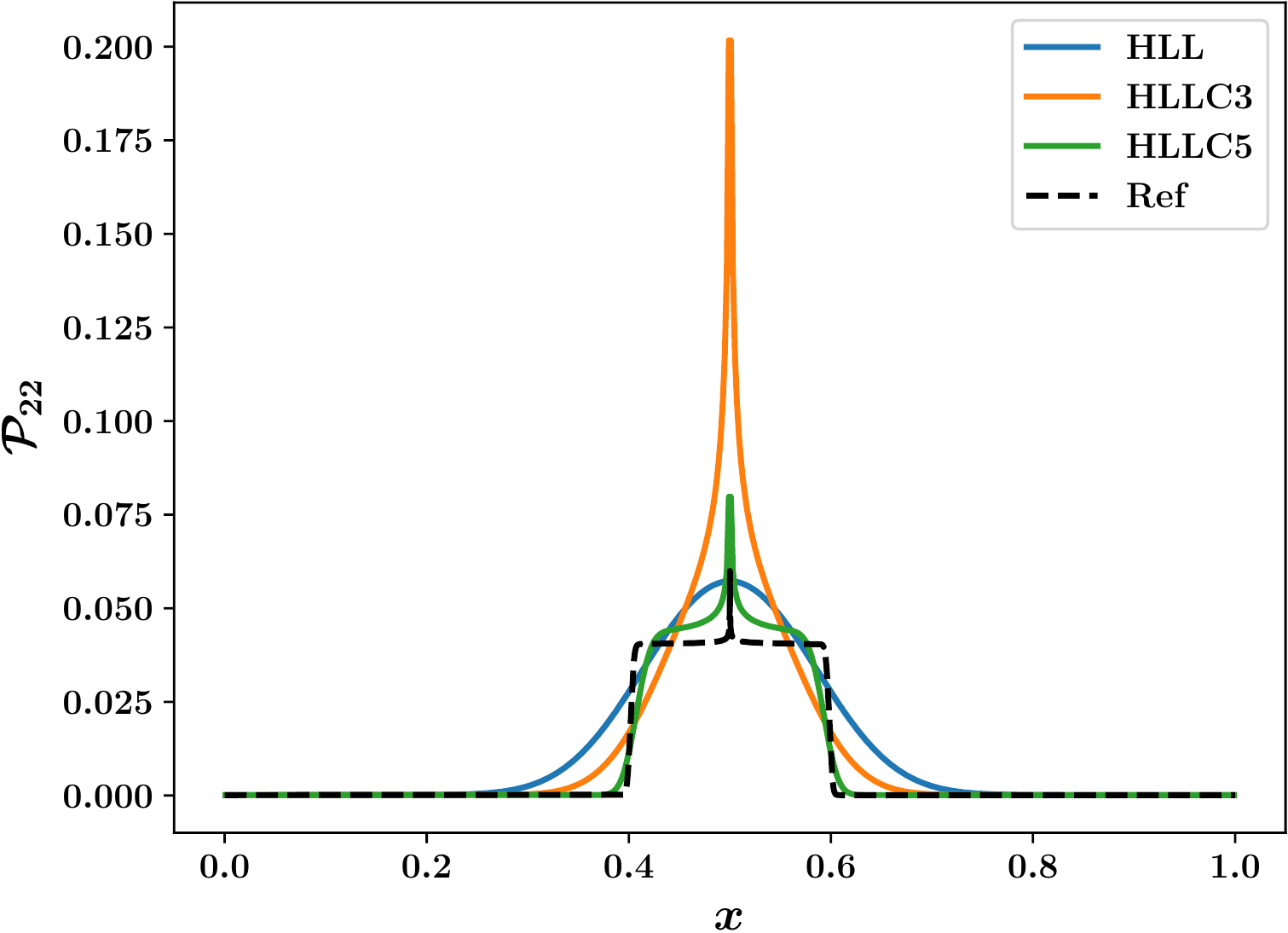} &
         \includegraphics[width=0.40\textwidth]{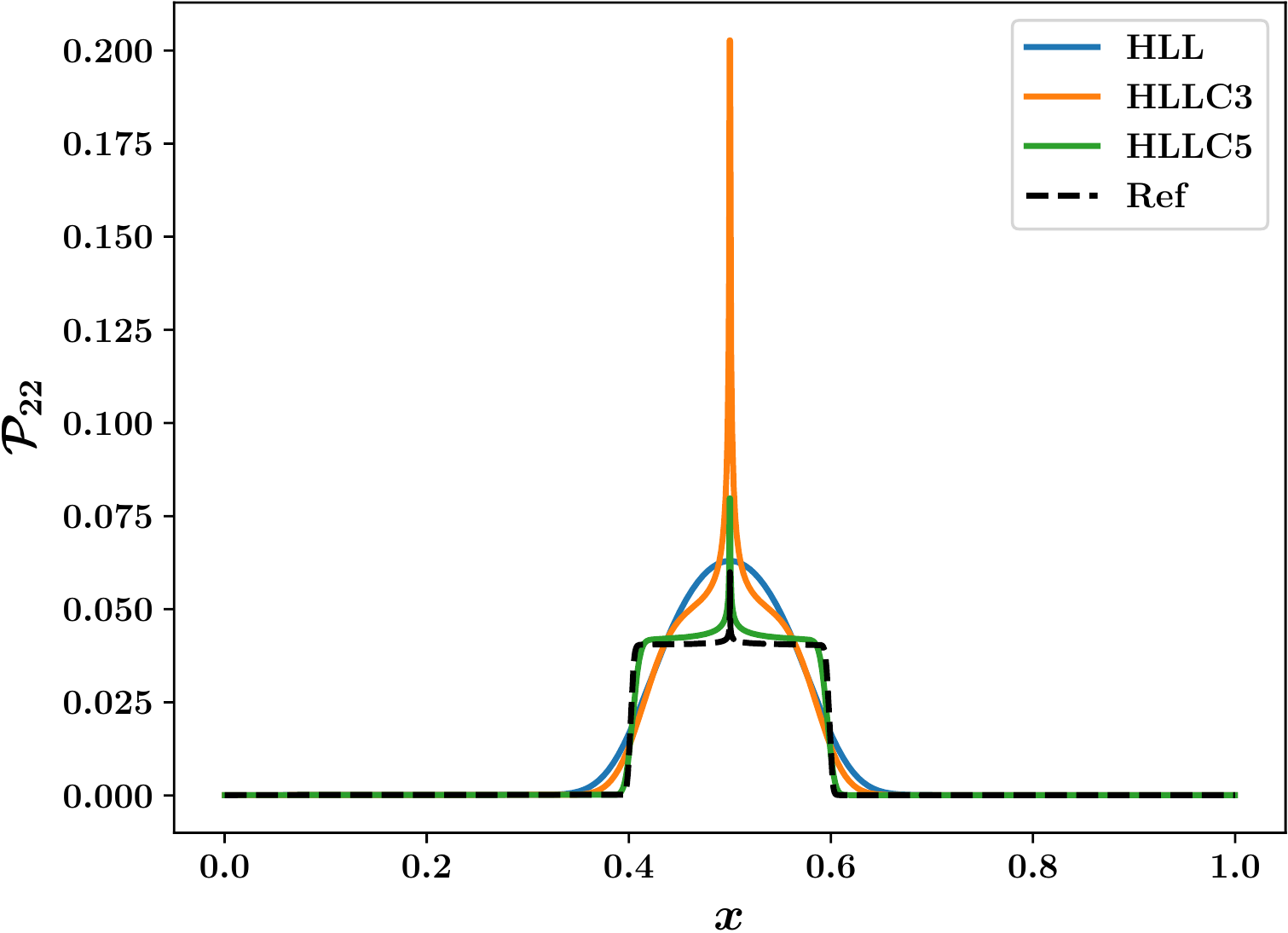} \\
			(a) $\p_{22}$, 500 cells & (b) $\p_{22}$, 2000 cells
		\end{tabular}
	\end{center}
	\caption{1-D shear test case from Section~(\ref{sec:shear}). Plots of y-velocity, stress tensor component $\p_{12}$ and $\p_{22}$ obtained using first order scheme with HLL, HLLC3 and HLLC5 solvers for 500 and 2000 cells  compared with reference solutions from~\cite{Bhole2019}.}
	\label{fig:shearO1b}
\end{figure}

\begin{figure}
	\begin{center}
		\begin{tabular}{ccc}
			\includegraphics[width=0.40\textwidth]{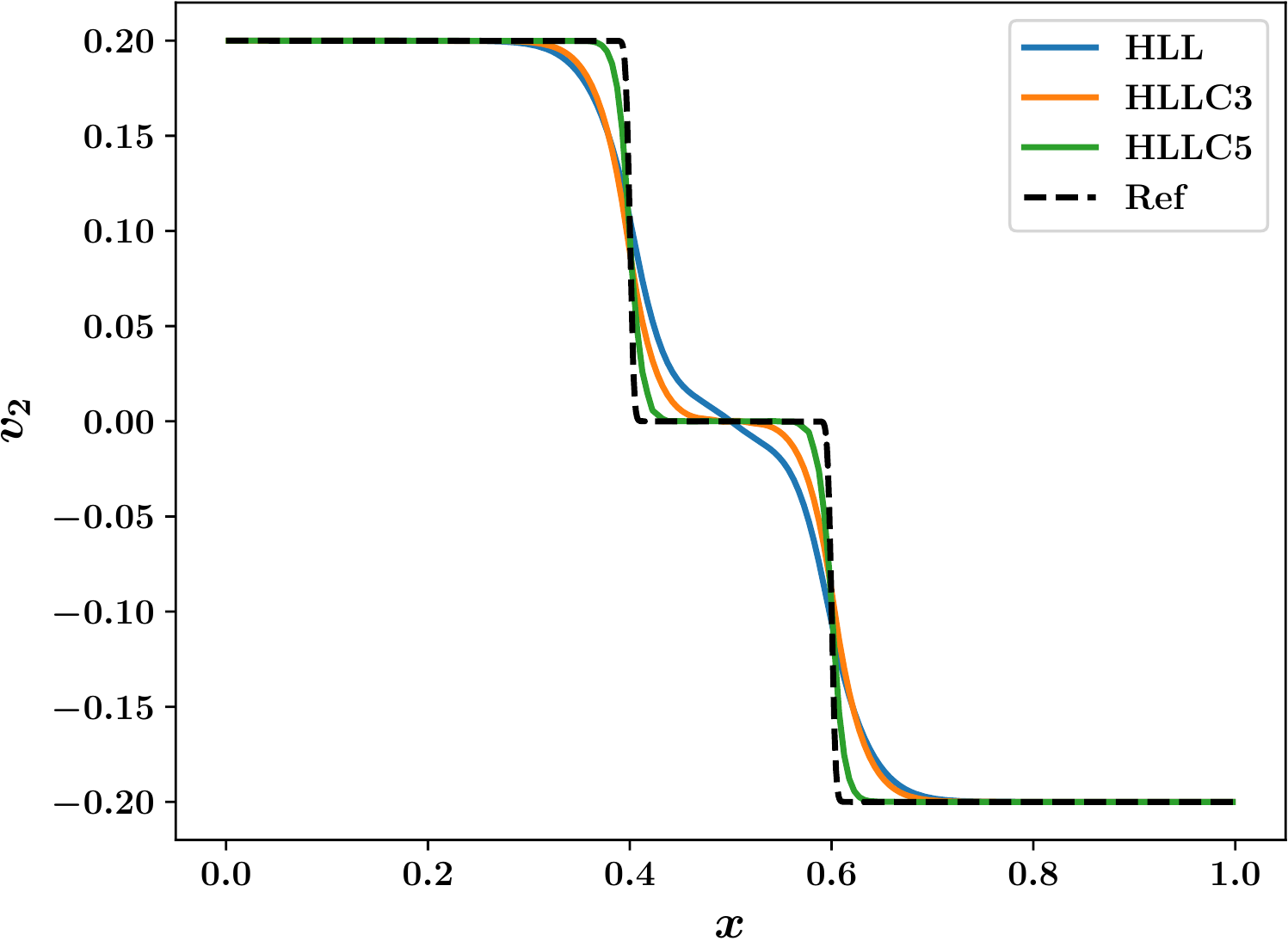} &
         \includegraphics[width=0.40\textwidth]{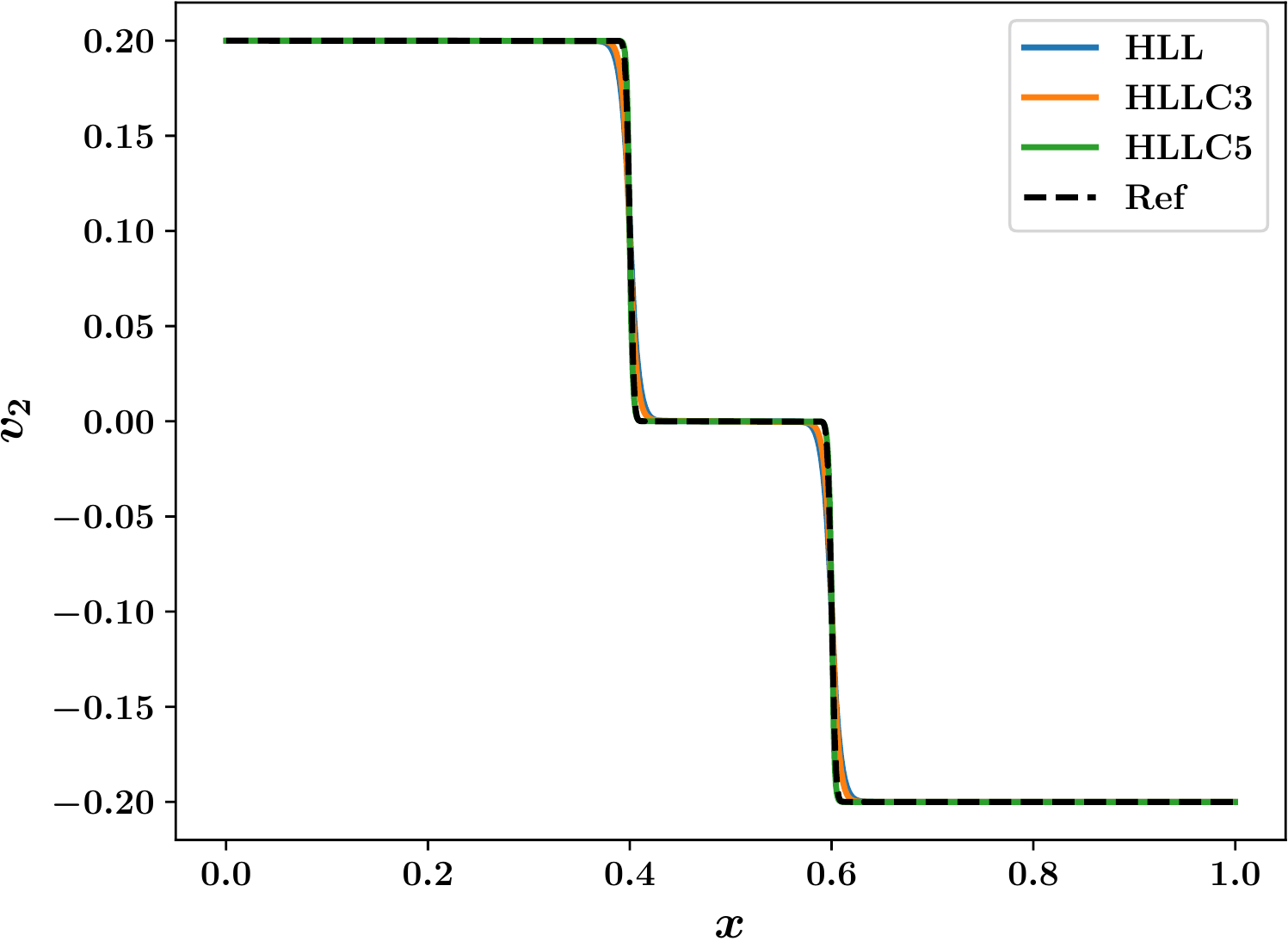} \\
			(a) $v_{2}$, 200 cells & (b) $v_{2}$, 2000 cells \\
			\includegraphics[width=0.40\textwidth]{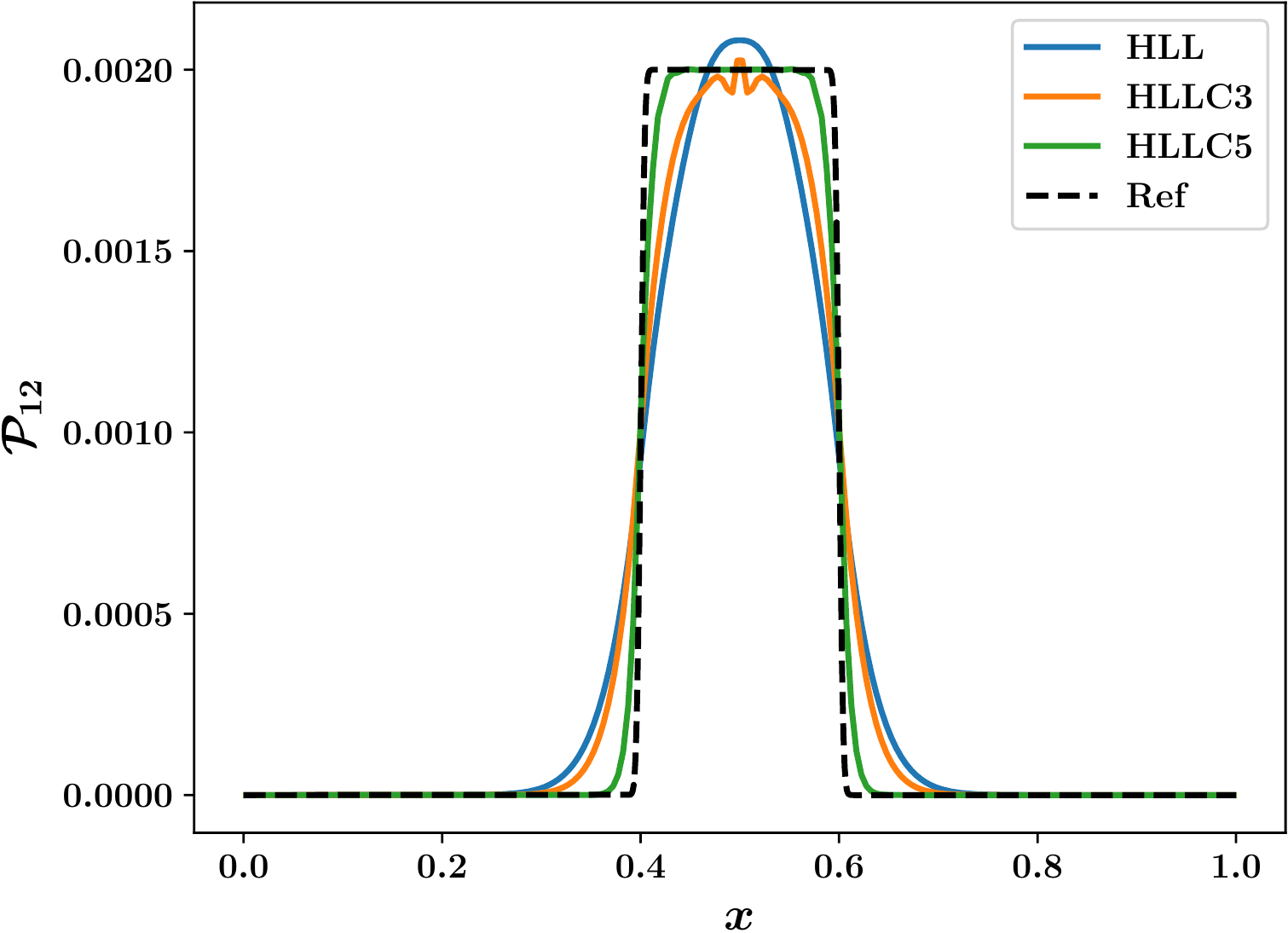} &
         \includegraphics[width=0.40\textwidth]{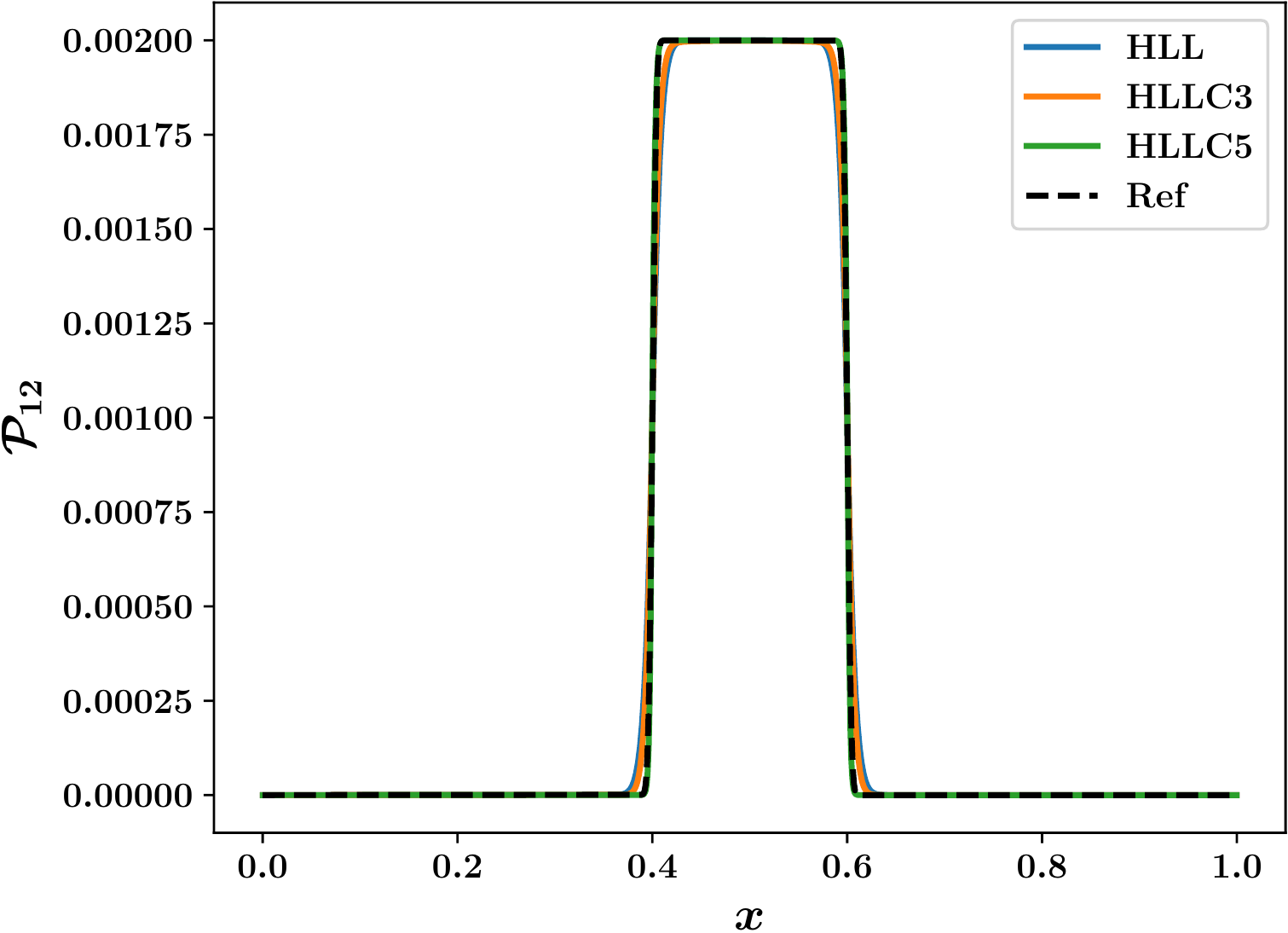} \\
			(a) $\p_{12}$, 200 cells & $\p_{12}$, 2000 cells \\
			\includegraphics[width=0.40\textwidth]{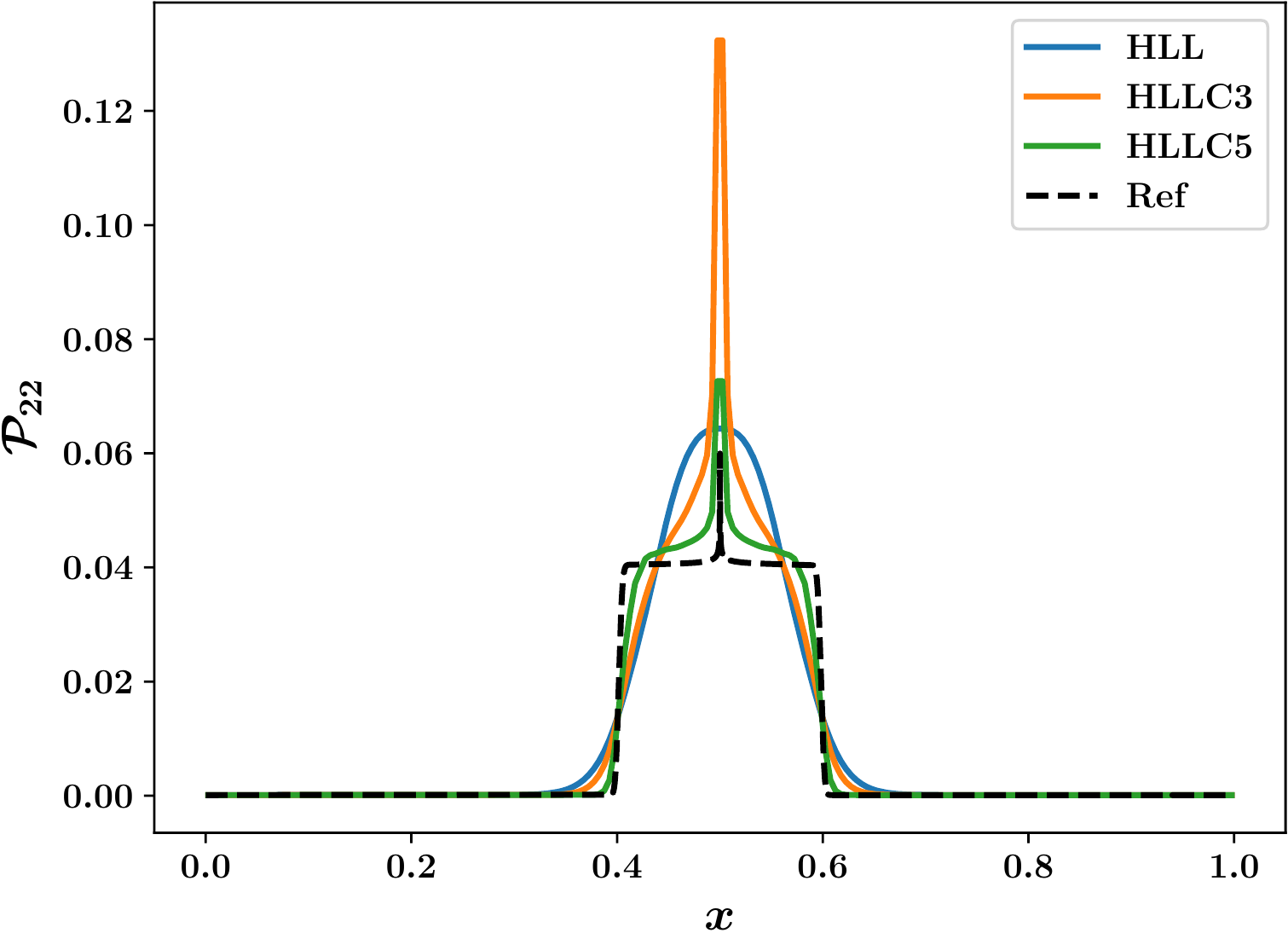} &
         \includegraphics[width=0.40\textwidth]{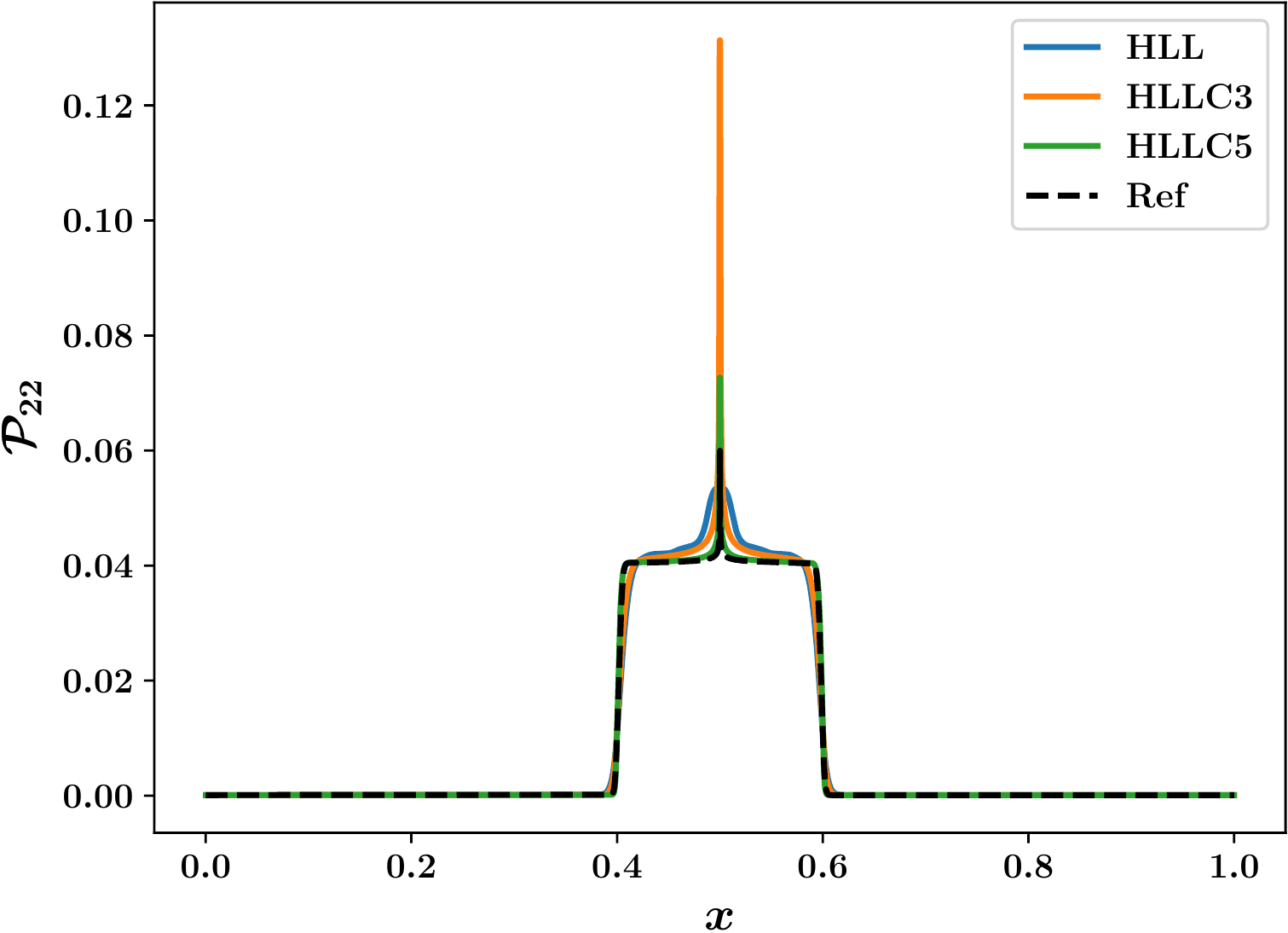} \\
			(a) $\p_{22}$, 200 cells & $\p_{22}$, 2000 cells 
		\end{tabular}
	\end{center}
	\caption{Plots of y-velocity, stress tensor components $\p_{12}$ and $\p_{22}$ obtained using  second-order scheme with HLL, HLLC3 and HLLC5 solvers for 200 and 2000 cells with the reference solution from~\cite{Bhole2019} for Example \ref{sec:shear}.}
	\label{fig:shearO2b}
\end{figure}

\subsection{1-D dam break problem}
\label{sec:dambreak}
The dam break problem models a situation where a dam gate is suddenly opened. The  initial condition has a jump only for the depth variable and all the remaining variables are constant in space. We solve the Riemann problem without source terms in which the initial velocity is zero everywhere and the components of the stress tensor are $\p_{11}  = \p_{22} = 10^{-4} \ m^2/s^2$, $\p_{12}=0 \ m^2/s^2$,  while the initial depth has a discontinuity given by
\[
h = \begin{cases}
0.02 \ m & x < 0.5 \ m \\
0.01 \ m & x > 0.5 \ m
\end{cases}
\]
The solution contains one rarefaction wave and one shock wave separated by a contact discontinuity. In this context, the stress tensor is initially very small and almost at the limit where there is no shear and det(P) is also small (here it is $10^{-8}$). The first order results are shown in Figure~(\ref{fig:dambreakO1b}) and second order results are shown in Figure~(\ref{fig:dambreakO2b}). The numerical solutions show some differences depending on the Riemann solver used. The HLLC3 and the HLLC5 solvers converge asymptotically to the same solution. However, the HLL solver converges asymptotically to a slightly different solution. This is clearly visible in the shape of $\p_{11}$. Moreover, none of these solutions matches with the solution obtained in previous works~\cite{Bhole2019}. Such differences can be expected since we are using different  jump conditions for the non-conservative system. There are also differences between the Riemann solvers developed in this work. The HLLC3 and HLLC5 solvers converge to the same solution under grid refinement while the HLL solver shows some differences compared to these two. In this sense, the inclusion of the intermediate linear waves in the Riemann solver seems to be important.

\begin{figure}
	\begin{center}
		\begin{tabular}{cc}
			\includegraphics[width=0.40\textwidth]{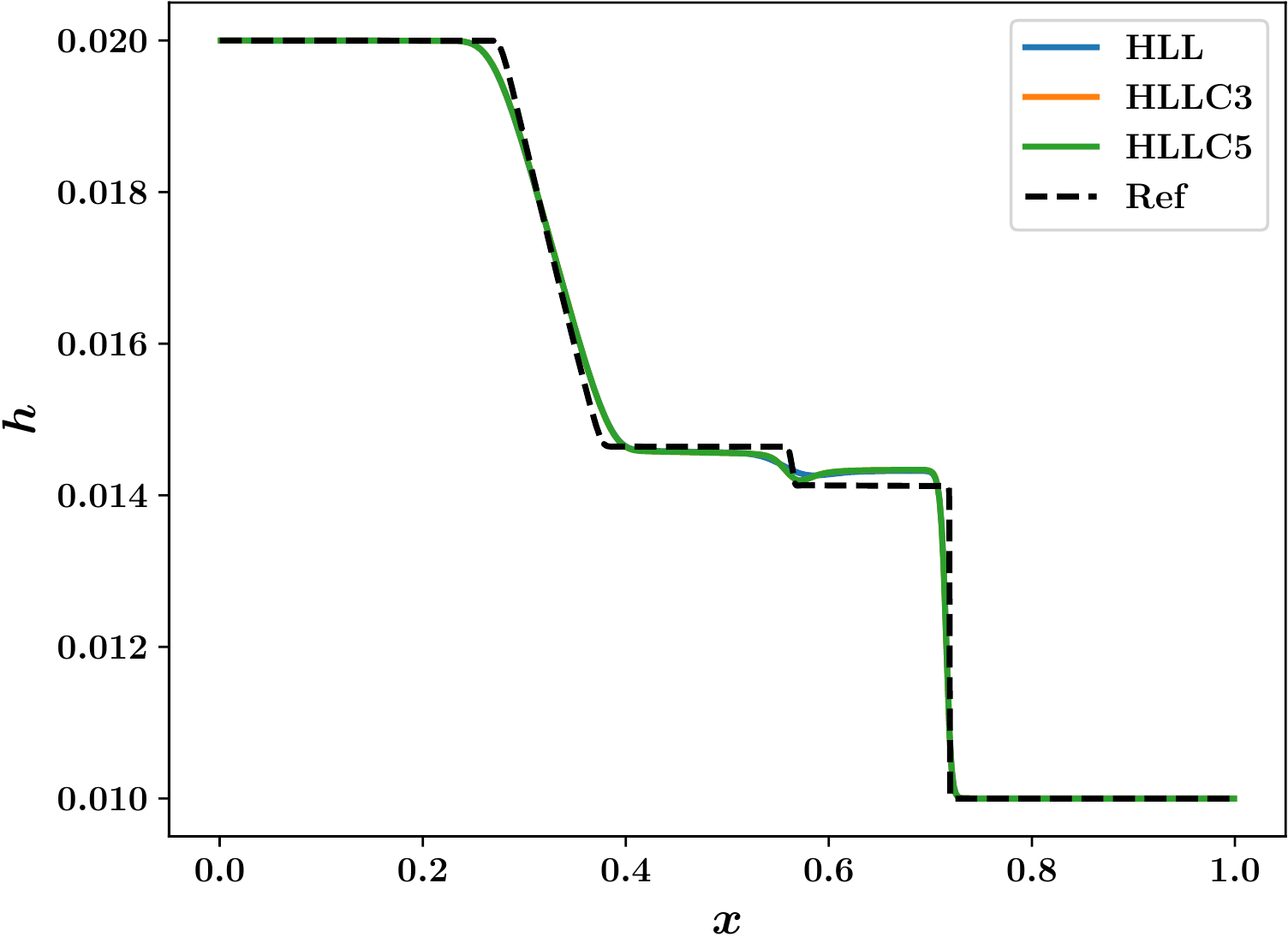} &
			\includegraphics[width=0.40\textwidth]{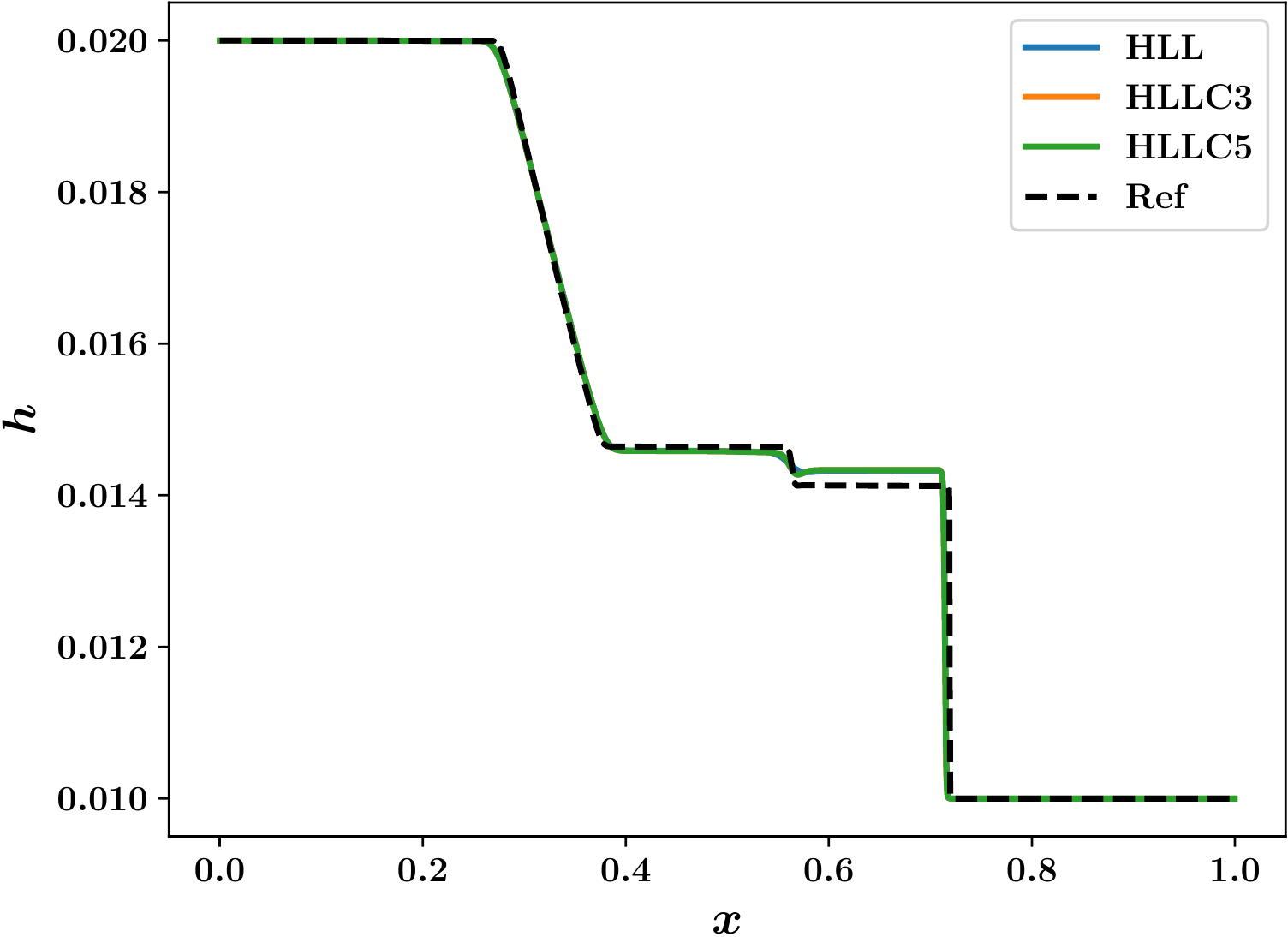} \\
			(a) $h$, 500 cells & (b) $h$, 2000 cells \\
			\includegraphics[width=0.40\textwidth]{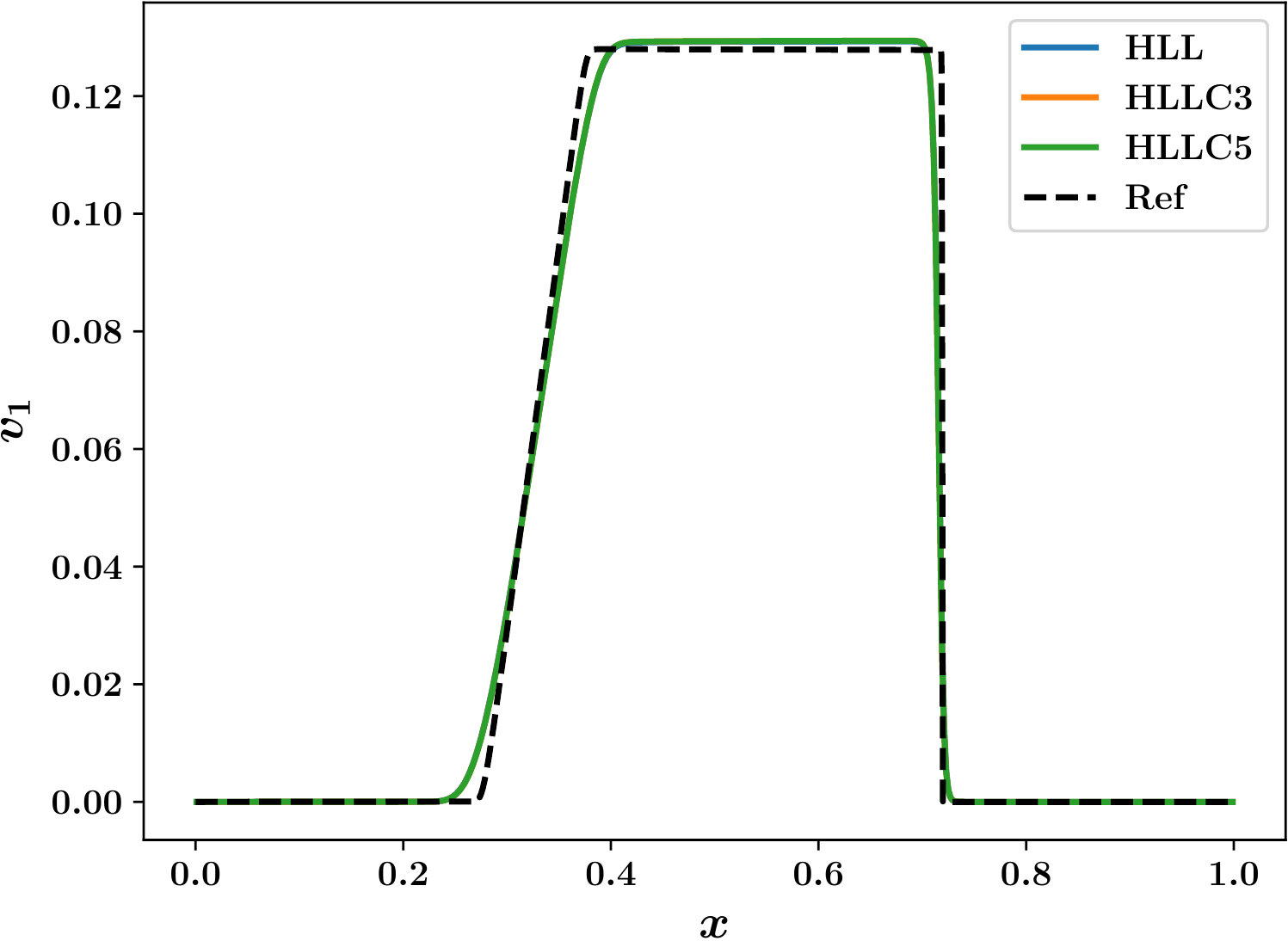} &
			\includegraphics[width=0.40\textwidth]{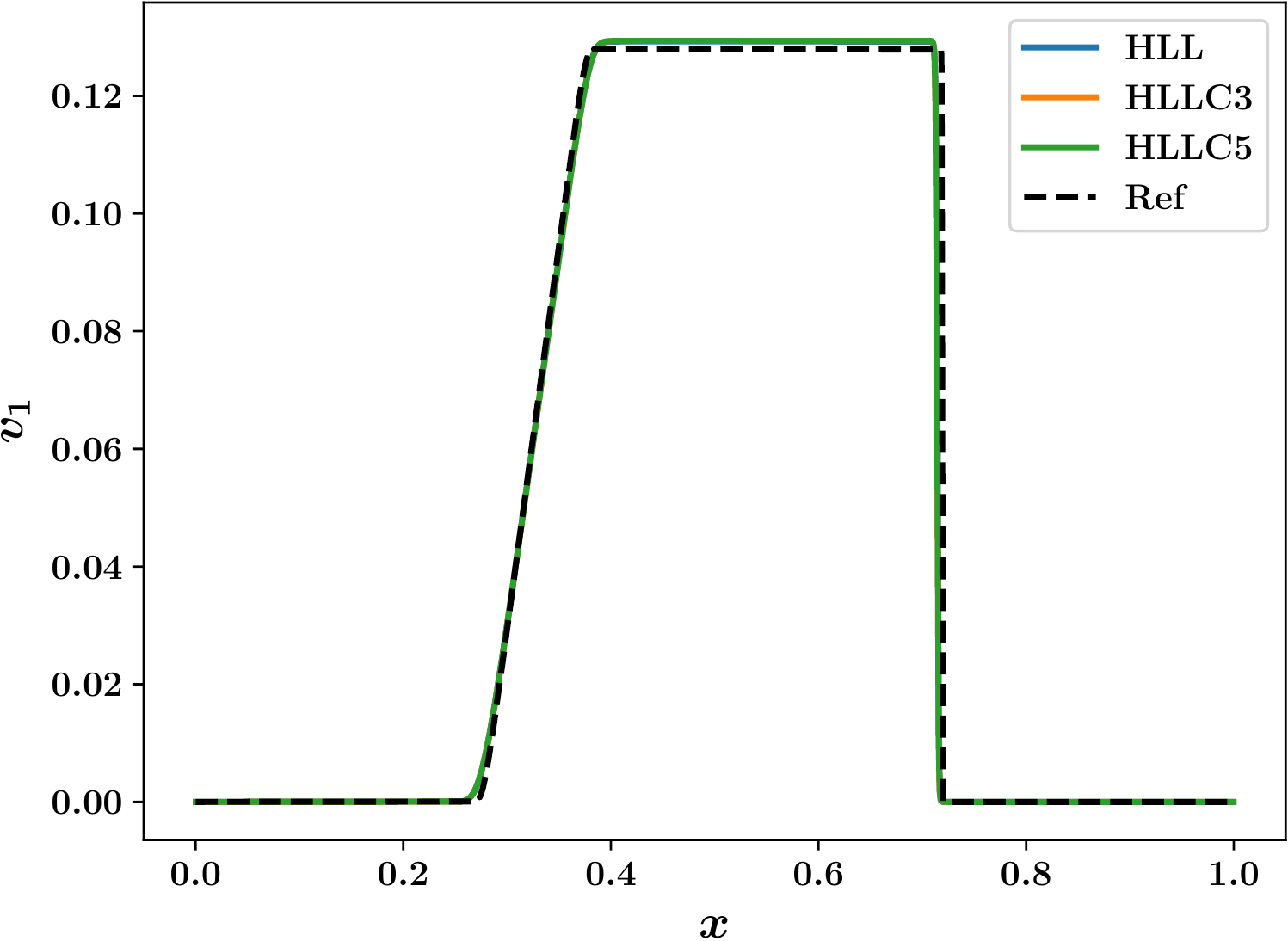} \\
			(a) $v_1$, 500 cells & (b) $v_1$, 2000 cells \\
			\includegraphics[width=0.40\textwidth]{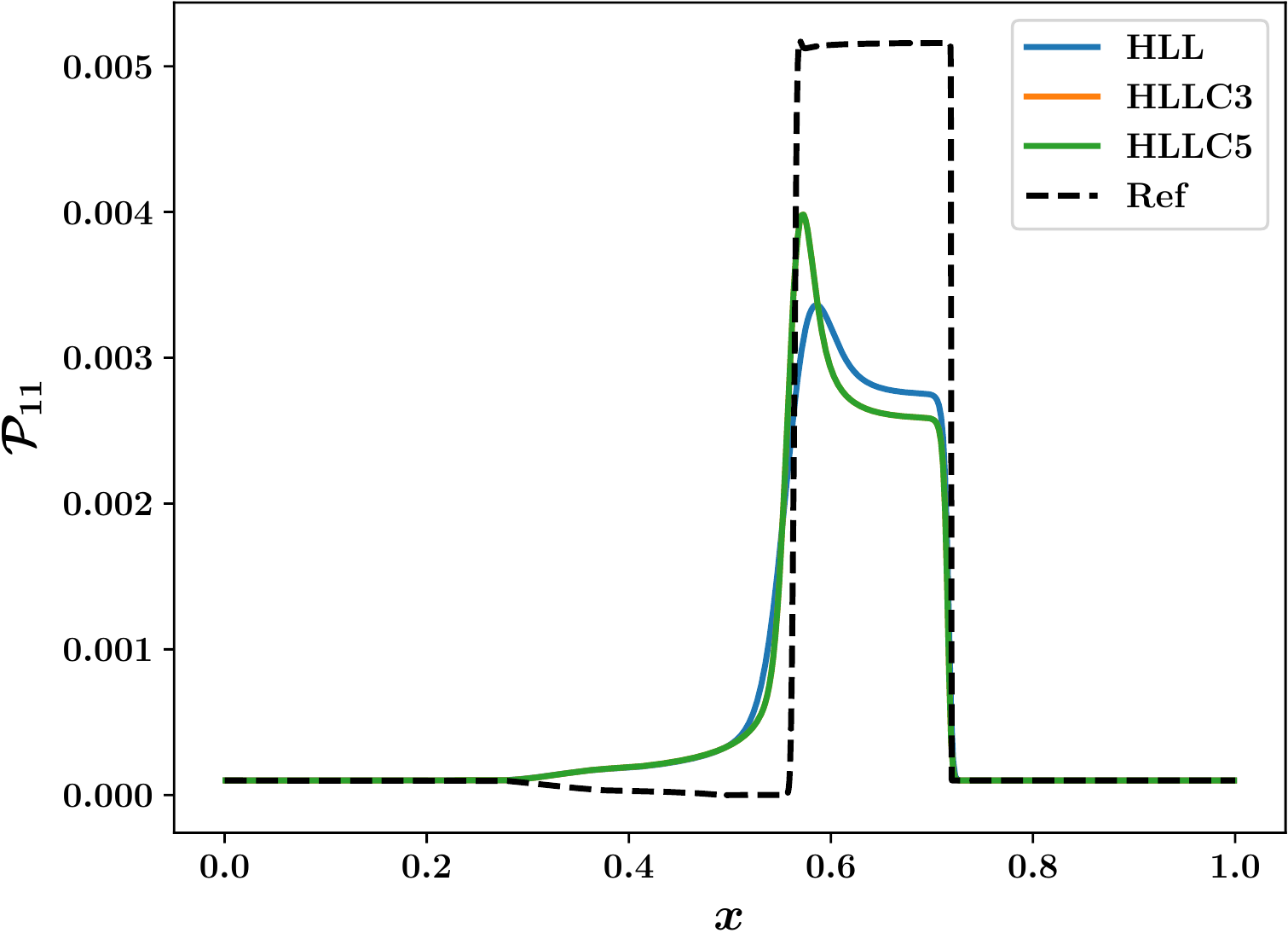} &
			\includegraphics[width=0.40\textwidth]{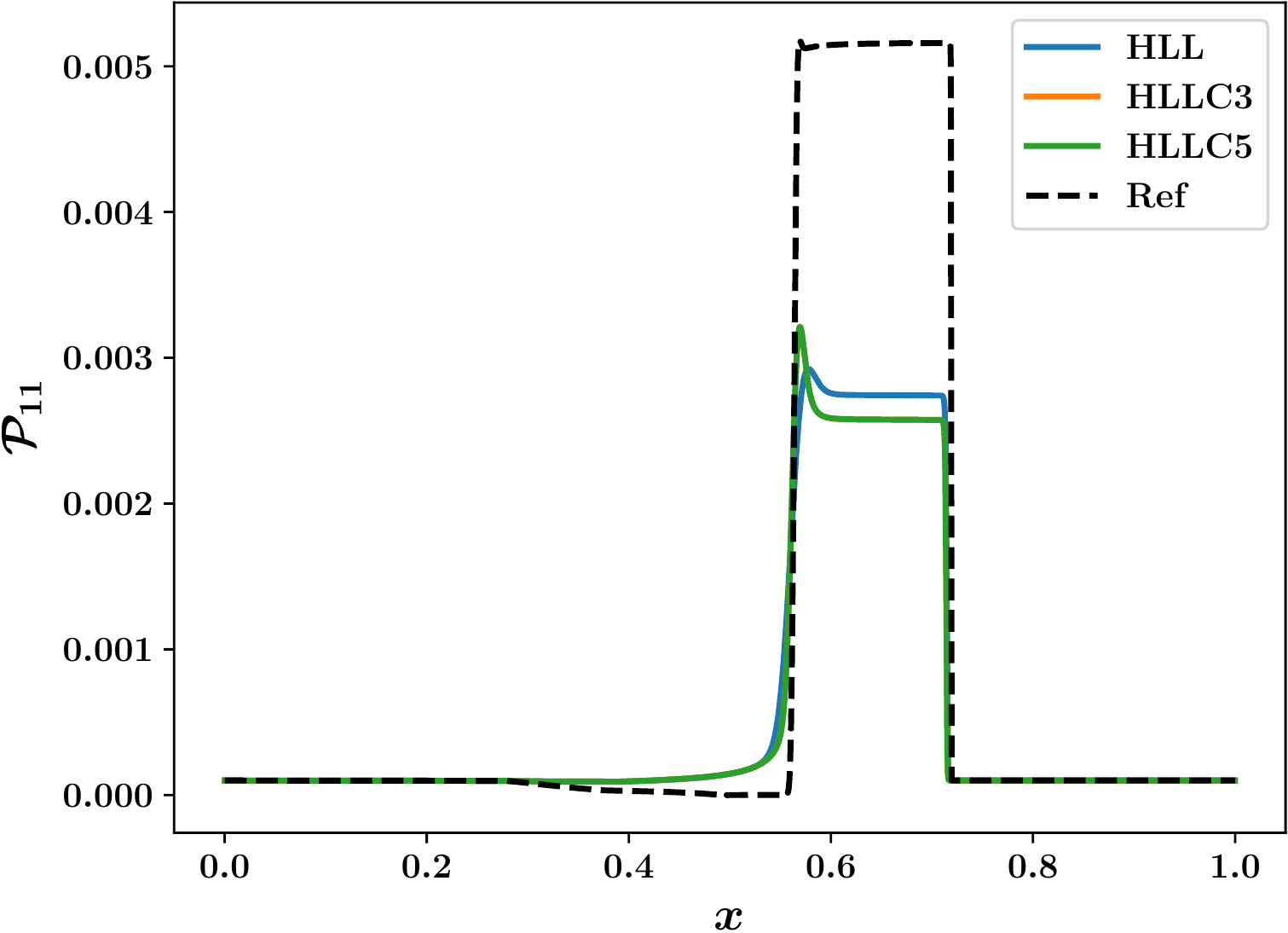}\\
			(a) $\p_{11}$, 500 cells & (b) $\p_{11}$, 2000 cells 
		\end{tabular}
	\end{center}
	\caption{1-D dam break problem from Section~(\ref{sec:dambreak}). Plots of water depth $h$, x-velocity $v_1$, and stress component $\p_{11}$ obtained using first order scheme with HLL, HLLC3 and HLLC5 solvers for 500 and 2000 cells compared with reference solution from~\cite{Bhole2019}.}
	\label{fig:dambreakO1b}
\end{figure}

\begin{figure}
	\begin{center}
		\begin{tabular}{ccc}
         \includegraphics[width=0.40\textwidth]{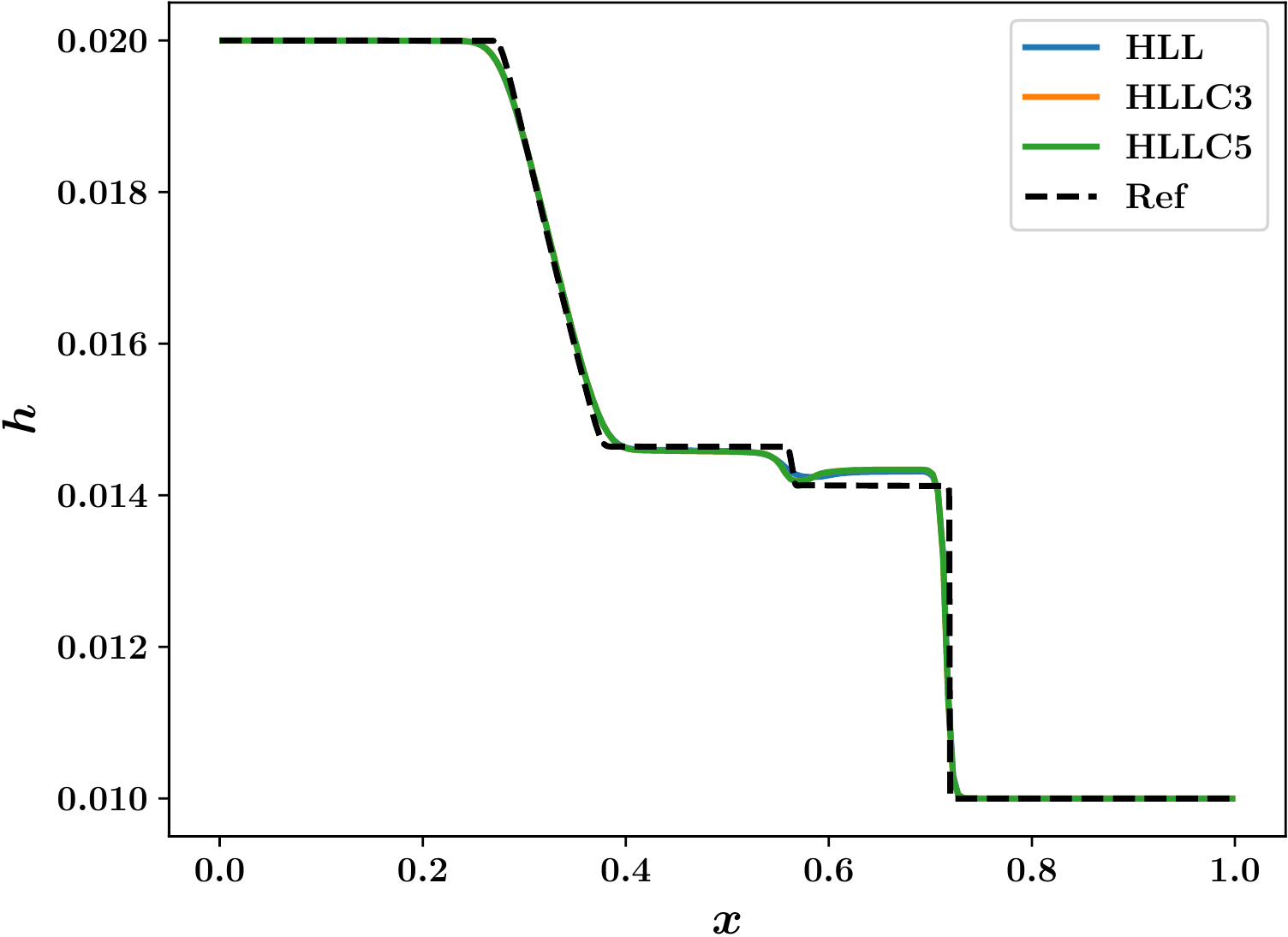} &
         \includegraphics[width=0.40\textwidth]{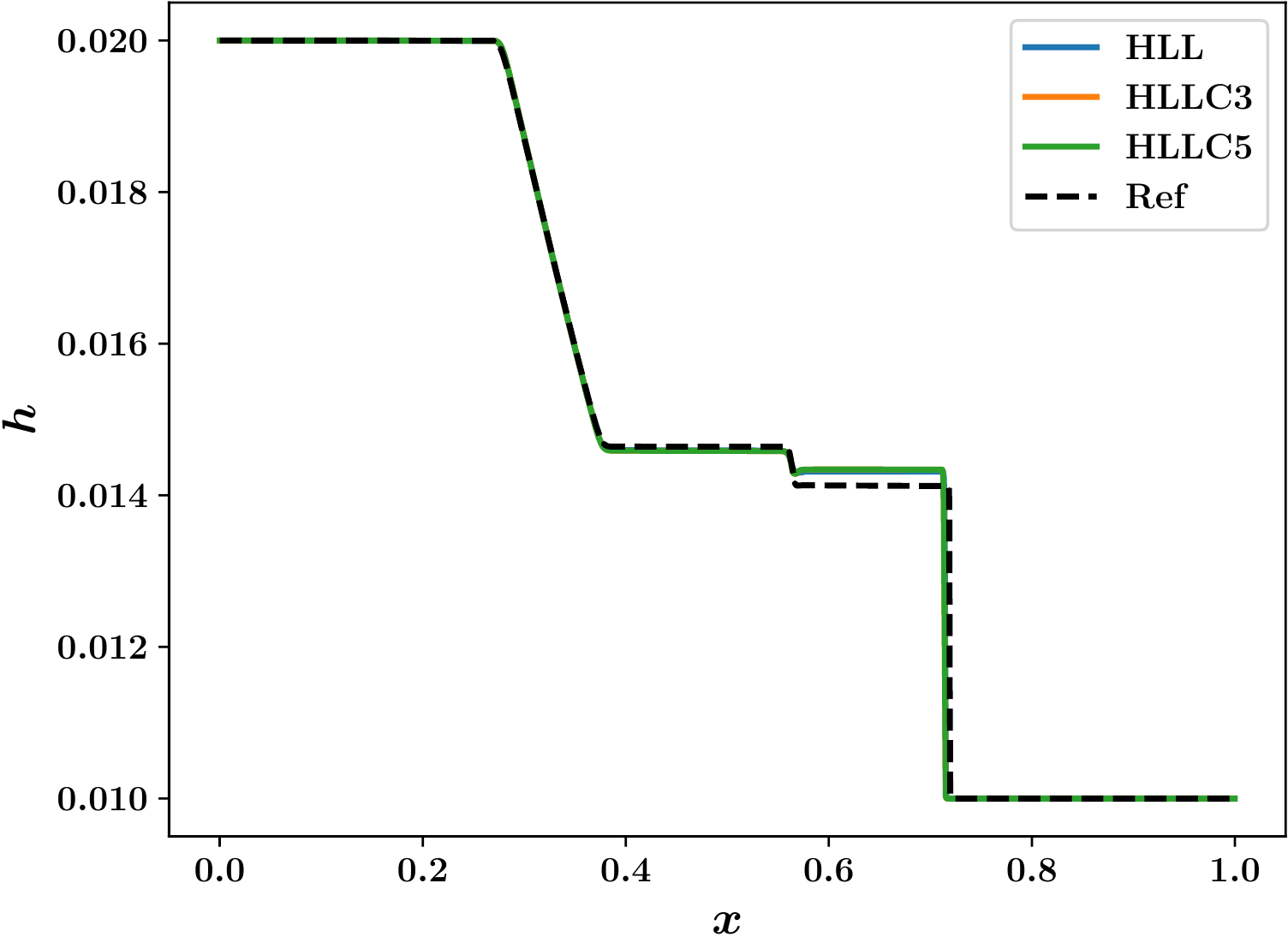} \\
         (a) $h$, 200 cells & (b) $h$, 2000 cells \\
         \includegraphics[width=0.40\textwidth]{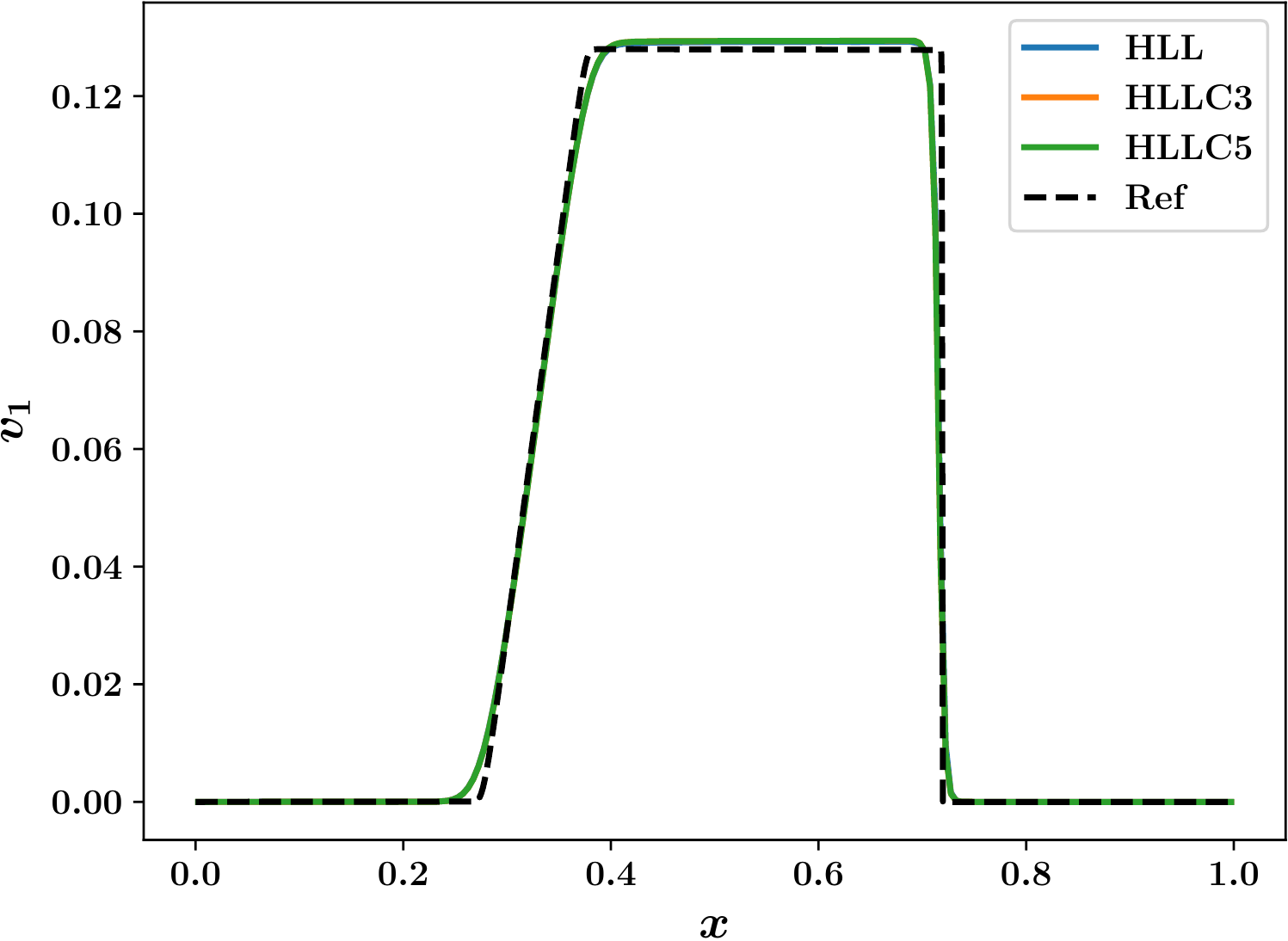} &
         \includegraphics[width=0.40\textwidth]{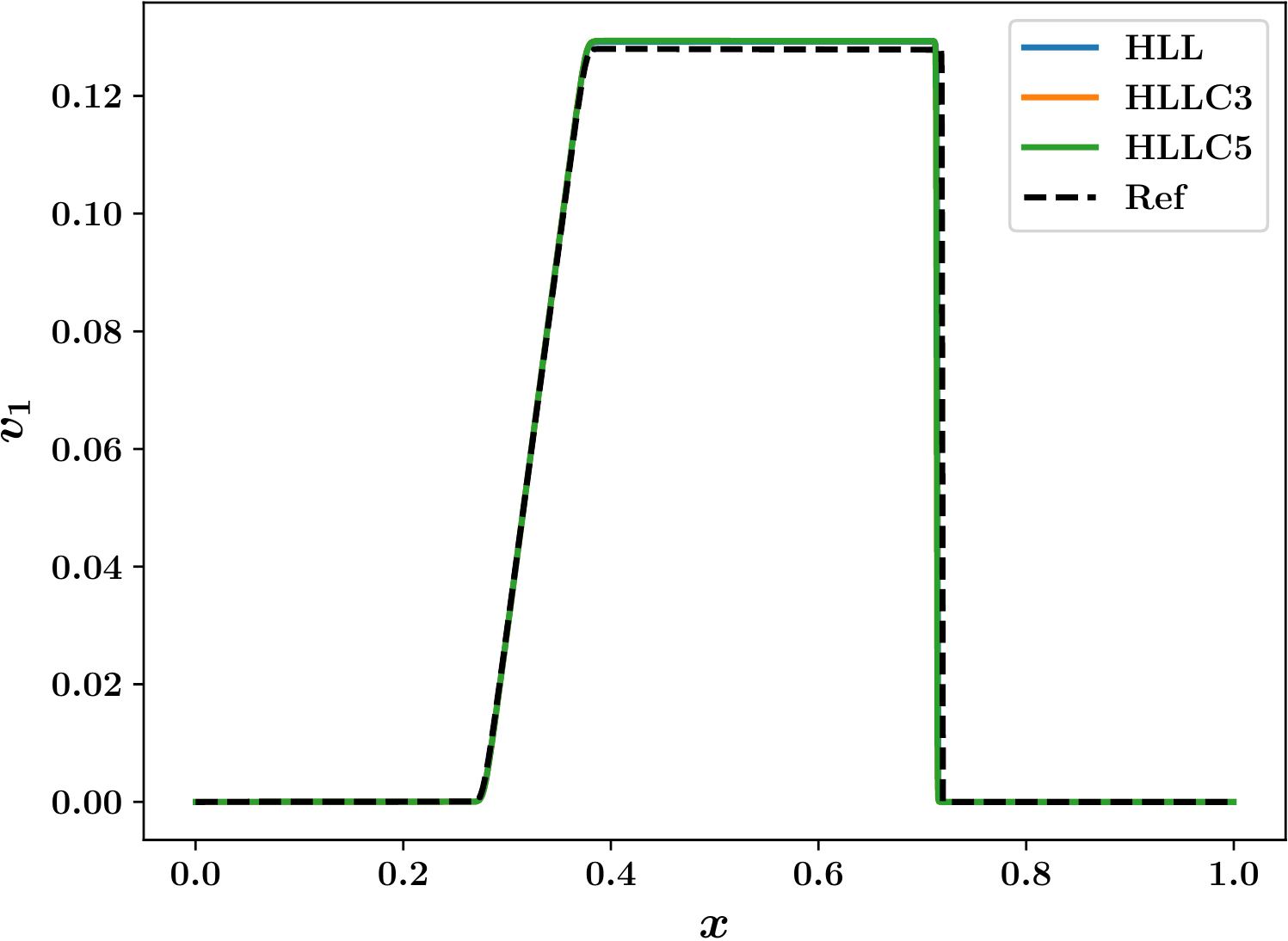} \\
         (a) $v_1$, 200 cells & (b) $v_1$, 2000 cells \\
         \includegraphics[width=0.40\textwidth]{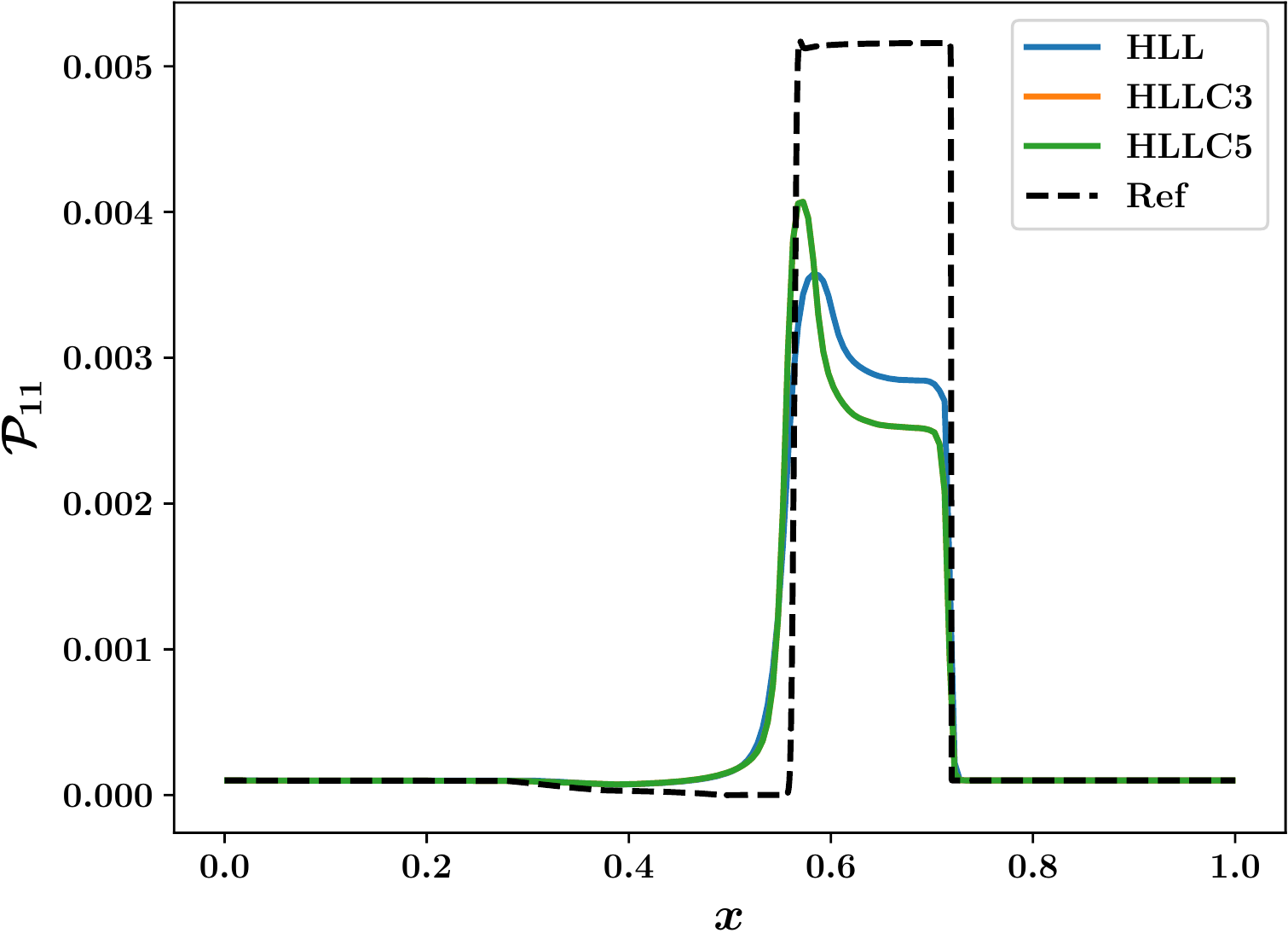} &
         \includegraphics[width=0.40\textwidth]{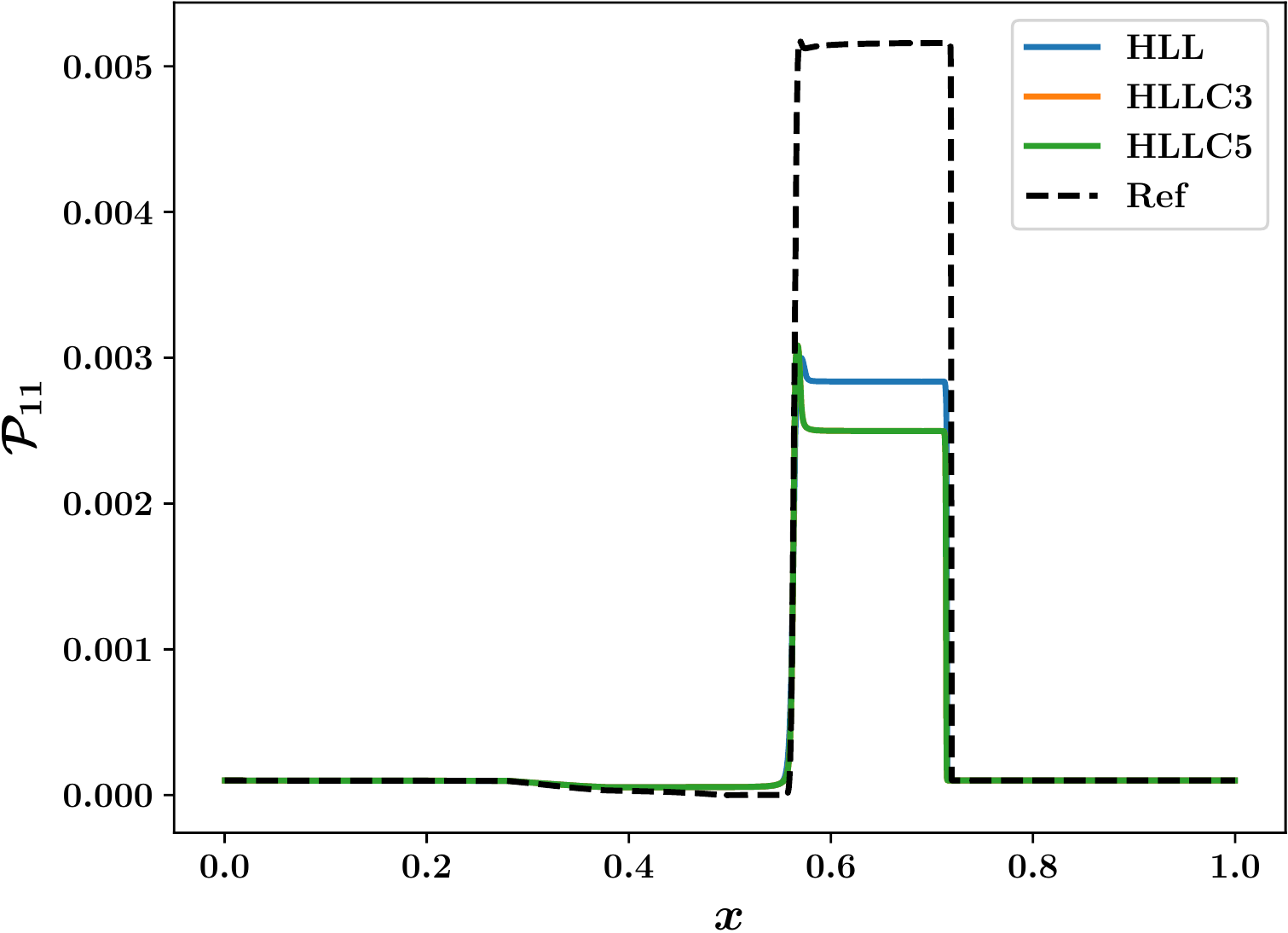} \\
         (a) $\p_{11}$, 200 cells & (b) $\p_{11}$, 2000 cells  
		\end{tabular}
	\end{center}
	\caption{1-D dam break problem from Section~(\ref{sec:dambreak}). Plots of water depth $h$ x-velocity $v_1$ and stress component $\p_{11}$ obtained using second-order scheme with HLL, HLLC3 and HLLC5 solvers for 200 and 2000 cells compared with the reference solution from~\cite{Bhole2019}.}
	\label{fig:dambreakO2b}
\end{figure}

\subsection{1-D modified dam break problem}
\label{sec:fivewave}
This is a modified dam break problem where additional jump is initially introduced in the transverse velocity
component. The Riemann problem is solved without any source terms and the initial condition is given by
\[
h(x,0) = \begin{cases}
0.01 & x < 0.5 \\
0.02 & x > 0.5 \end{cases}, \qquad v_2(x,0) = \begin{cases}
\ \ 0.2 & x < 0.5 \\
-0.2 & x > 0.5 \end{cases}
\]
while the remaining quantities are constant in space and are given by $v_1=0.1$, $\p_{11} = \p_{22} = 4 \times 10^{-2}$, $\p_{12} = 10^{-8}$. The solution shown in Figure~(\ref{fig:fivewaveO2b}) is obtained from the second order scheme and behaves similar to the dam break problem with two additional contact discontinuities. However, contact discontinuities move slowly than for the dam break problem and  det($\p$) is less close to zero (of order $10^{-3}$). With this modified problem, all the Riemann solvers asymptotically converge to the same solution. As expected, the resolution is improved when the approximate solver contains more physical waves.

\begin{figure}
	\begin{center}
		\begin{tabular}{ccc}
			\includegraphics[width=0.40\textwidth]{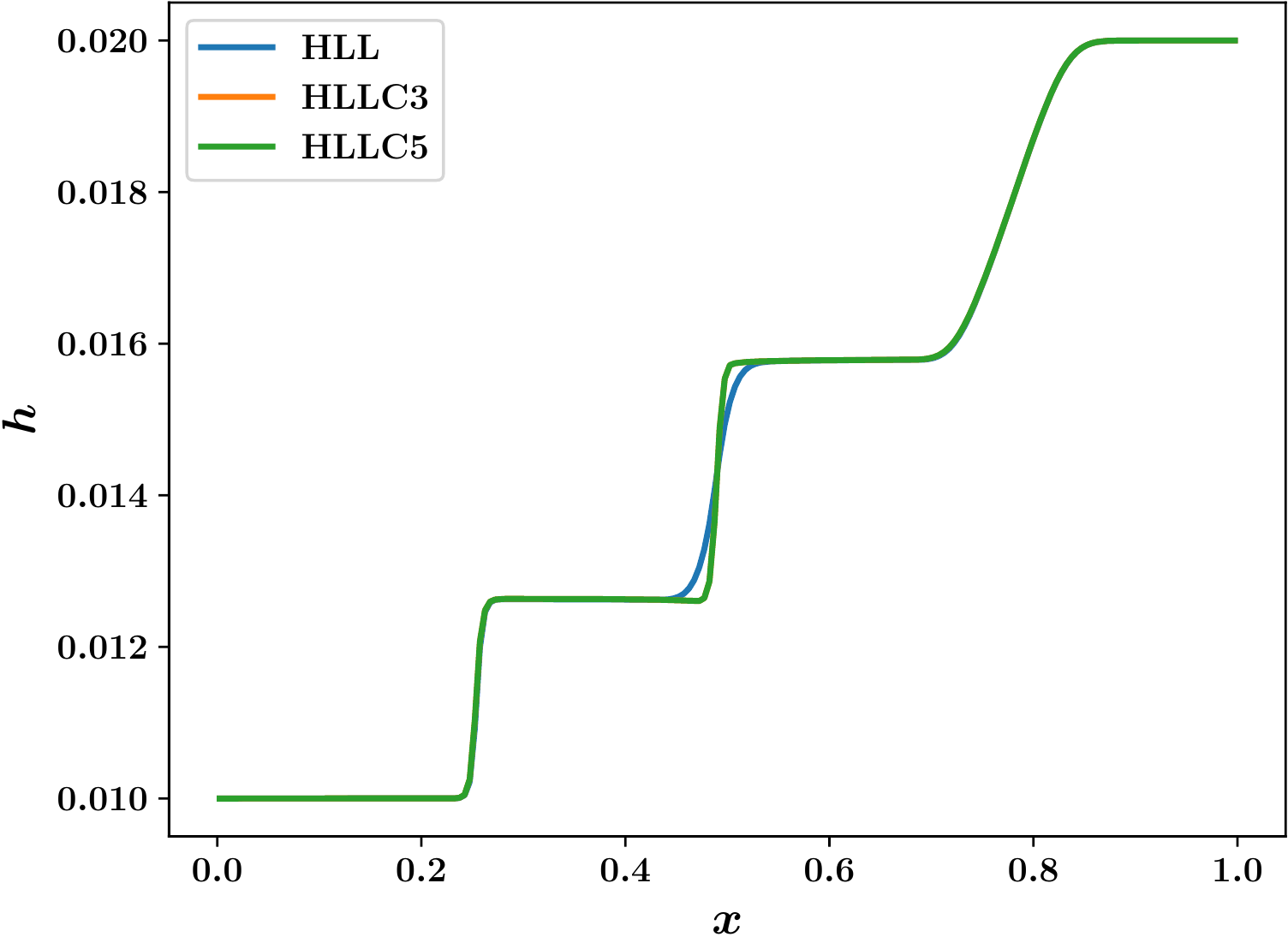} &
         \includegraphics[width=0.40\textwidth]{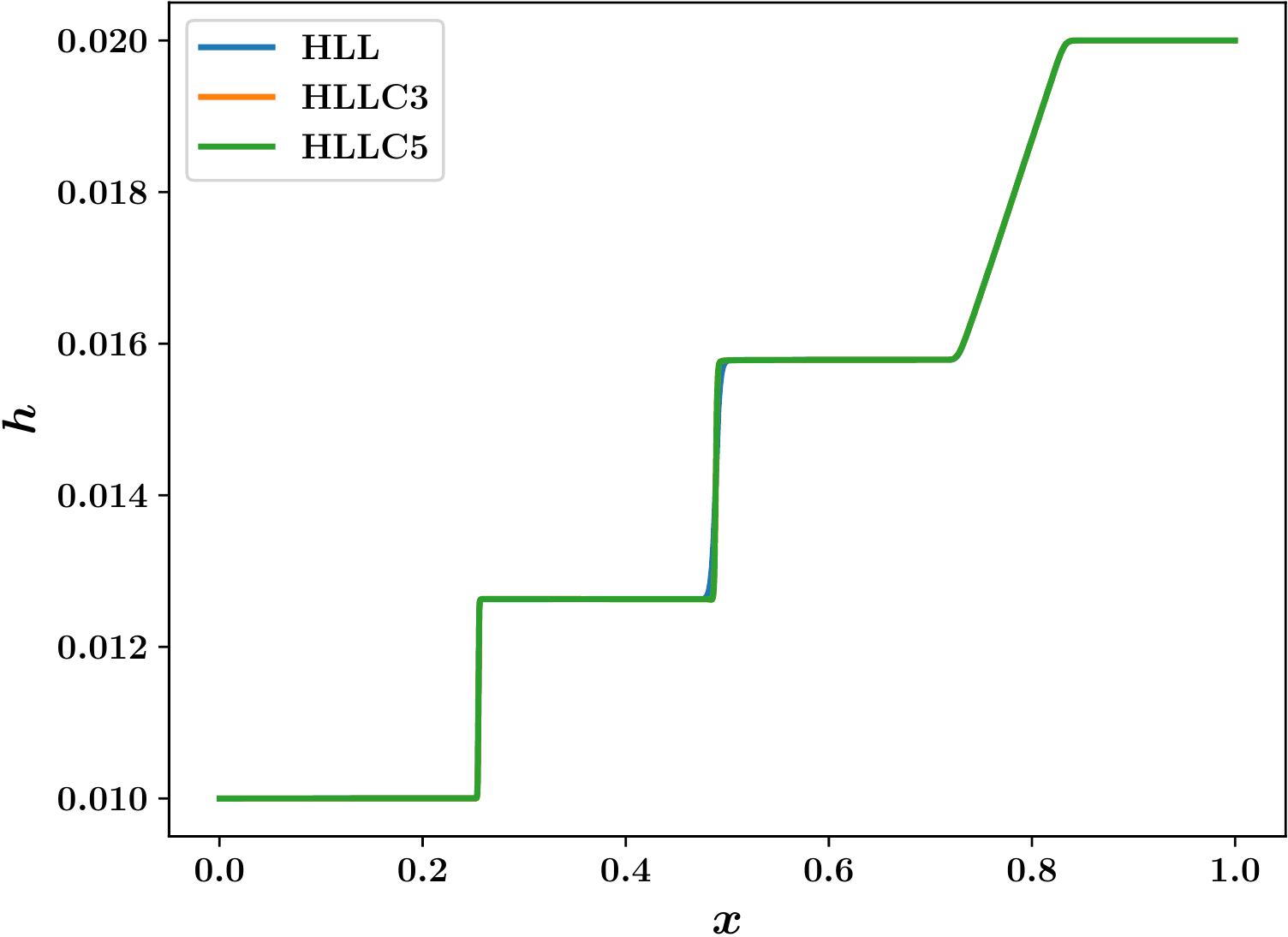} \\
			(a) $h$, 200 cells & (b) $h$, 2000 cells \\
			\includegraphics[width=0.40\textwidth]{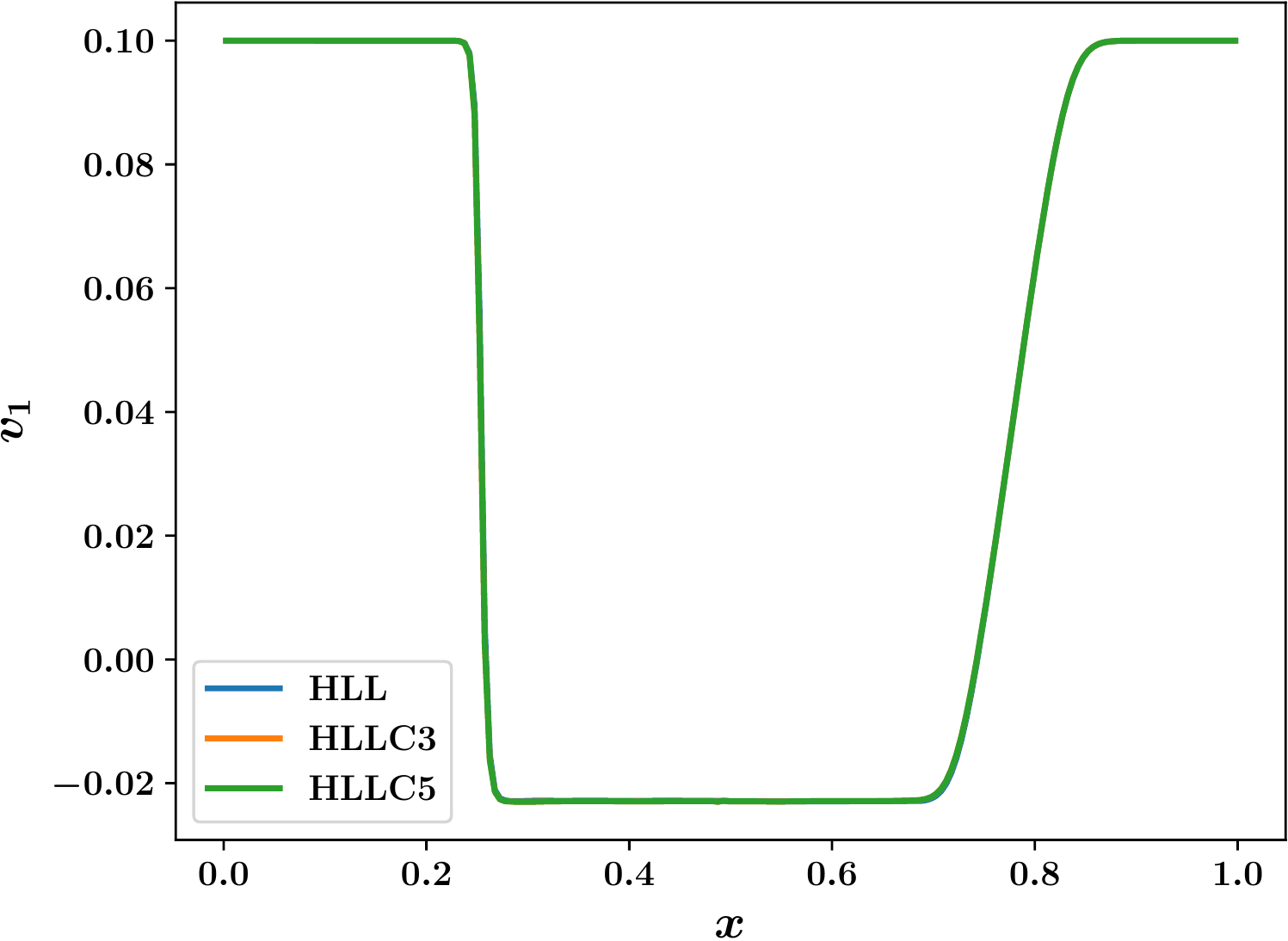} &
         \includegraphics[width=0.40\textwidth]{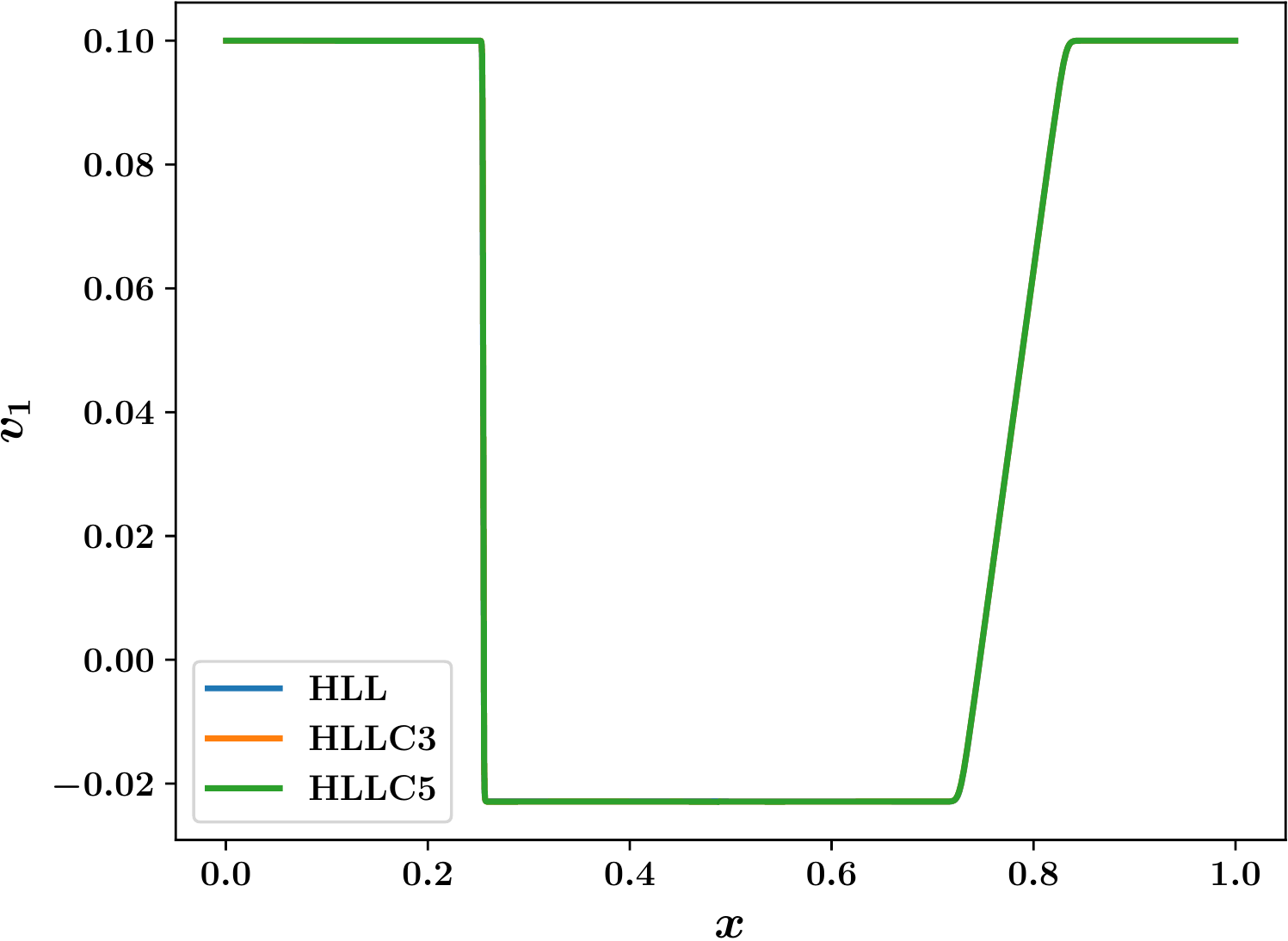} \\
			(a) $v_1$, 200 cells & (b) $v_1$, 2000 cells \\
			\includegraphics[width=0.40\textwidth]{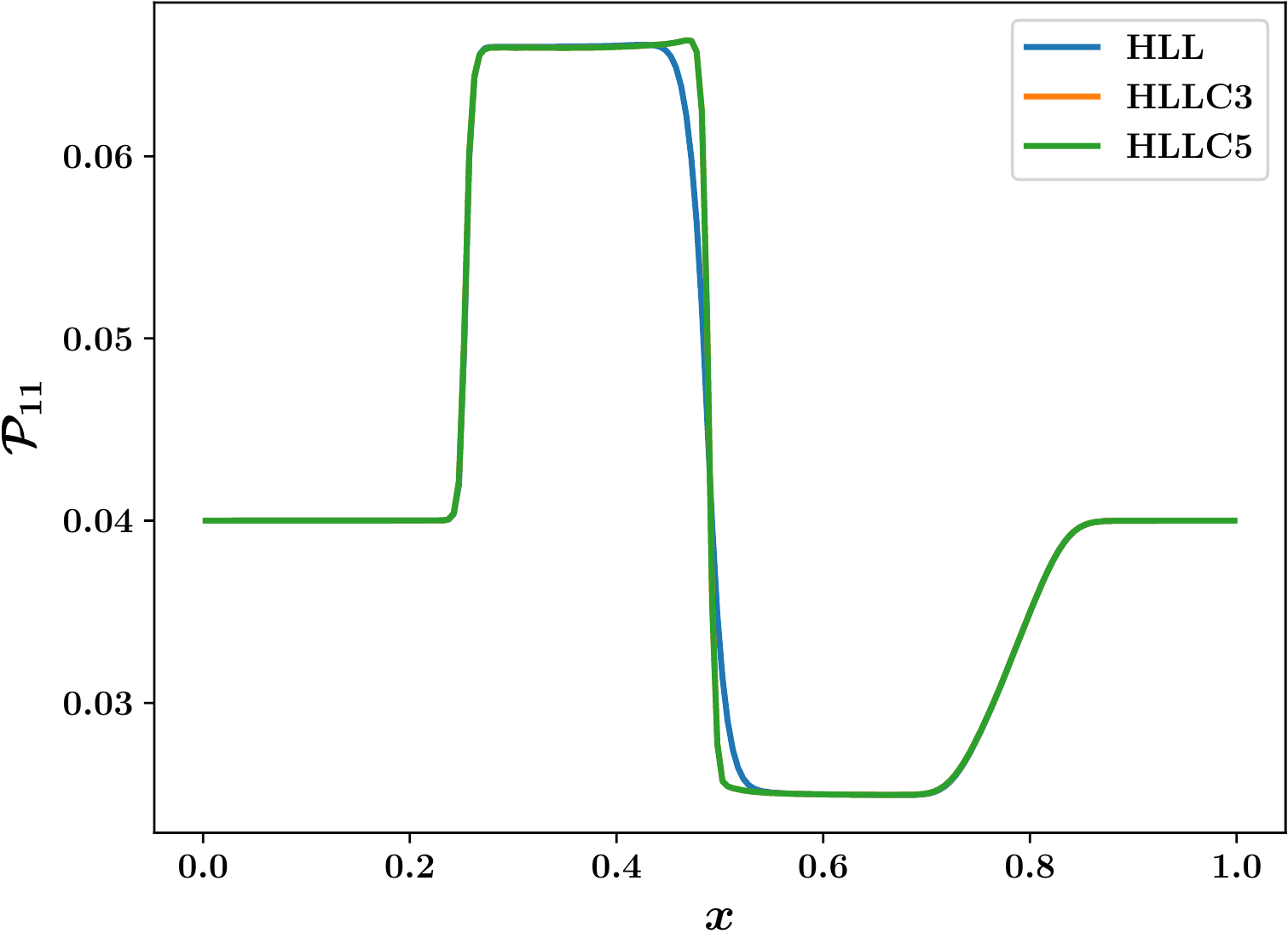} &
         \includegraphics[width=0.40\textwidth]{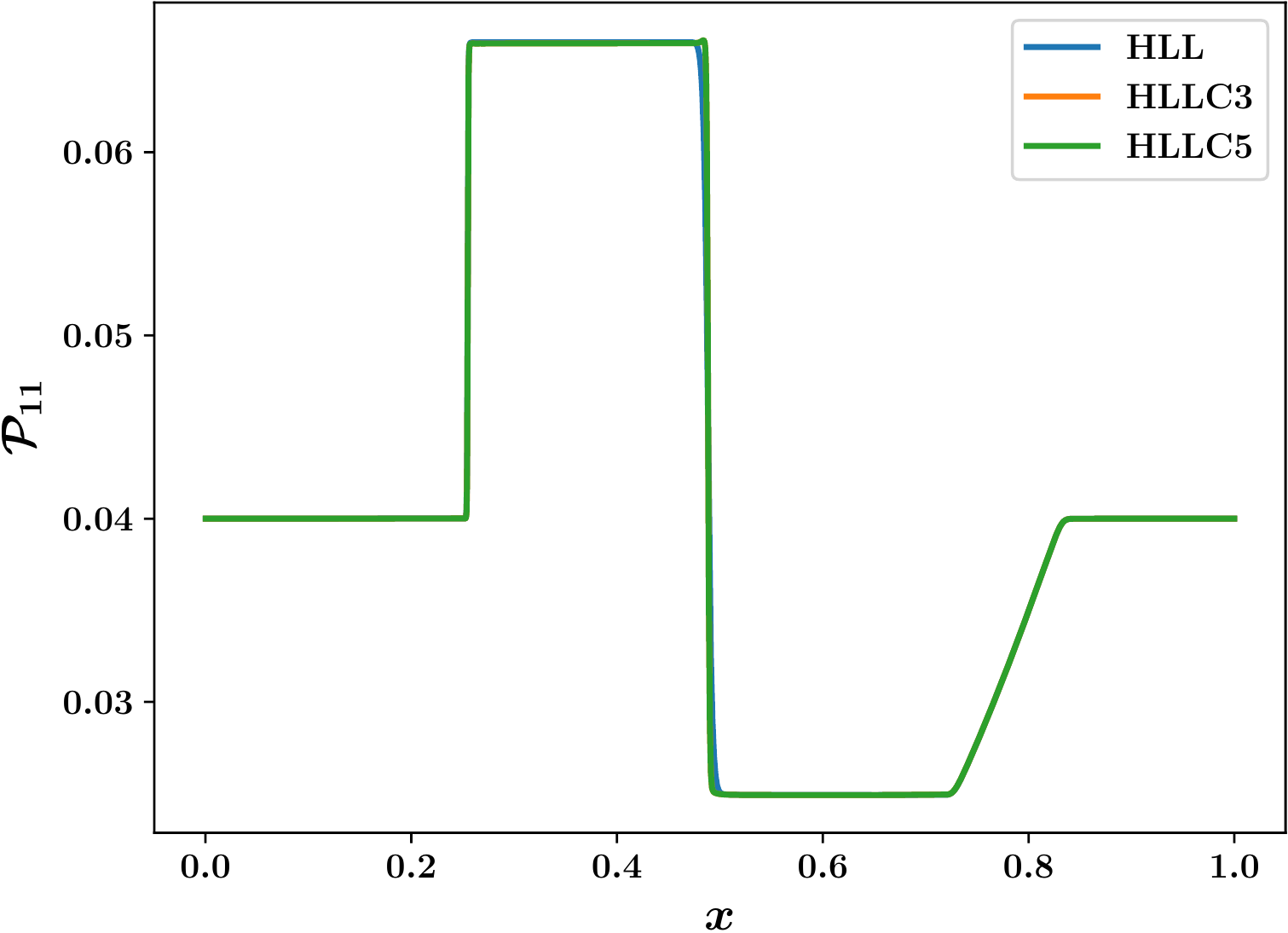} \\
			(a) $\p_{11}$, 200 cells & (b) $\p_{11}$, 2000 cells 
		\end{tabular}
	\end{center}
	\caption{1-D modified dam break problem from Section~(\ref{sec:fivewave}). Plots of water depth $h$ x-velocity $v_1$ and stress component $\p_{11}$ obtained using second order scheme with HLL, HLLC3 and HLLC5 solvers for 200 and 2000 cells.}
	\label{fig:fivewaveO2b}
\end{figure}

\subsection{1-D roll wave problem}
\label{sec:roll1d}
This problem models the flow of a thin layer of liquid down an inclined bottom and we include all the source terms. The initial condition is given by
\[
h(x,0) = h_0 \left[ 1 + a \sin(2\pi x/L_x) \right], \qquad v_1(x,0) = \sqrt{ \frac{g h_0 \tan\theta}{C_f} }, \qquad v_2(x,0) = 0
\]
\[
\p_{11}(x,0) = \p_{22}(x,0) = \half \phi h^2(x,0), \qquad \p_{12}(x,0) = 0
\]
\reva{
The bottom topography is given by $b = - x \tan\theta$ and the boundary conditions are periodic. We consider two sets of parameters as in~\cite{Ivanova2017}. In Case 1, the parameters are $\theta = 0.05011 \ rad$ the inclination angle, $C_f = 0.0036$, $h_0 = 7.98 \times 10^{-3} \ m$, $a = 0.05$, $\phi = 22.76 s^{-2}$, $C_r = 0.00035$, $L_x = 1.3 \ m$.  In Case 2, the parameters are $\theta = 0.119528 \ rad$ the inclination angle, $C_f = 0.0038$, $h_0 = 5.33 \times 10^{-3} \ m$, $a = 0.05$, $\phi = 153.501 s^{-2}$, $C_r = 0.002$, $L_x = 1.8 \ m$

The chosen parameters  lead to the formation of 1-D roll waves starting from a uniform flow which has the same structure as in Brock's experiments~\cite{Brock1969,Brock1970}. The water depth is compared with Brock's experimental data in Figure~(\ref{fig:roll1da}) for both Case 1 and Case 2, which shows a hydraulic jump around $x/L_x=1$ but also a  smooth profile immediately behind it which is the signature of a roll wave. The classical shallow water model does not predict this roll wave profile but only gives rise to the hydraulic jump. All the three Riemann solvers yield essentially similar results and the comparison with the experimental data is good. Figure~(\ref{fig:roll1db}) shows a comparison of some other quantities for Case 1 and we observe that all solvers yield essentially the same solution. A 2-D version of this test case is discussed in Section~(\ref{sec:roll2d}).
}
\begin{figure}
\begin{center}
\begin{tabular}{cc}
	\includegraphics[width=0.48\textwidth]{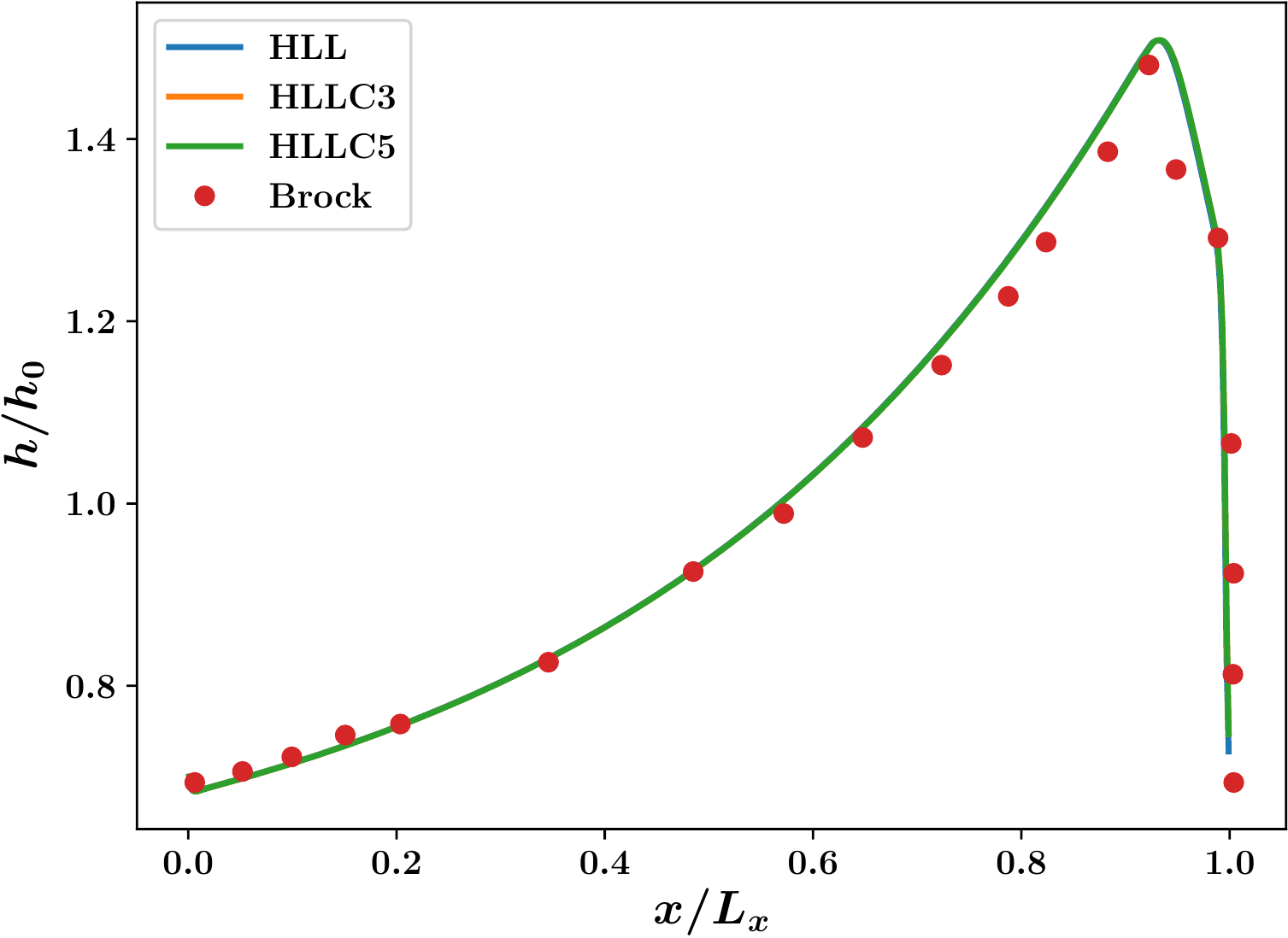} &
	\includegraphics[width=0.48\textwidth]{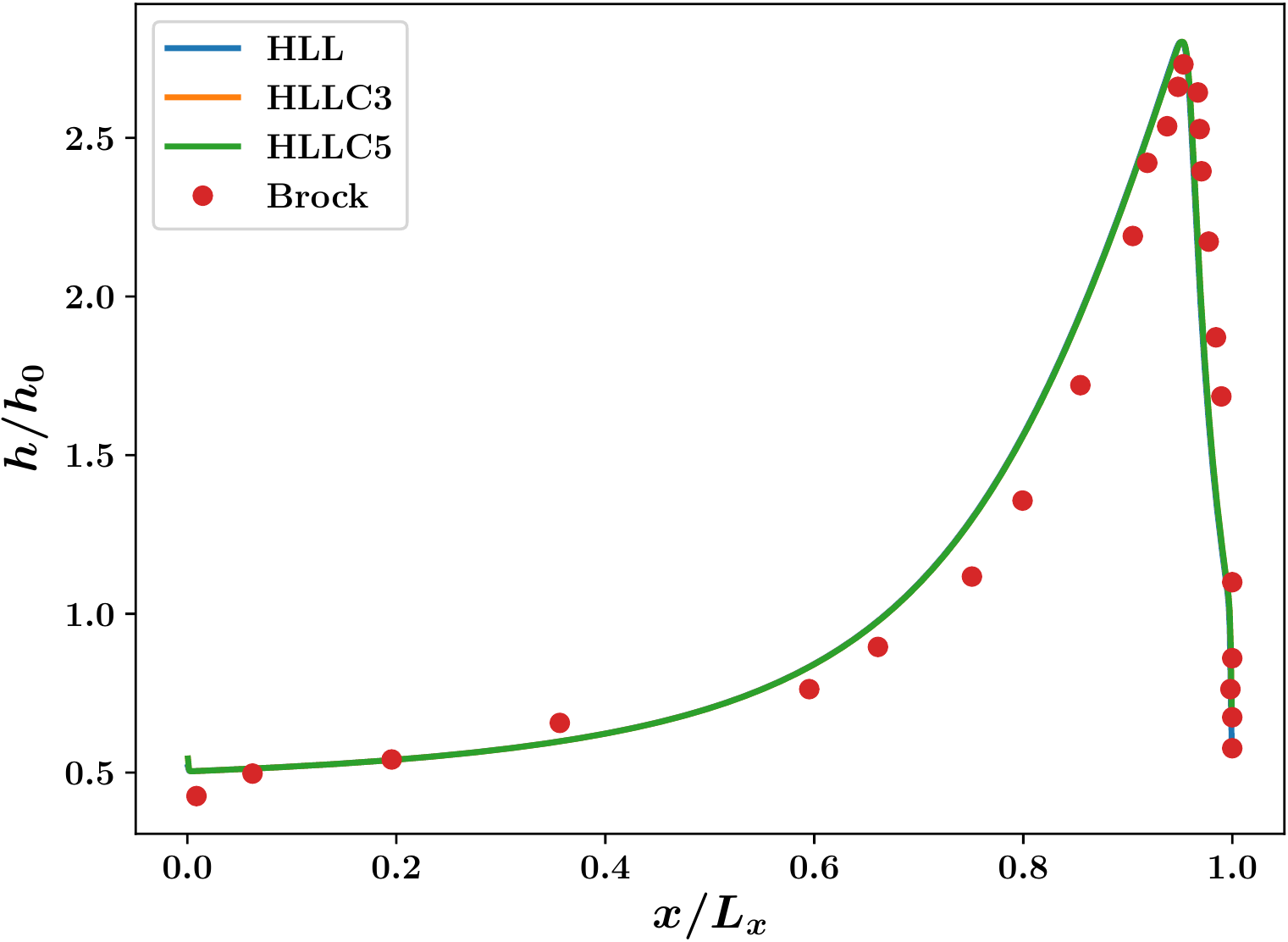} \\
(a) Case 1 & (b) Case 2
\end{tabular}
\end{center}
\caption{\reva{1-D roll wave problem from Section~(\ref{sec:roll1d}). Comparison of water depth obtained using  second-order scheme with HLL, HLLC3 and HLLC5 solvers for 500 cells with Brock's experimental data. (a) Case 1 at time $t=26.99$, (b) Case 2 at time $t=26.35185$}}
\label{fig:roll1da}
\end{figure}

\begin{figure}
	\begin{center}
		\begin{tabular}{cc}
			\includegraphics[width=0.48\textwidth]{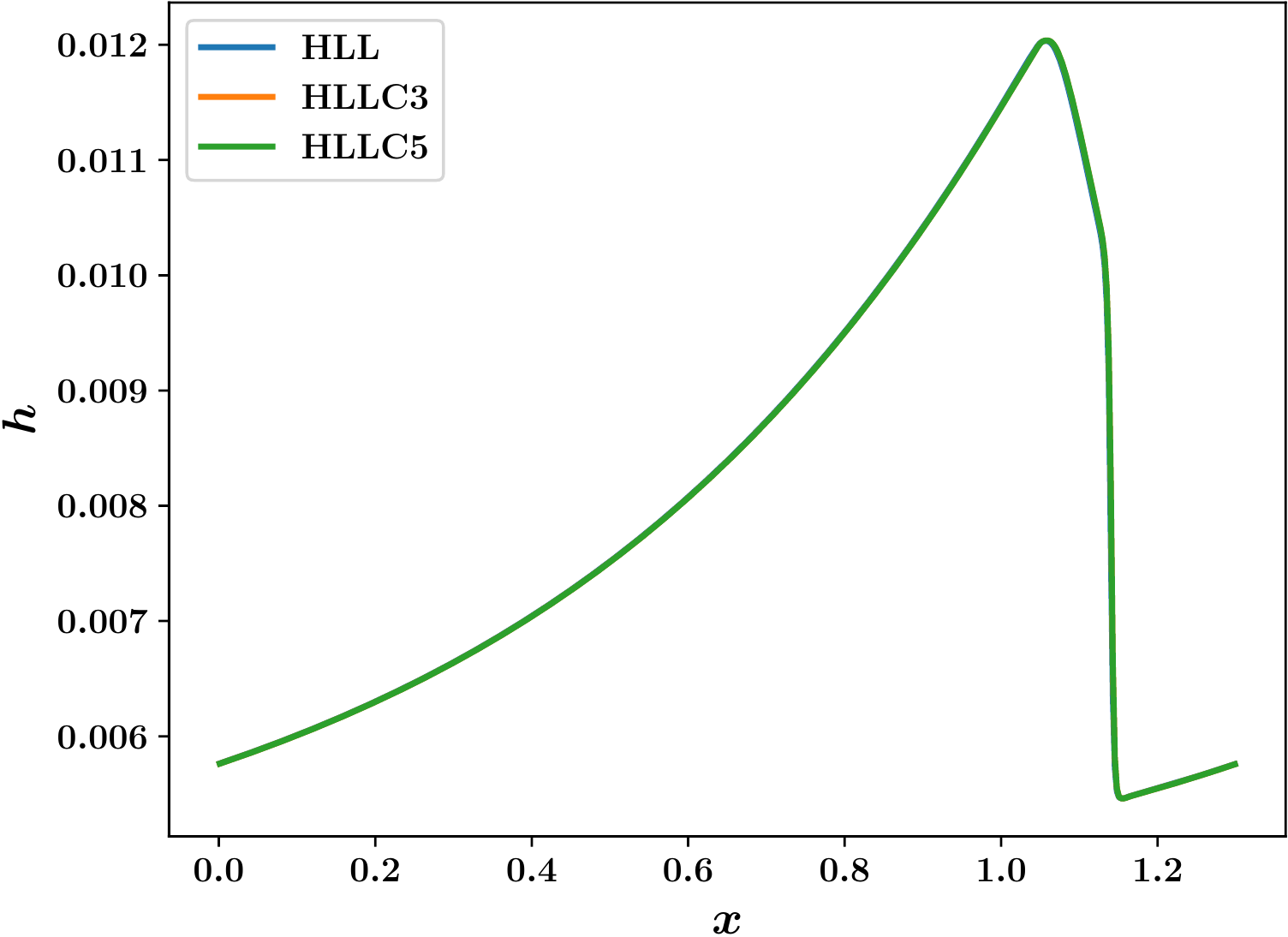} &
			\includegraphics[width=0.48\textwidth]{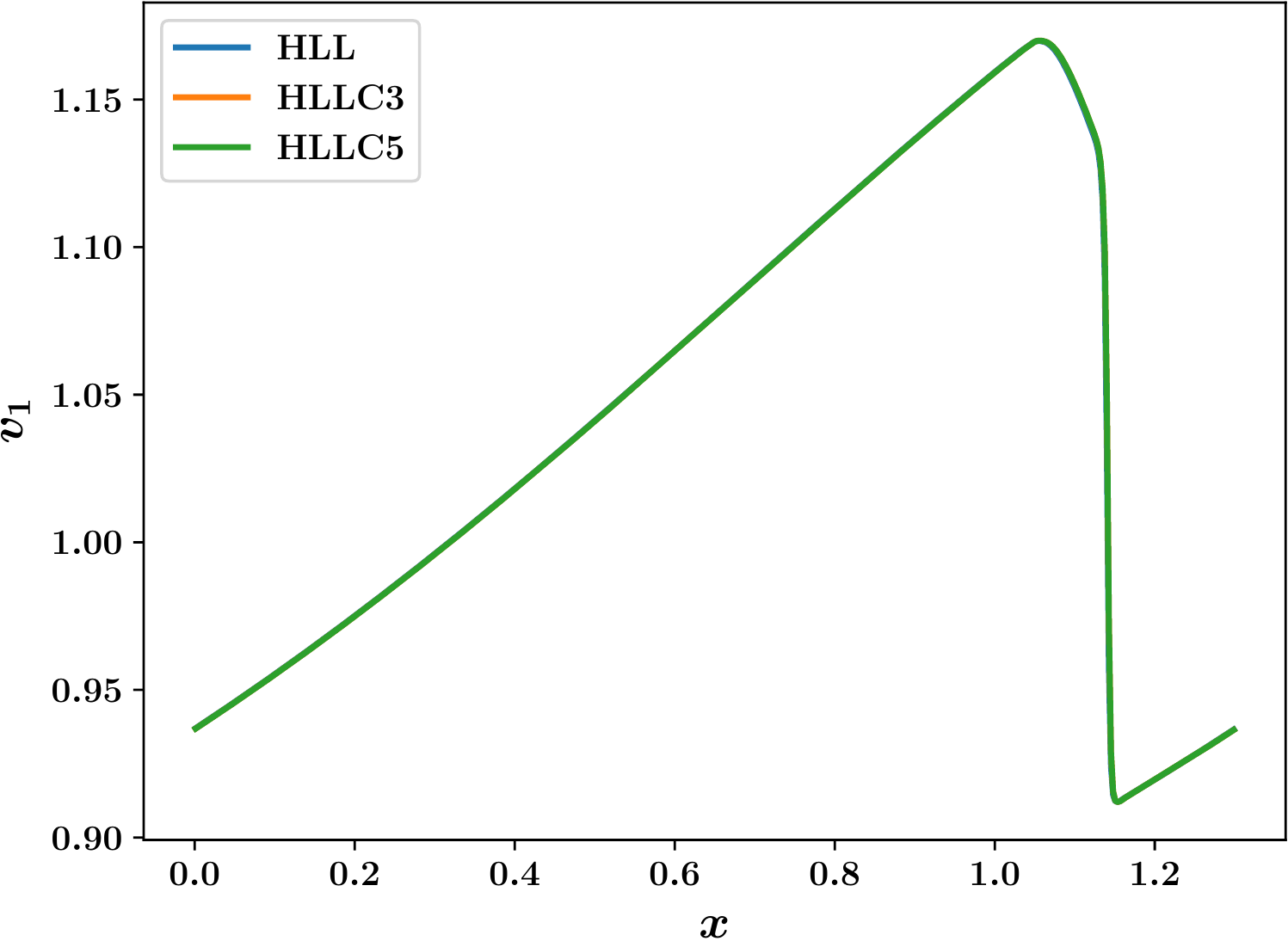} \\
			$h$ & $v_1$  \\
			\includegraphics[width=0.48\textwidth]{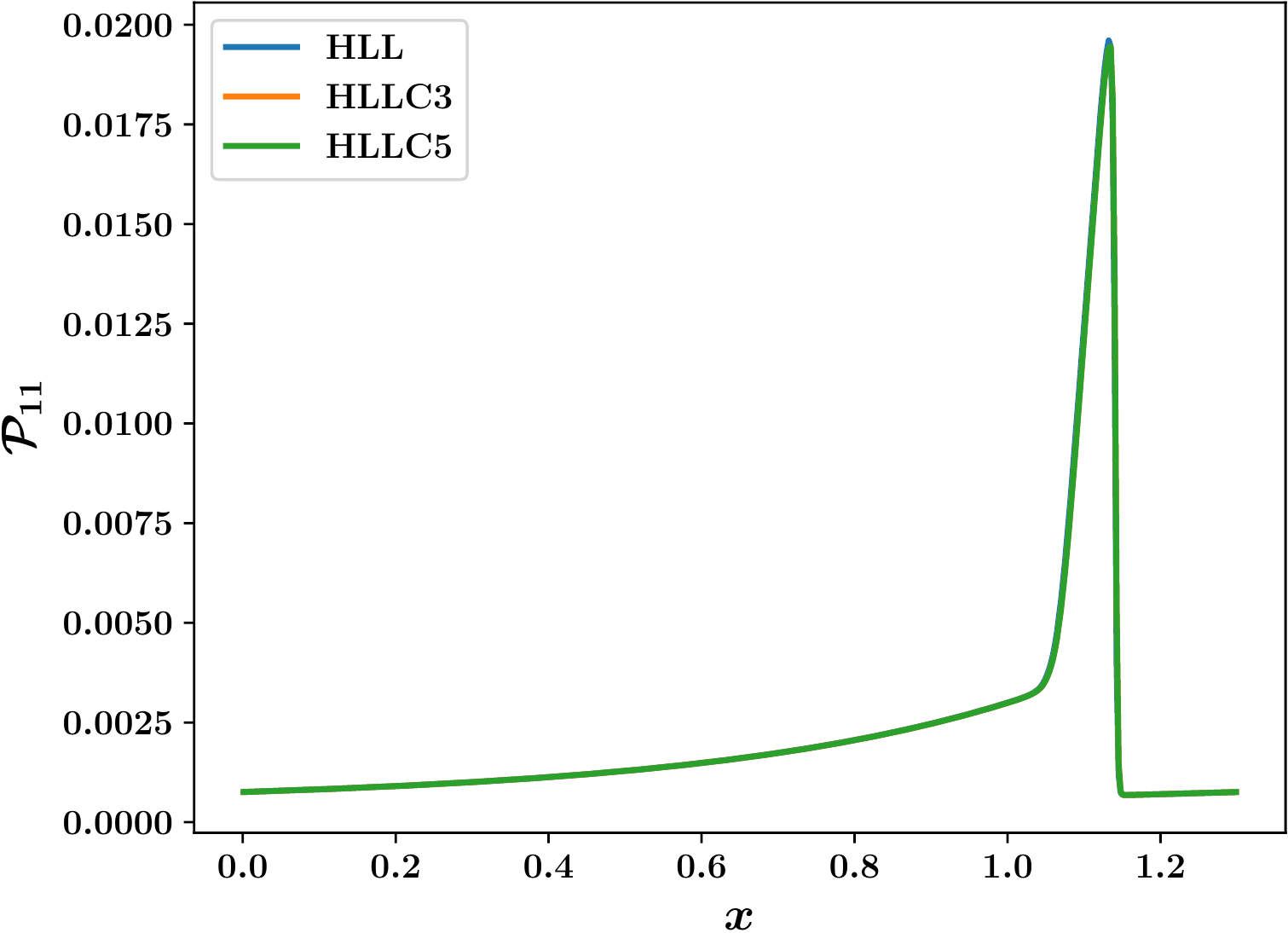} &
			\includegraphics[width=0.48\textwidth]{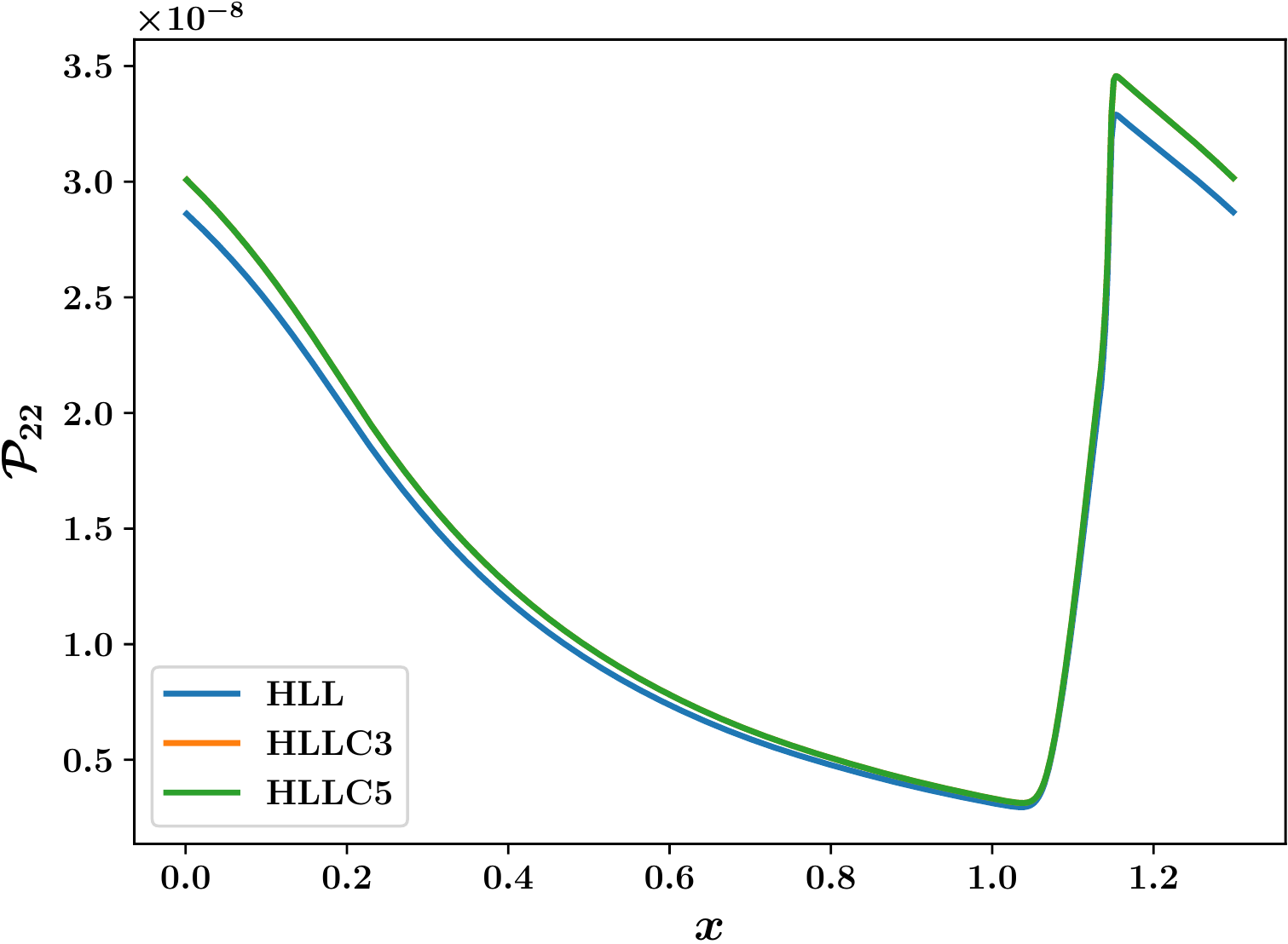} \\
			$\p_{11}$ & $\p_{22}$ 
		\end{tabular}
	\end{center}
	\caption{1-D roll wave problem from Section~(\ref{sec:roll1d}). Plots of water depth $h$, x-velocity $v_1$ and stress components $\p_{11}$ and $\p_{22}$ obtained using second-order scheme with HLL, HLLC3 and HLLC5 solvers on 500 cells at time $t=25$.}
	\label{fig:roll1db}
\end{figure}

\subsection{2-D analytical solution}
The rate of convergence can be estimated from a 2-D exact solution developed in~\cite{Gavrilyuk2018} which is linear in space and non-linear in time. The exact solution is given by
\[
h = \frac{h_0}{1 + \beta^2 t^2}, \quad \vel = \frac{\beta}{1 + \beta^2 t^2} \begin{bmatrix}
\beta t x + y \\
-x + \beta t y \end{bmatrix}, \quad \p = \frac{1}{(1 + \beta^2 t^2)^2} \begin{bmatrix}
\lambda + \gamma \beta^2 t^2 & (\lambda - \gamma) \beta t \\
(\lambda-\gamma) \beta t & \gamma + \lambda \beta^2 t^2 \end{bmatrix}
\]
where the parameters are taken to be $h_0 = 1 \ m$, $\lambda = 0.1 \ m^2/s^2$, $\gamma = 0.01 \ m^2/s^2$, $\beta = 10^{-3} /s$. Previous works \cite{Gavrilyuk2018}, \cite{Bhole2019} have shown first order convergence of the error norm for this test case. We perform the computations on the domain $[0,10] \times [0,10]$ upto the time $t=50$ seconds using Dirichlet boundary conditions where the solution in ghost cells is set to the exact solution. The HLL solver is run on meshes $10^2, 20^2, 40^2, 80^2, 160^2$ while HLLC3 and HLLC5 solvers are run on meshes $10^2, 20^2, 40^2, 80^2$. The convergence of the error norm is shown in Figures~(\ref{fig:exact}) for all three Riemann solvers, which shows second order convergence of our method. The HLL solver requires finer meshes before we see the second order convergence probably due to its larger numerical dissipation as a consequence of not resolving the linear waves.
\begin{figure}
\begin{center}
\begin{tabular}{cc}
\includegraphics[width=0.48\textwidth]{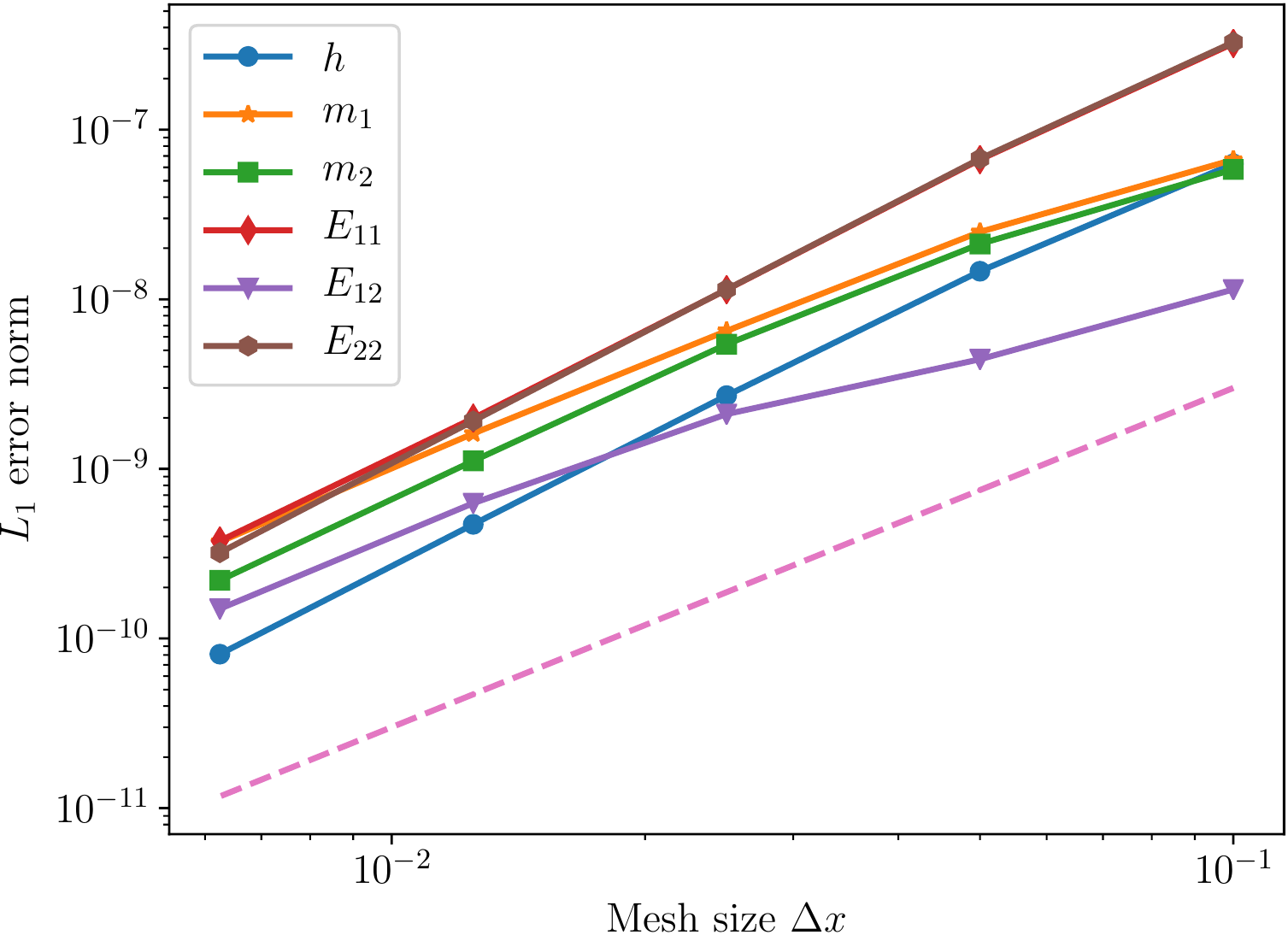} &
\includegraphics[width=0.48\textwidth]{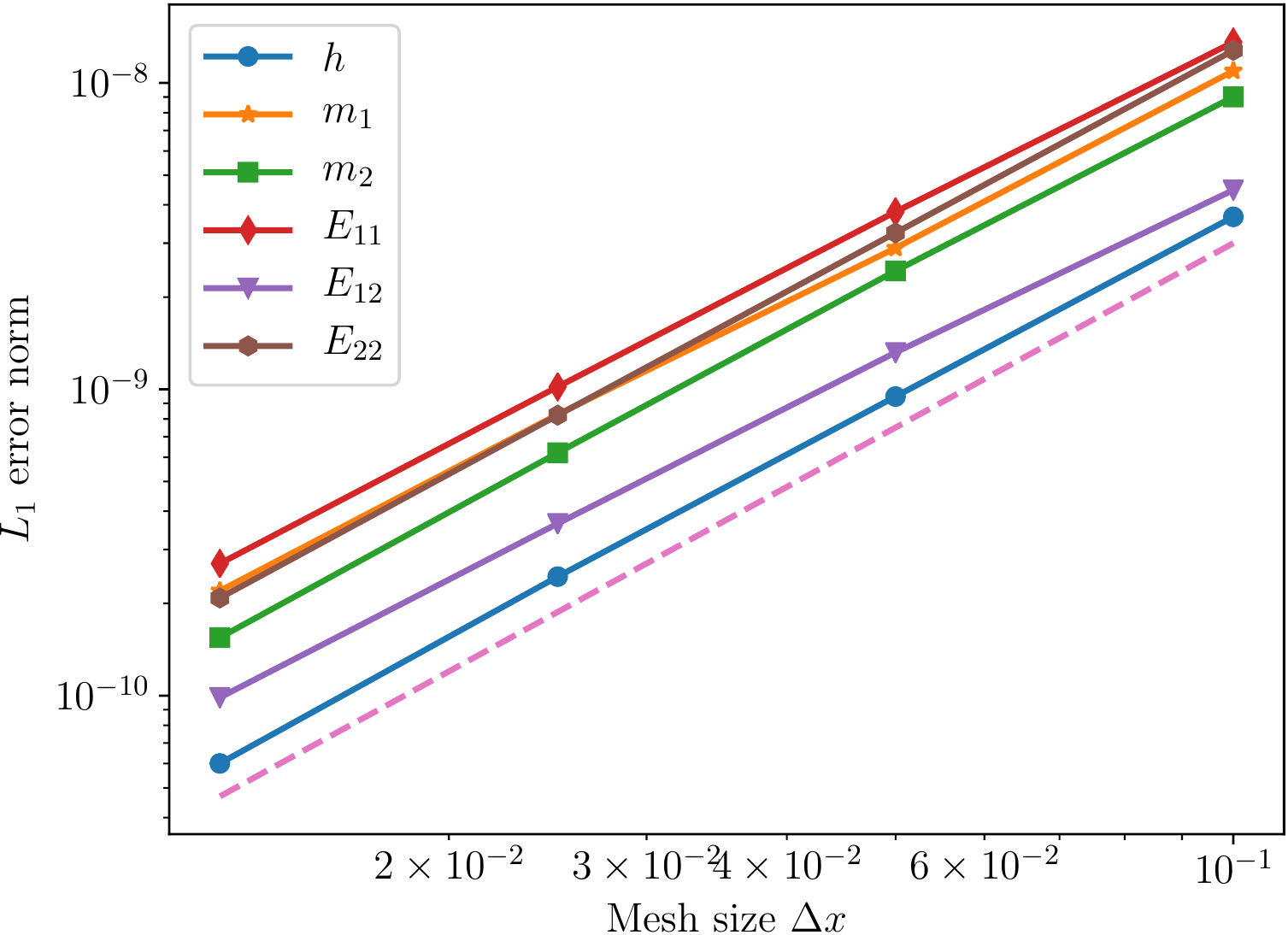} \\
(a) HLL & (b) HLLC3
\end{tabular}
\includegraphics[width=0.48\textwidth]{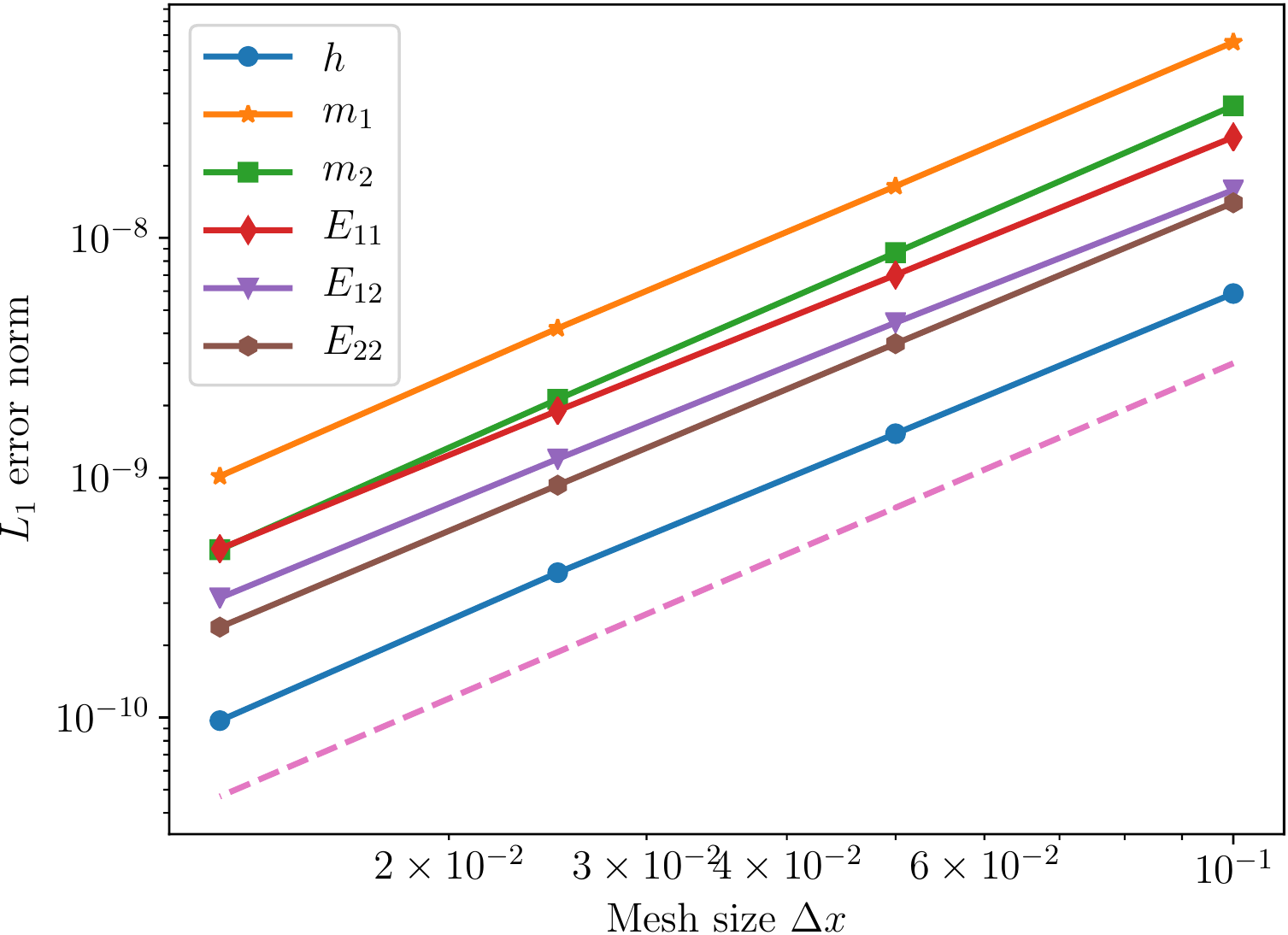} \\
(c) HLLC5
\end{center}
\caption{2-D analytical solution: convergence of error norm at time $t=50$ for all variables with mesh refinement. The broken line shows second order convergence.}
\label{fig:exact}
\end{figure}
\subsection{2-D roll wave problem}
\label{sec:roll2d}
This problem is a 2-D extension of the 1-D roll wave problem described in Section~(\ref{sec:roll1d}) and models the flow of a thin layer of liquid down an inclined plane. The initial conditions are given by
\[
h(x,y,0) = h_0 \left[ 1 + a \sin(2\pi  x/L_x) + a \sin(2\pi  y/L_y) \right], \qquad v_1(x,y,0) = \sqrt{ \frac{g h_0 \tan\theta}{C_f} }, \qquad v_2(x,y,0) = 0
\]
\[
\p_{11}(x,y,0) = \p_{22}(x,y,0) = \half \phi h^2(x,y,0), \qquad \p_{12}(x,y,0) = 0
\]
Here $\theta = 0.05011 \ rad$ is the inclination angle of the bottom surface, $C_f = 0.0036$, $h_0 = 7.98 \times 10^{-3} \ m$, $a = 0.05$, $\phi = 22.76 s^{-2}$, $C_r = 0.00035$, $L_x = 1.3 \ m$, $L_y = 0.5 \ m$. The chosen parameters correspond to those considered in~\cite{Ivanova2017} which leads to the formation of 1-D roll waves starting from a uniform flow which has the same structure as in Brock's experiments~\cite{Brock1969,Brock1970}. A small two dimensional perturbation is added to the water depth whose amplitude is controlled by the coefficient $a$ in the initial condition. The boundary conditions are periodic in both directions. Figure~(\ref{fig:roll2d1}) shows the surface plots of the depth field on a $1040 \times 400$ mesh while Figure~(\ref{fig:roll2d2}) shows similar plots on a finer mesh of $2080 \times 800$ cells. The first order scheme shows rather smooth solutions due to higher numerical dissipation while second order scheme shows more features in the solution. The HLL and HLLC3 solvers also show somewhat smooth solutions while HLLC5 shows a more turbulent-type of solution. However, even the HLL/HLLC3 solvers show transverse wave structures in the second order scheme. The first order results from HLLC5 solver look similar to the corresponding results in the literature~\cite{Gavrilyuk2018},~\cite{Bhole2019}. 

The development of the water depth profile with time using HLLC5 solver can be seen in Figure~(\ref{fig:roll2d3}), which shows three distinct phases. In the first phase, the profile develops a 1-D structure similar to the 1-D computation in Section~(\ref{sec:roll1d}) and is similar to Brock's experimental profile. At around time $t=10$, the profile develops some perturbations at the hydraulic jump which generates transverse structures and subsequently this spreads to the whole domain. In Figure~(\ref{fig:roll2d4}), we compute the $y$-average of the solution and also show the fluctuations around the $y$-average at time $t=36$. The average profile still resembles the 1-D profile with fluctuations superimposed on top of it. These fluctuations are more prominent in the HLLC5 solver which is the most sophisticated solver as it includes all the waves in the Riemann solution. Figure~(\ref{fig:roll2d5}) shows a zoomed view in the region $[0,0.5] \times [0,0.5]$ \revb{where we show the fluctuations of the height field relative to the $y$-average}; the HLLC3 solver yields somewhat regular wave patterns while the HLLC5 solver shows a more dis-organized behaviour in the solution. 

\revb{We compute the spectrum $E(k)$ of the kinetic energy of velocity fluctuations by first computing the two dimensional discrete Fourier transform of each velocity component to obtain $\hat{u}(k_x,k_y)$, $\hat{v}(k_x,k_y)$, and then integrating $E(k_x,k_y) = \half ( |\hat{u}|^2 + |\hat{v}|^2)$ over the circle of radius $k = \sqrt{k_x^2 + k_y^2}$. The DFT is computed using the Fast Fourier Transform method available in Numpy. The fluctuation energy distribution over the wave numbers is shown in Figure~(\ref{fig:tkespec}).} Firstly, we see that the spectra at times $t=36$ and $t=60$ are very similar which indicates that a statistically steady state has been reached. All the three solvers show a $k^{-4}$ spectrum in intermediate wave-numbers. This is in contrast to two dimensional hydrodynamic turbulence which exhibits a $k^{-3}$ energy spectrum~\cite{Boffetta2012}.  We also show the $k^{-5/3}$ line which is observed in three dimensional turbulence. With the present computations, we are not able to definitely conclude on the structure of the energy spectrum in the small and intermediate wave-numbers and a more refined computation with a higher order scheme may be necessary to clearly identify the precise scaling law in those wave-numbers. These results show that a higher order scheme is essential to capture the turbulent-like structure in the solution. Also, there seems to be a critical dependence on the wave structure included in the Riemann solver and the type of solutions observed in the computations. Only the five wave solver that models all the waves and combined with a higher order scheme yields solutions which seem to correspond better with observations in experiments~\cite{Aranda2016}.
\begin{figure}
\begin{center}
\begin{tabular}{cc}
\includegraphics[width=0.35\textwidth]{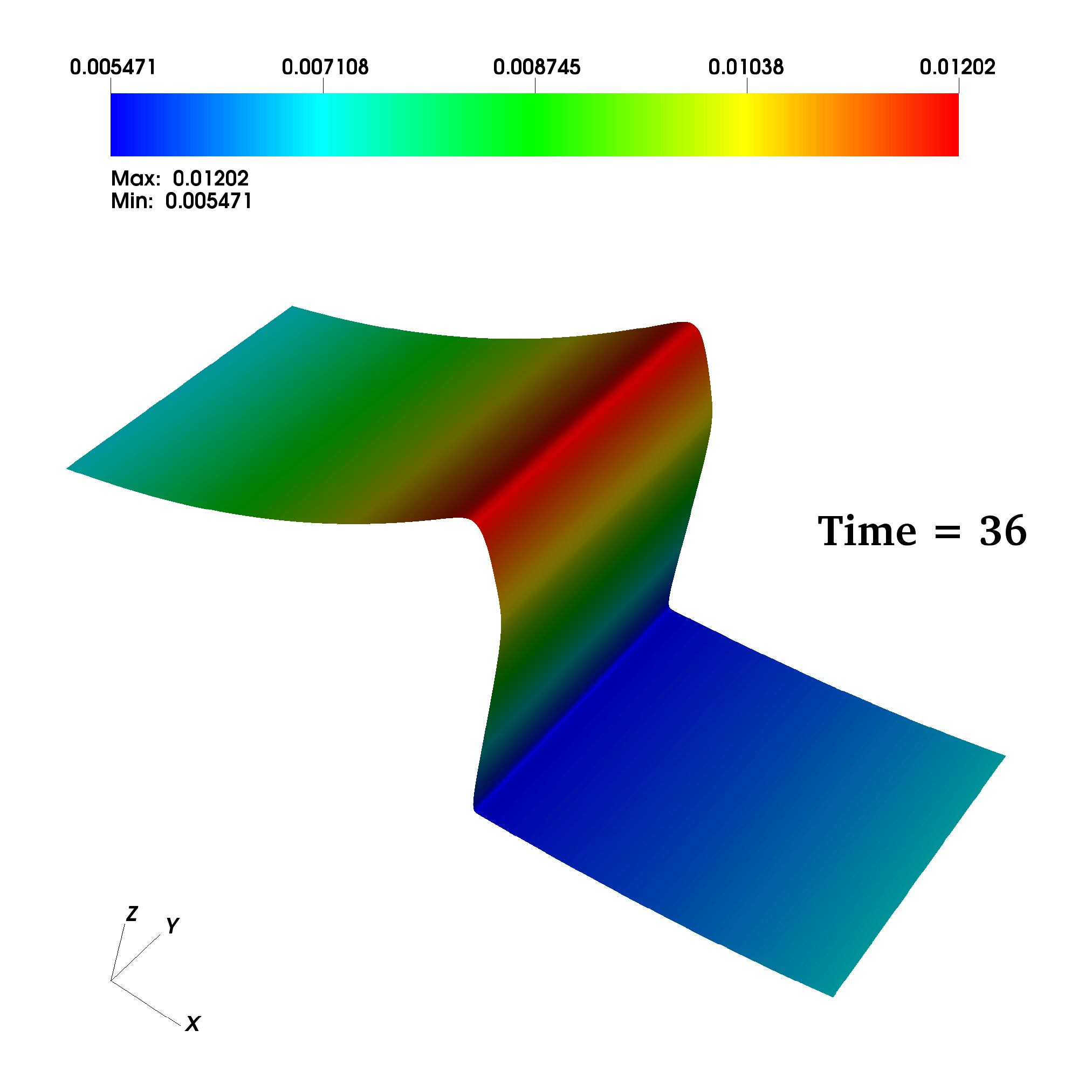} &
\includegraphics[width=0.35\textwidth]{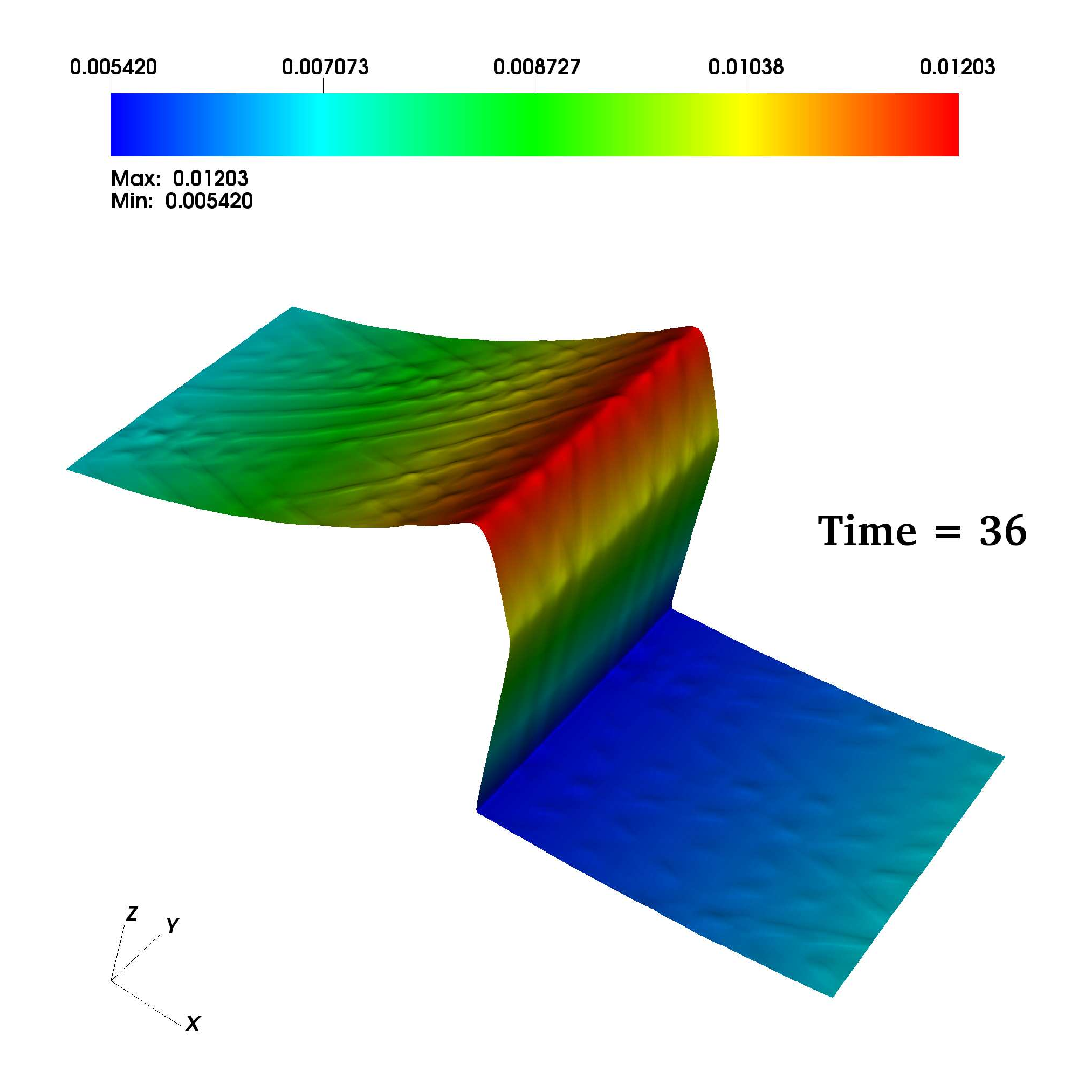} \\
HLL, first order & HLL, second order \\
\includegraphics[width=0.35\textwidth]{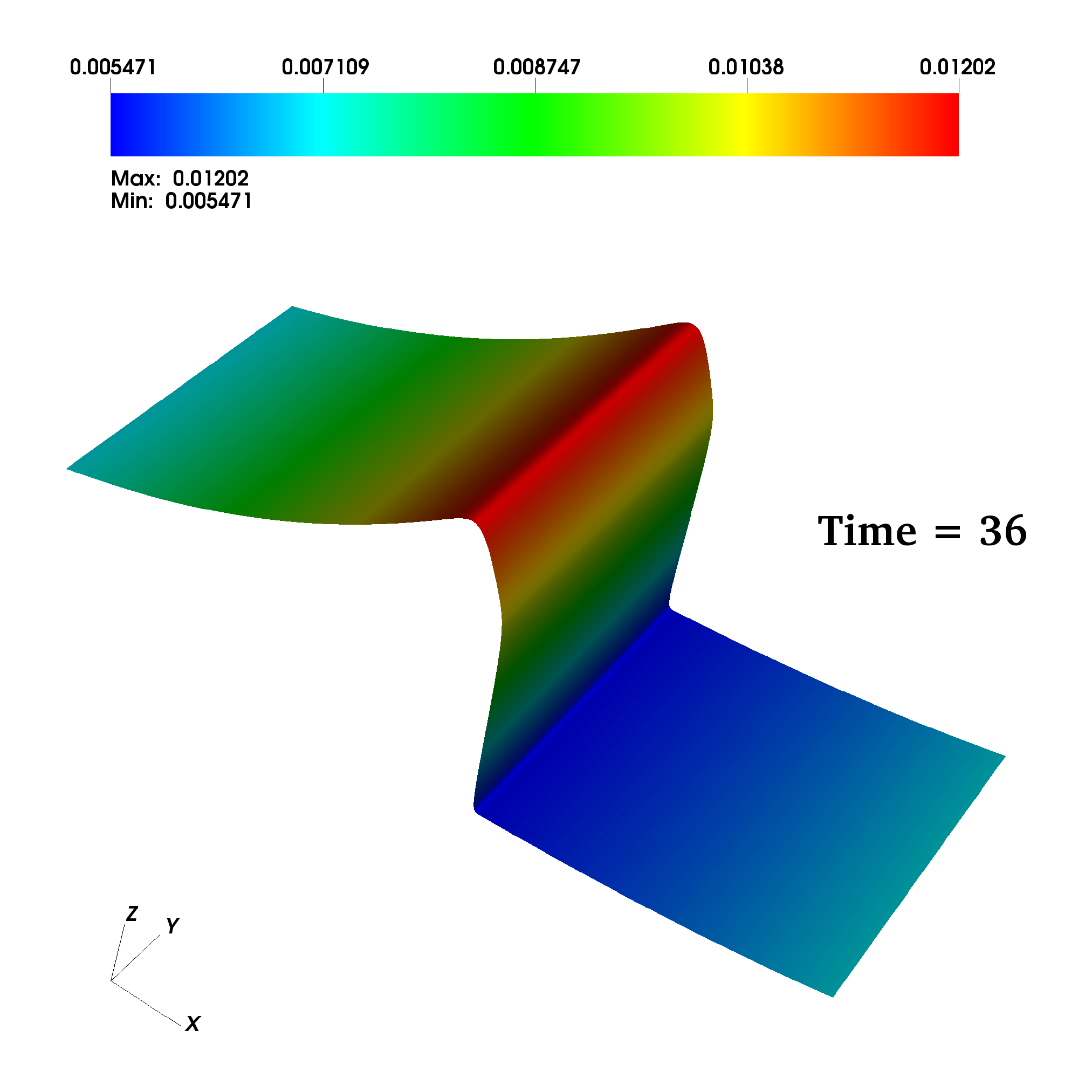} &
\includegraphics[width=0.35\textwidth]{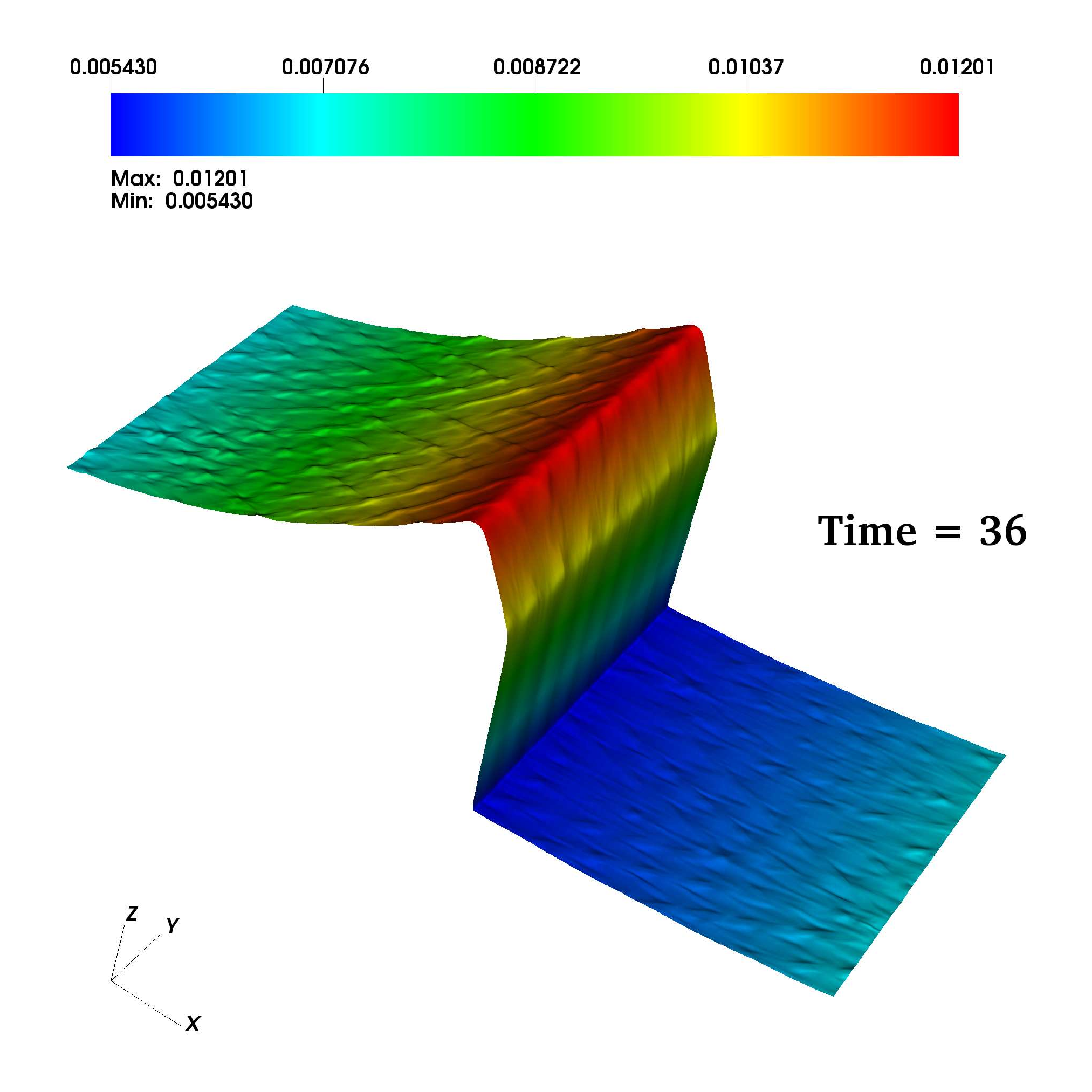} \\
HLLC3, first order & HLLC3, second order \\
\includegraphics[width=0.35\textwidth]{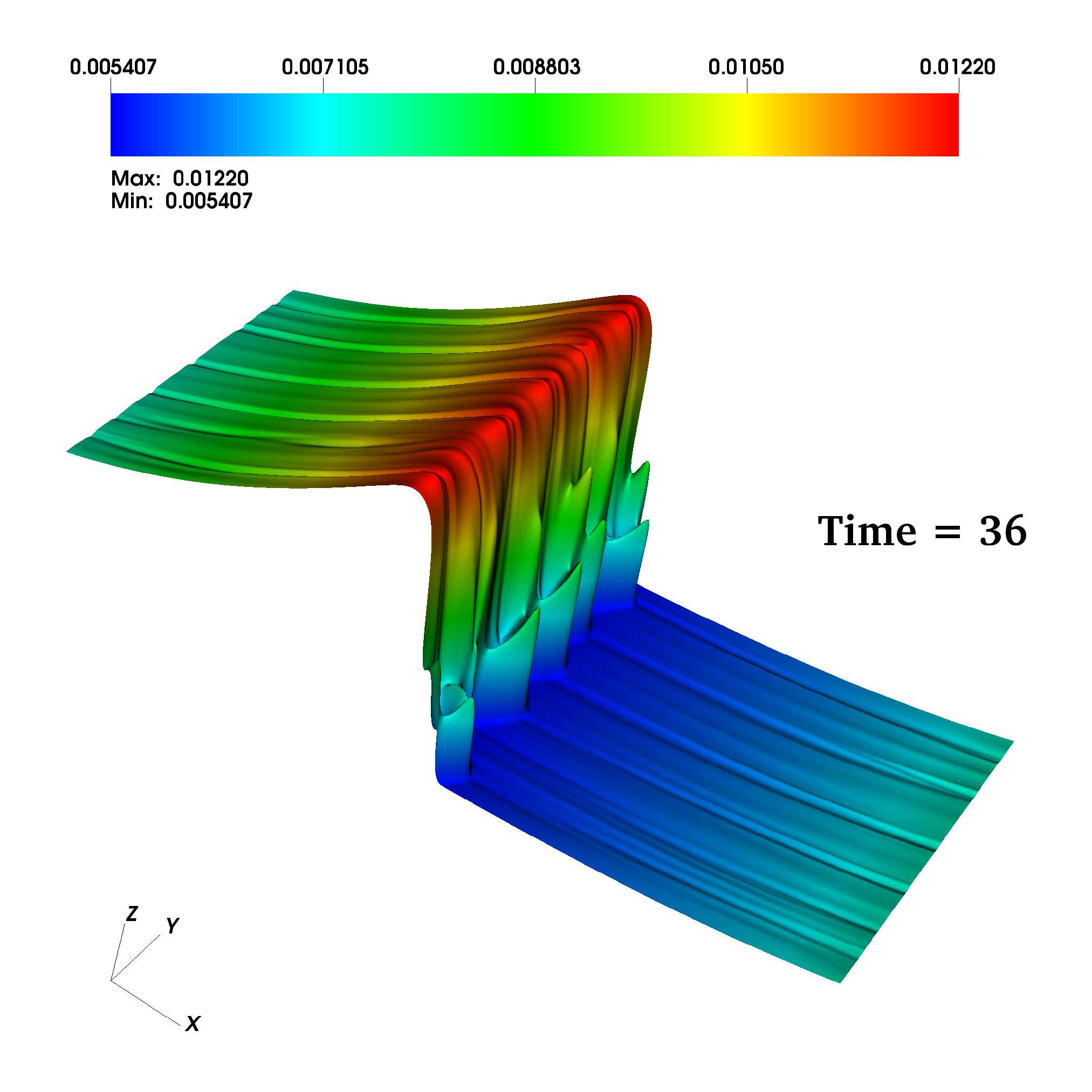} &
\includegraphics[width=0.35\textwidth]{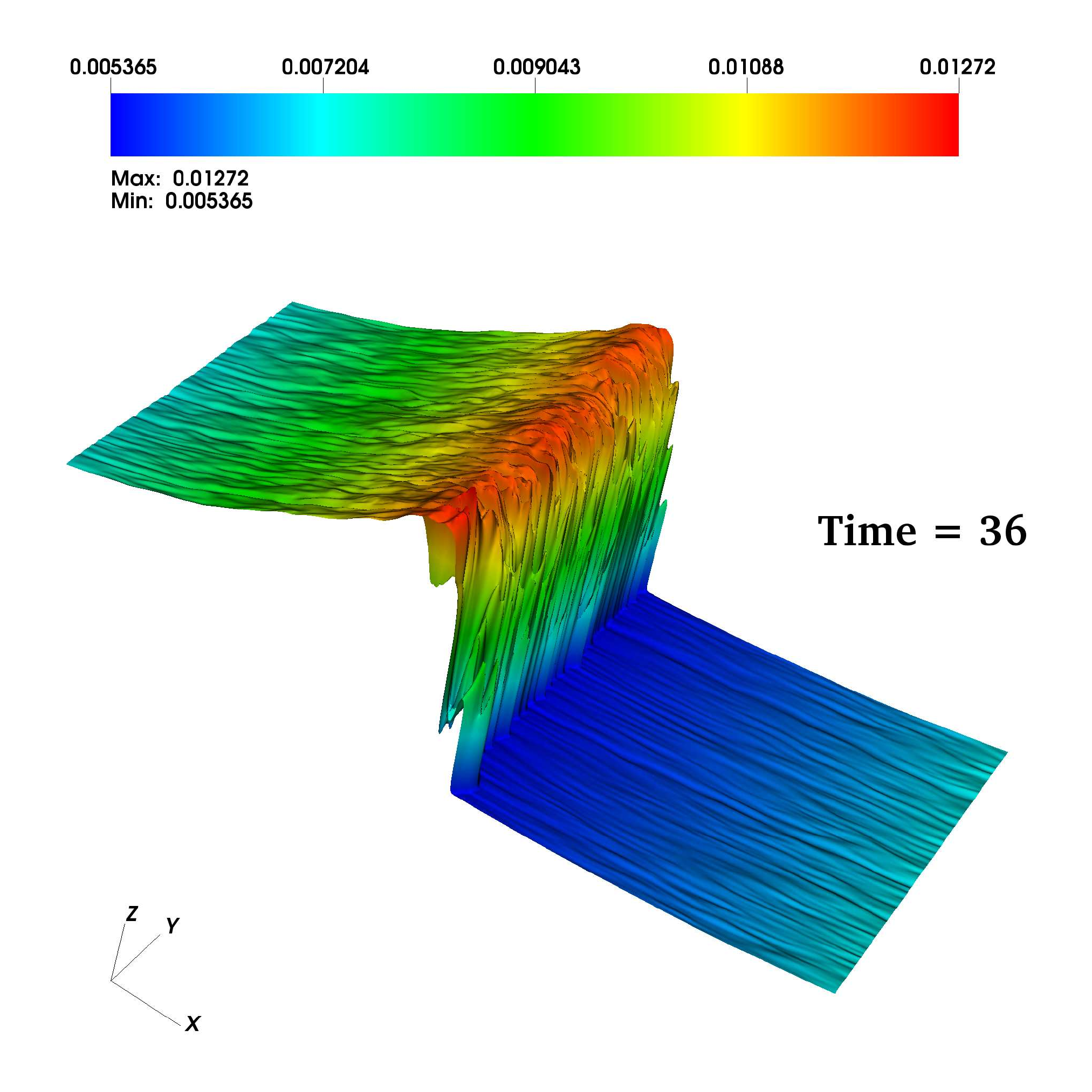} \\
HLLC5, first order & HLLC5, second order \\
\end{tabular}
\end{center}
\caption{2-D roll wave problem using $1040 \times 400$ mesh. Depth field at time $t=36$ units.}
\label{fig:roll2d1}
\end{figure}

\begin{figure}
\begin{center}
\begin{tabular}{cc}
\includegraphics[width=0.35\textwidth]{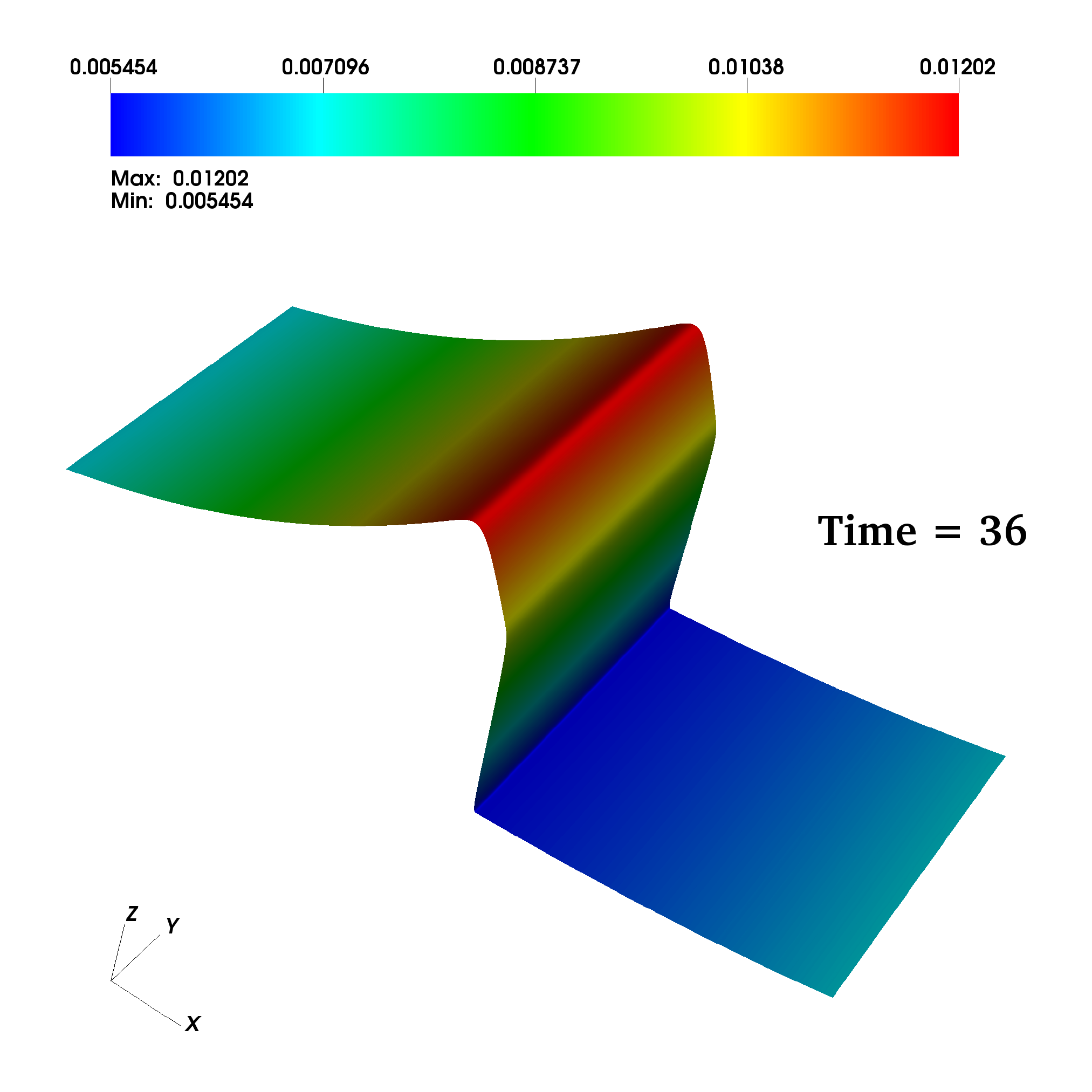} &
\includegraphics[width=0.35\textwidth]{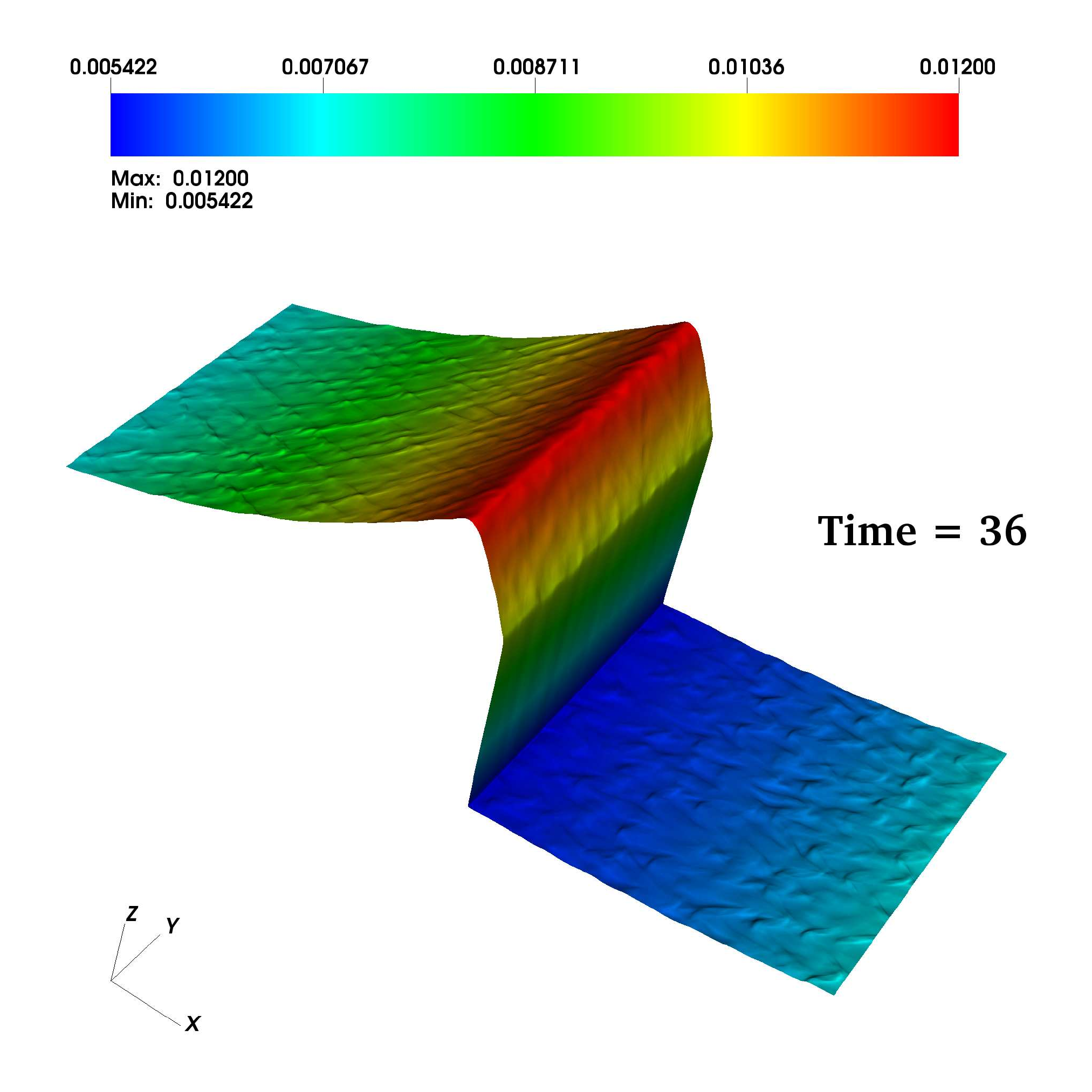} \\
HLL, first order & HLL, second order \\
\includegraphics[width=0.35\textwidth]{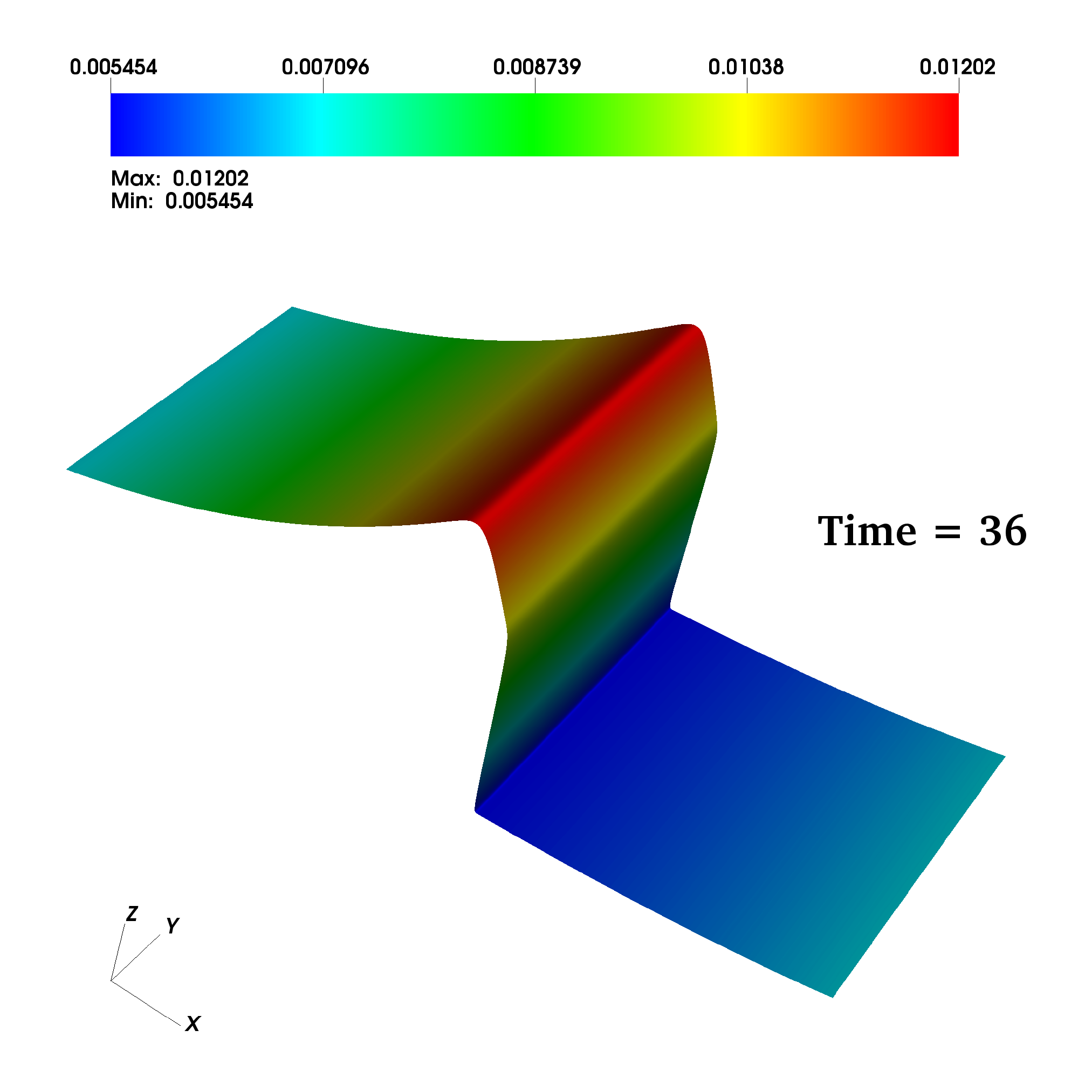} &
\includegraphics[width=0.35\textwidth]{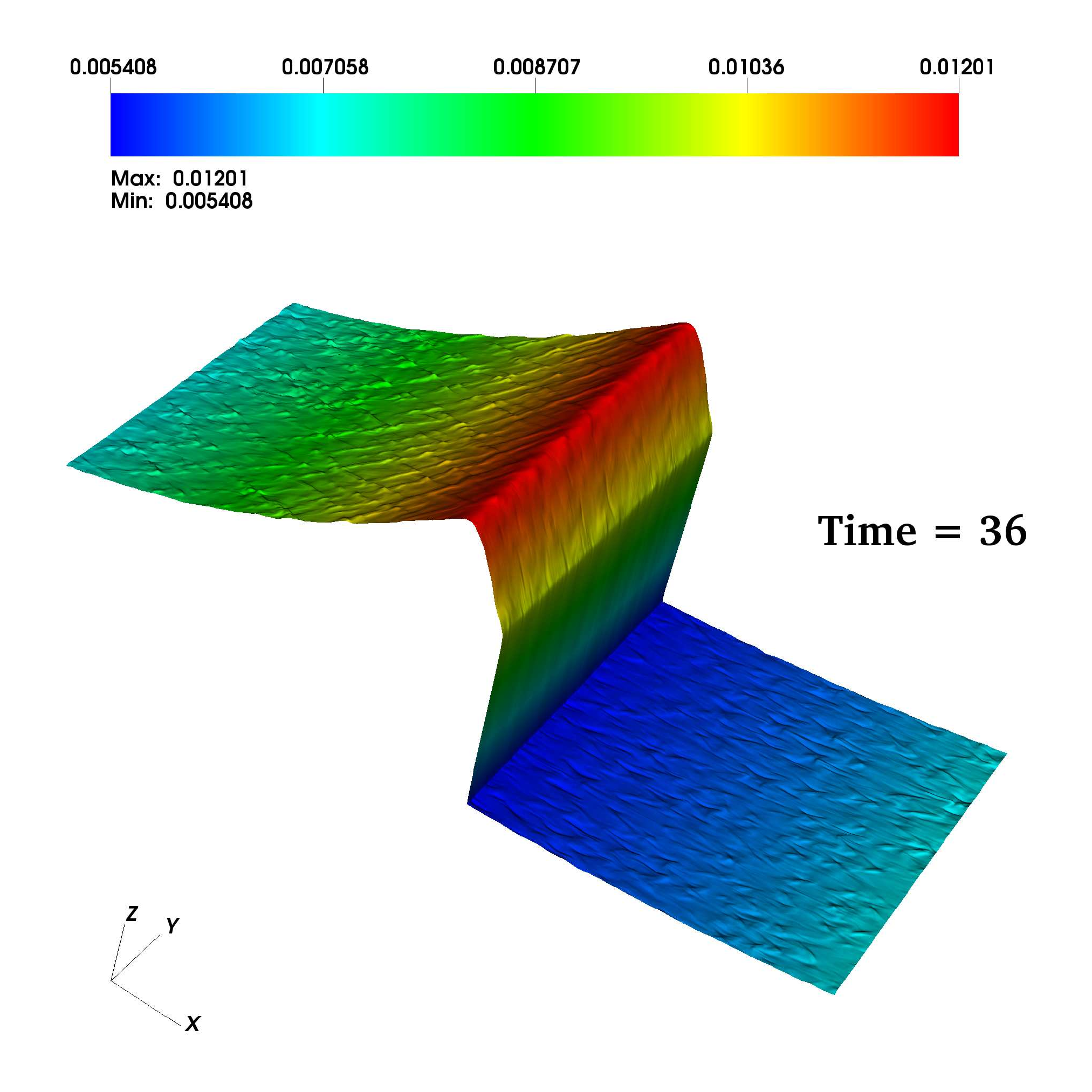} \\
HLLC3, first order & HLLC3, second order \\
\includegraphics[width=0.35\textwidth]{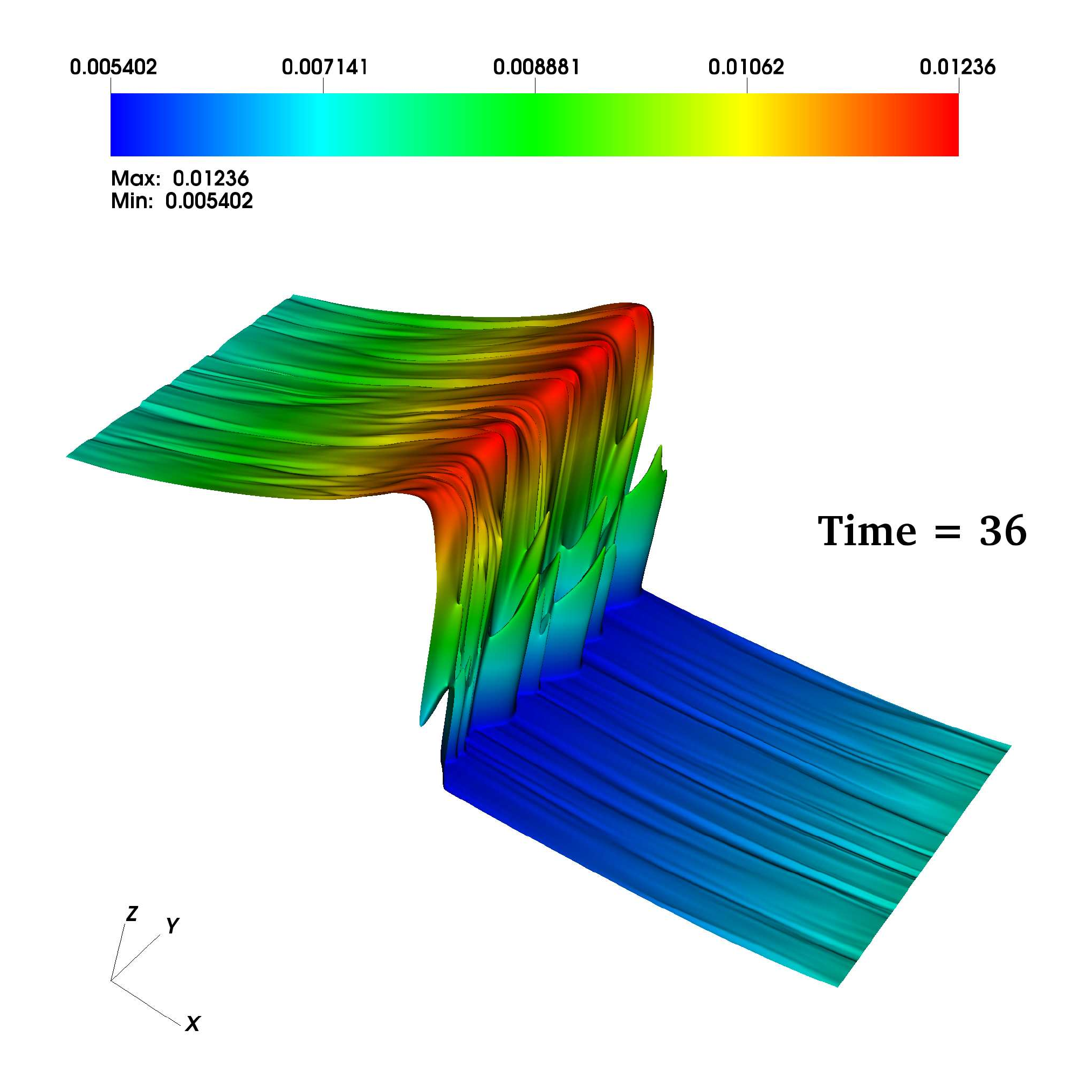} &
\includegraphics[width=0.35\textwidth]{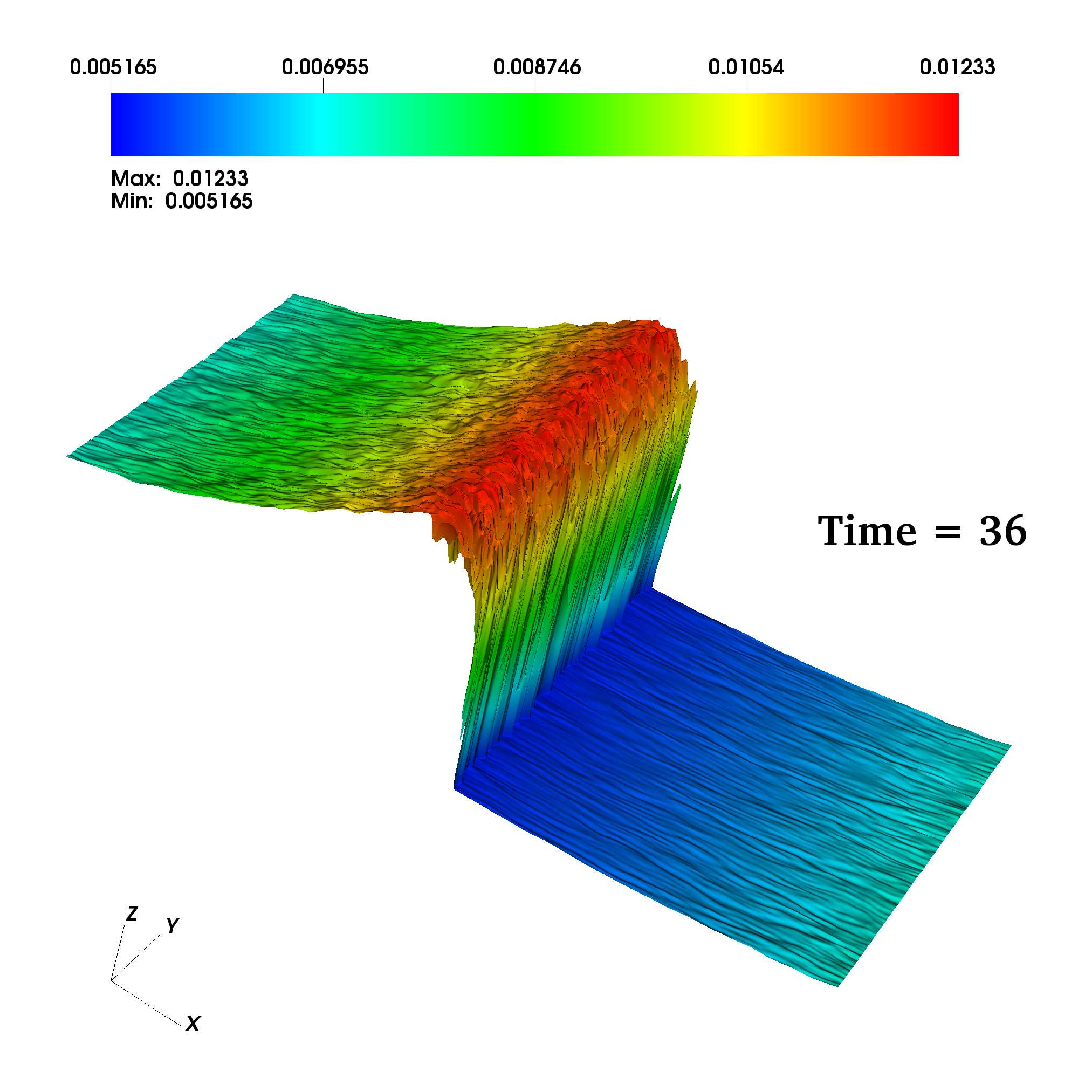} \\
HLLC5, first order & HLLC5, second order \\
\end{tabular}
\end{center}
\caption{2-D roll wave problem using $2080 \times 800$ mesh. Depth field at time $t=36$ units.}
\label{fig:roll2d2}
\end{figure}

\begin{figure}
\begin{center}
\begin{tabular}{cc}
\includegraphics[width=0.48\textwidth]{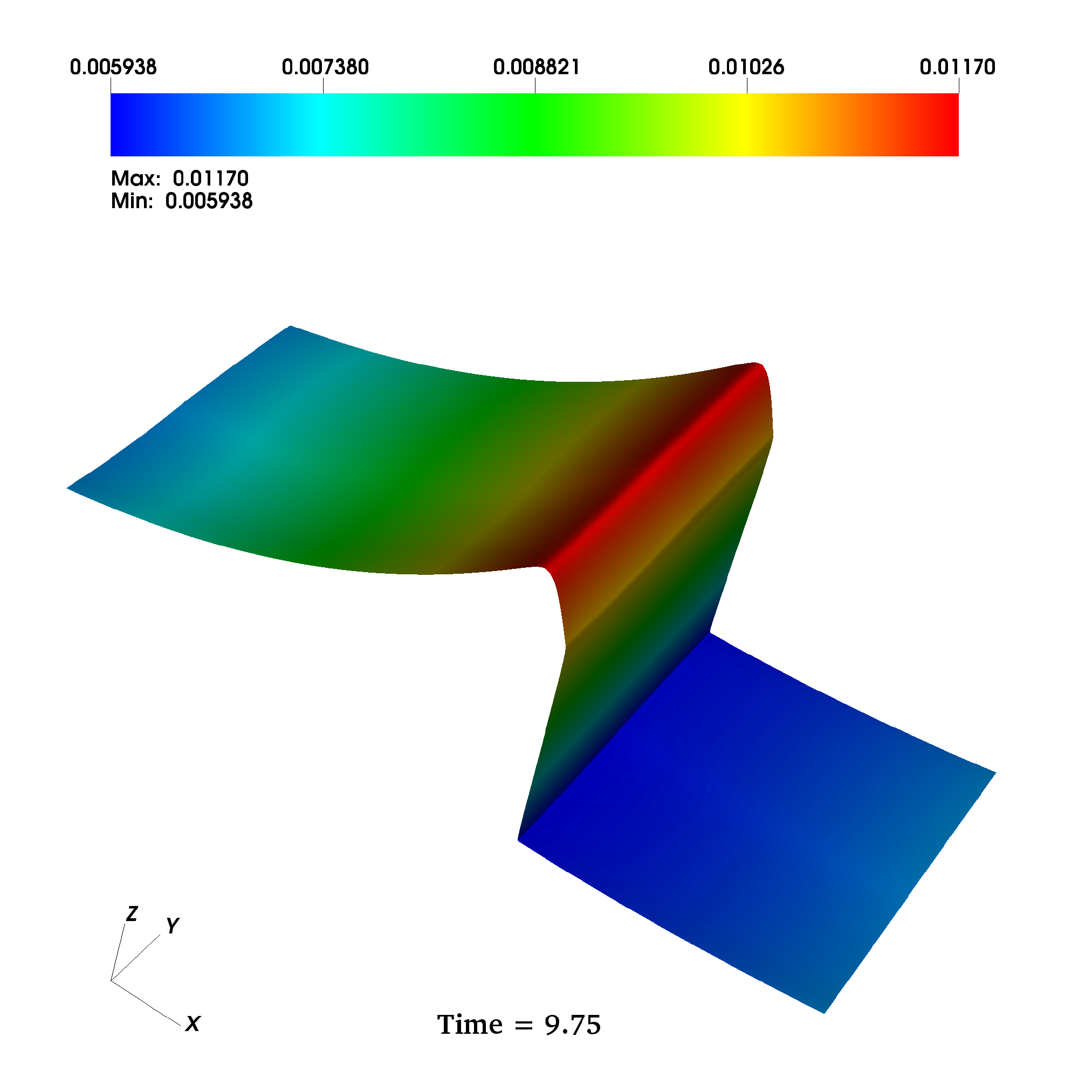} &
\includegraphics[width=0.48\textwidth]{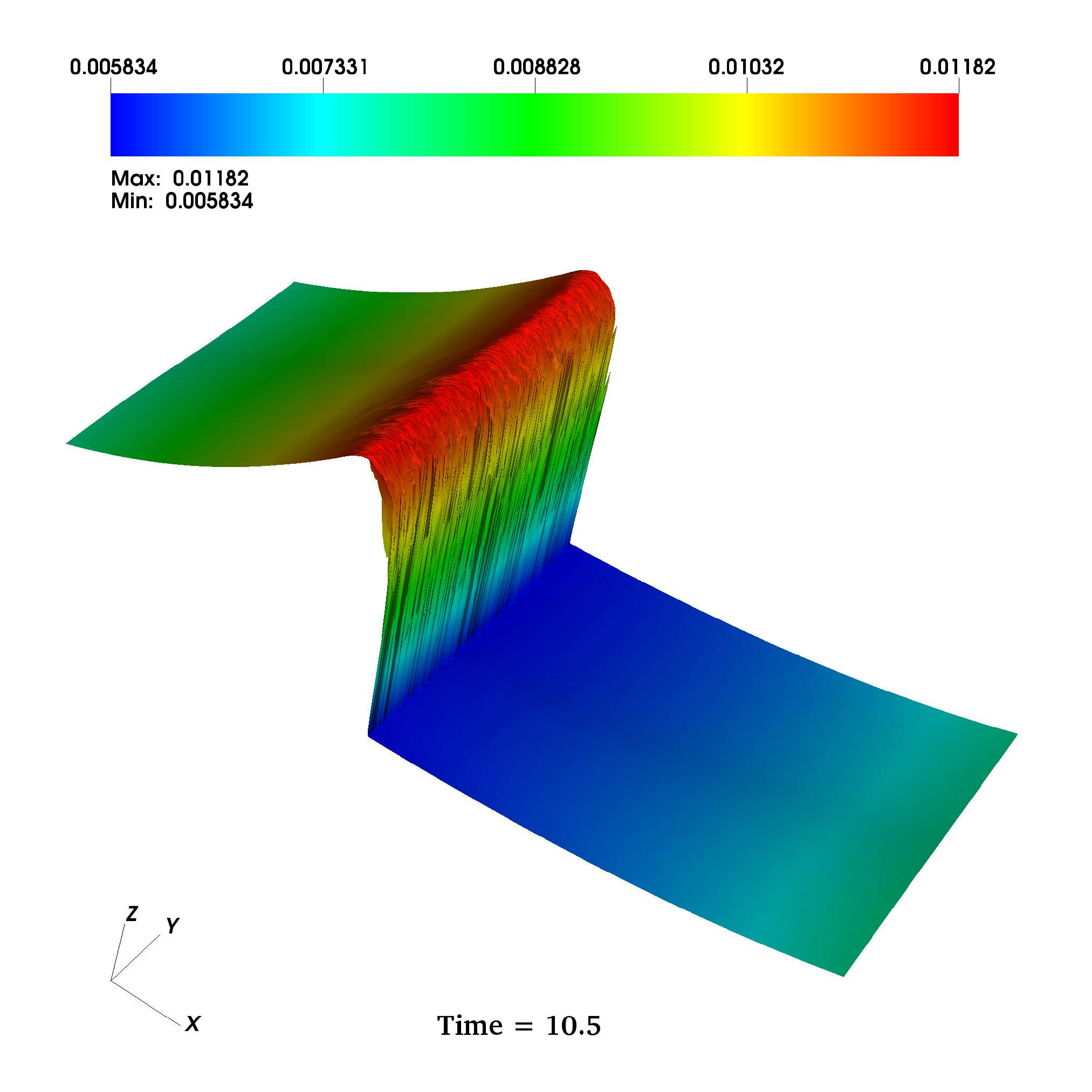} \\
\includegraphics[width=0.48\textwidth]{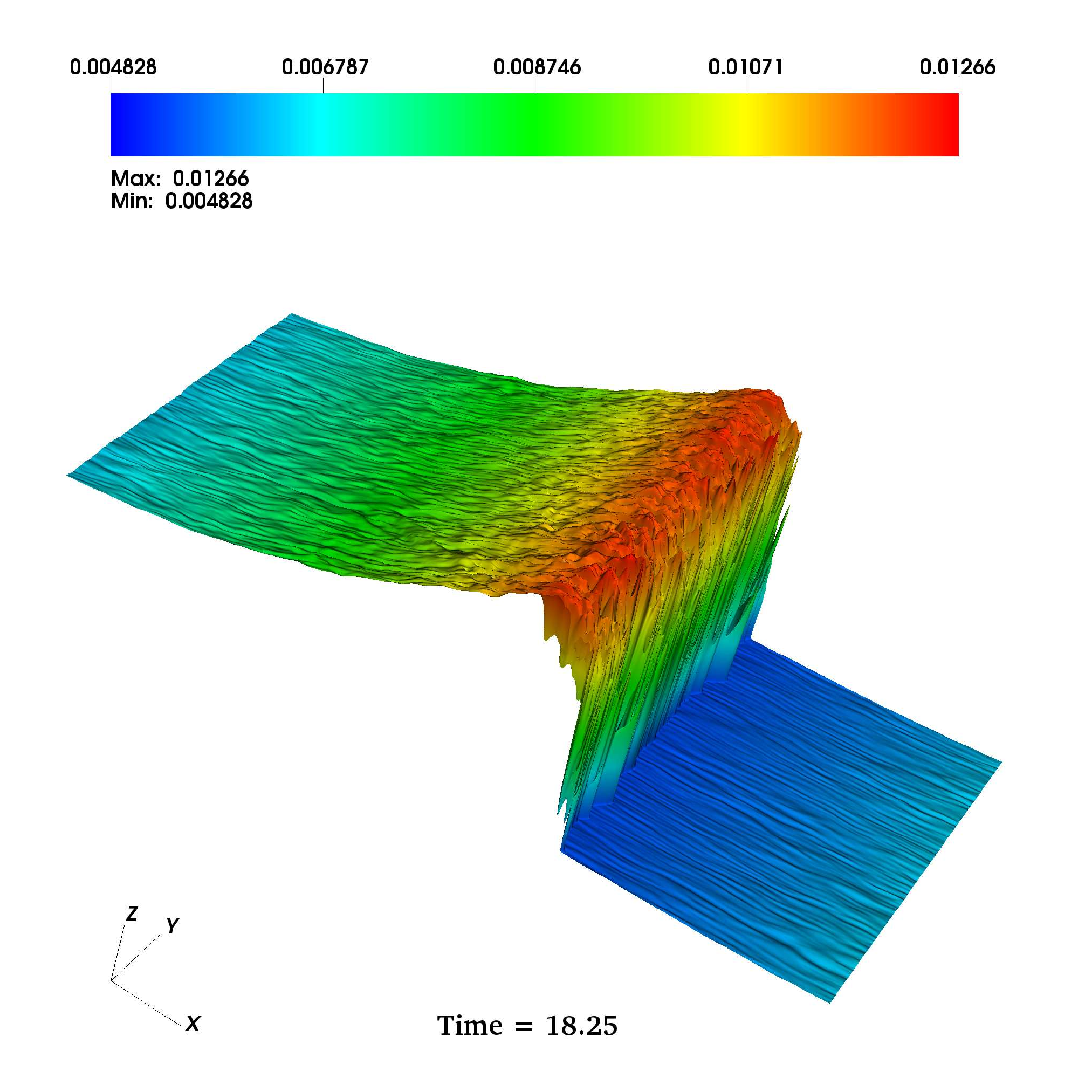} &
\includegraphics[width=0.48\textwidth]{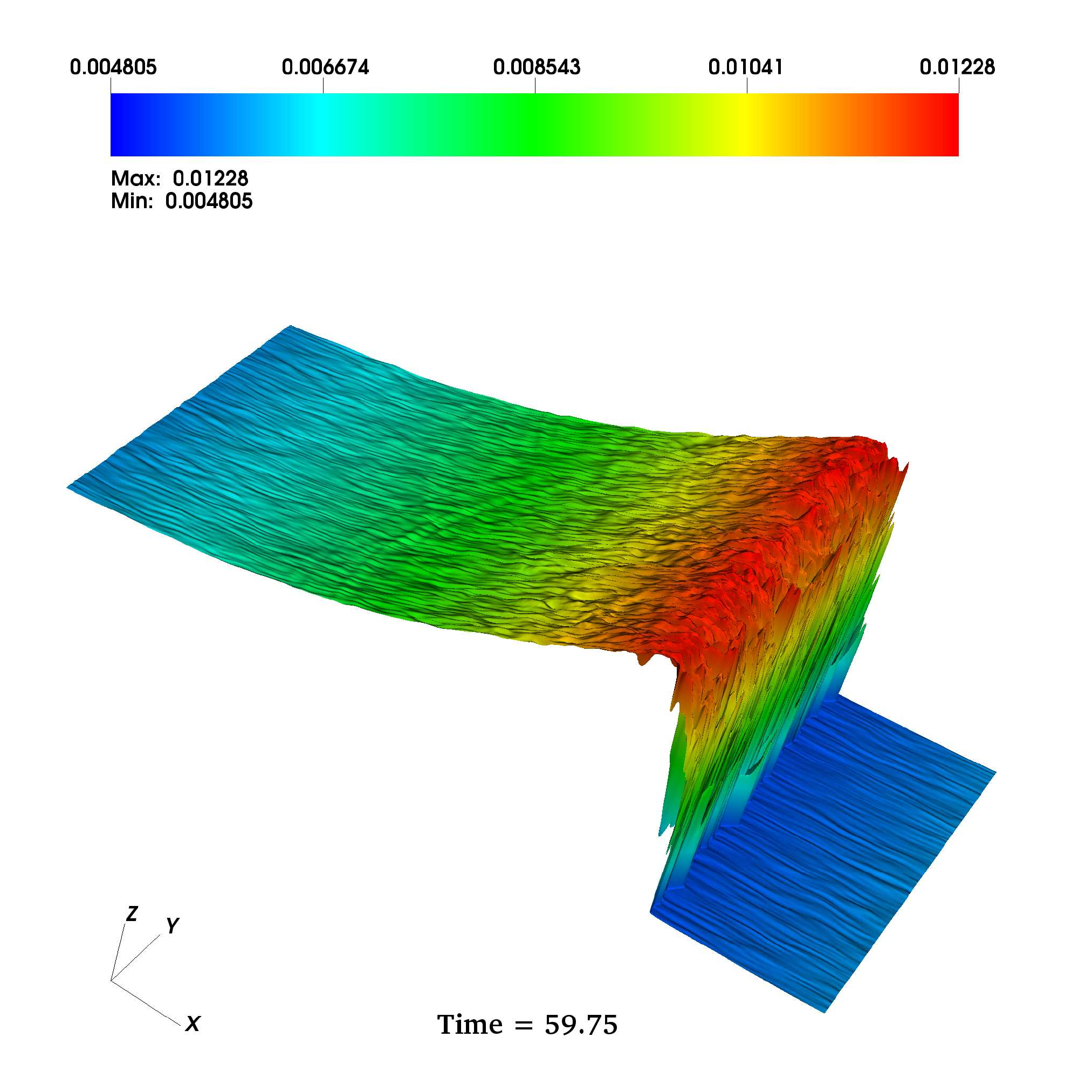}
\end{tabular}
\end{center}
\caption{2-D roll wave problem using second order HLLC5 scheme on $2080 \times 800$ mesh. Depth field at  different times.}
\label{fig:roll2d3}
\end{figure}

\begin{figure}
\begin{center}
\begin{tabular}{cc}
\includegraphics[width=0.48\textwidth]{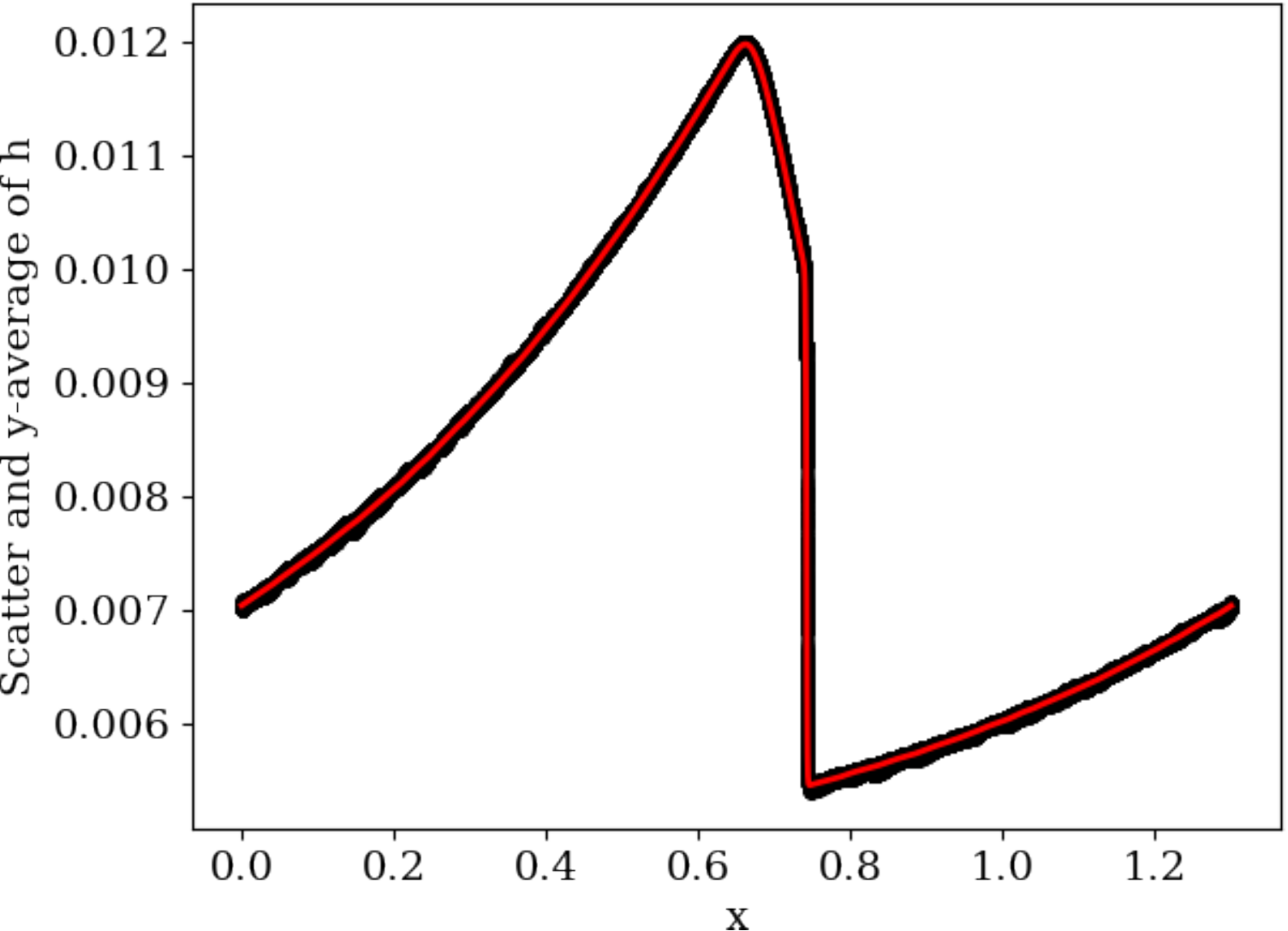} &
\includegraphics[width=0.48\textwidth]{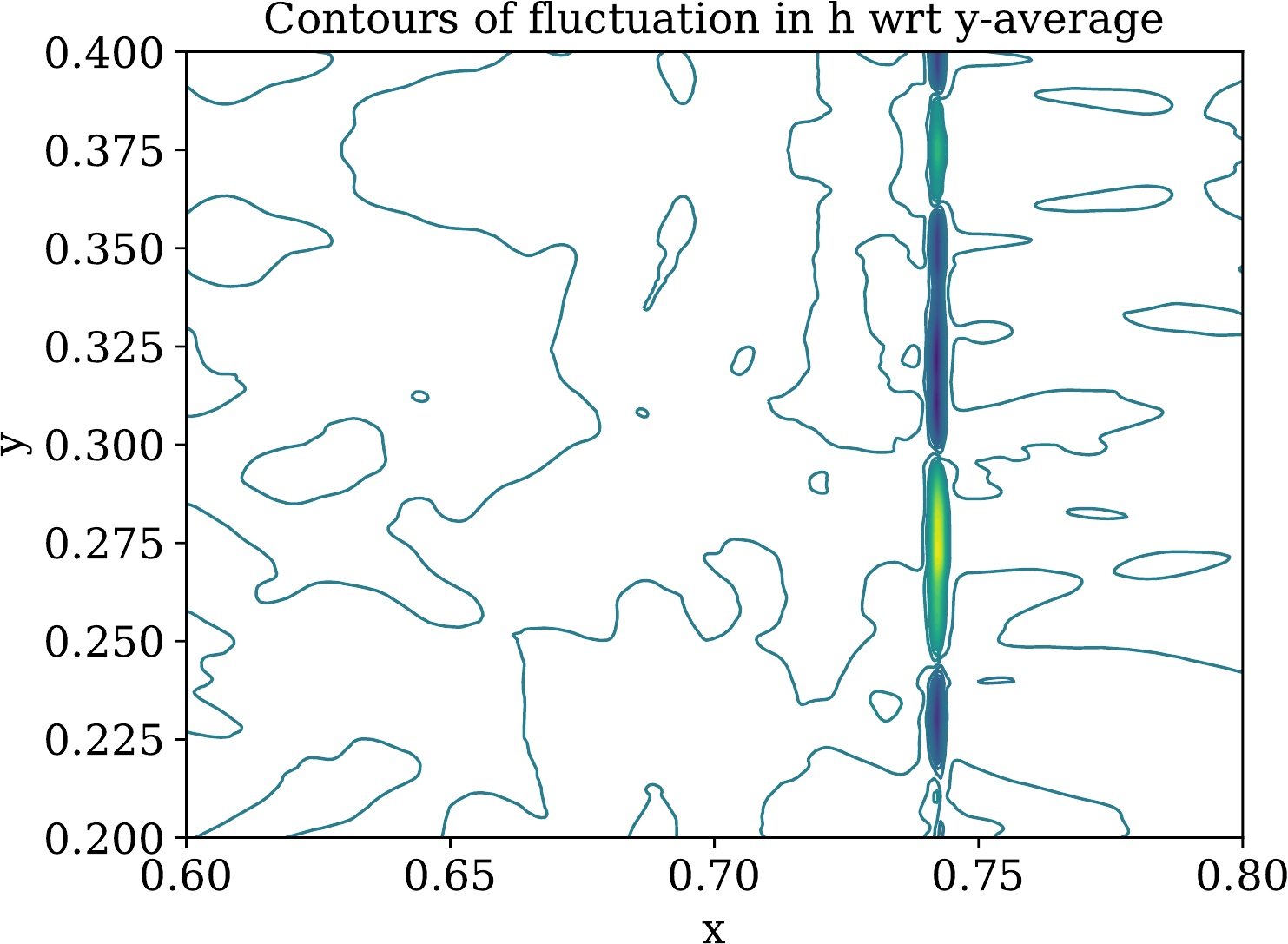} \\
\includegraphics[width=0.48\textwidth]{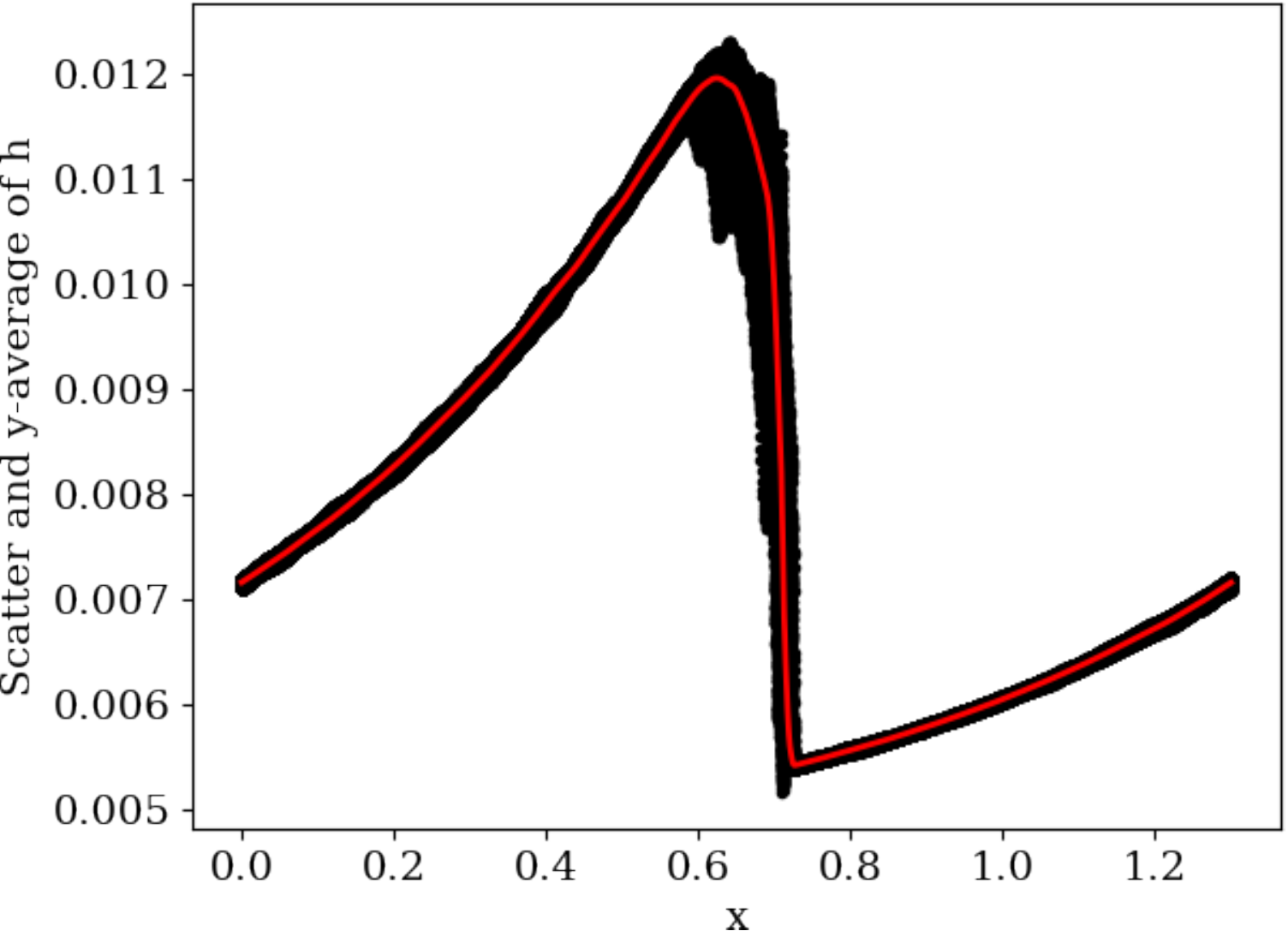} &
\includegraphics[width=0.48\textwidth]{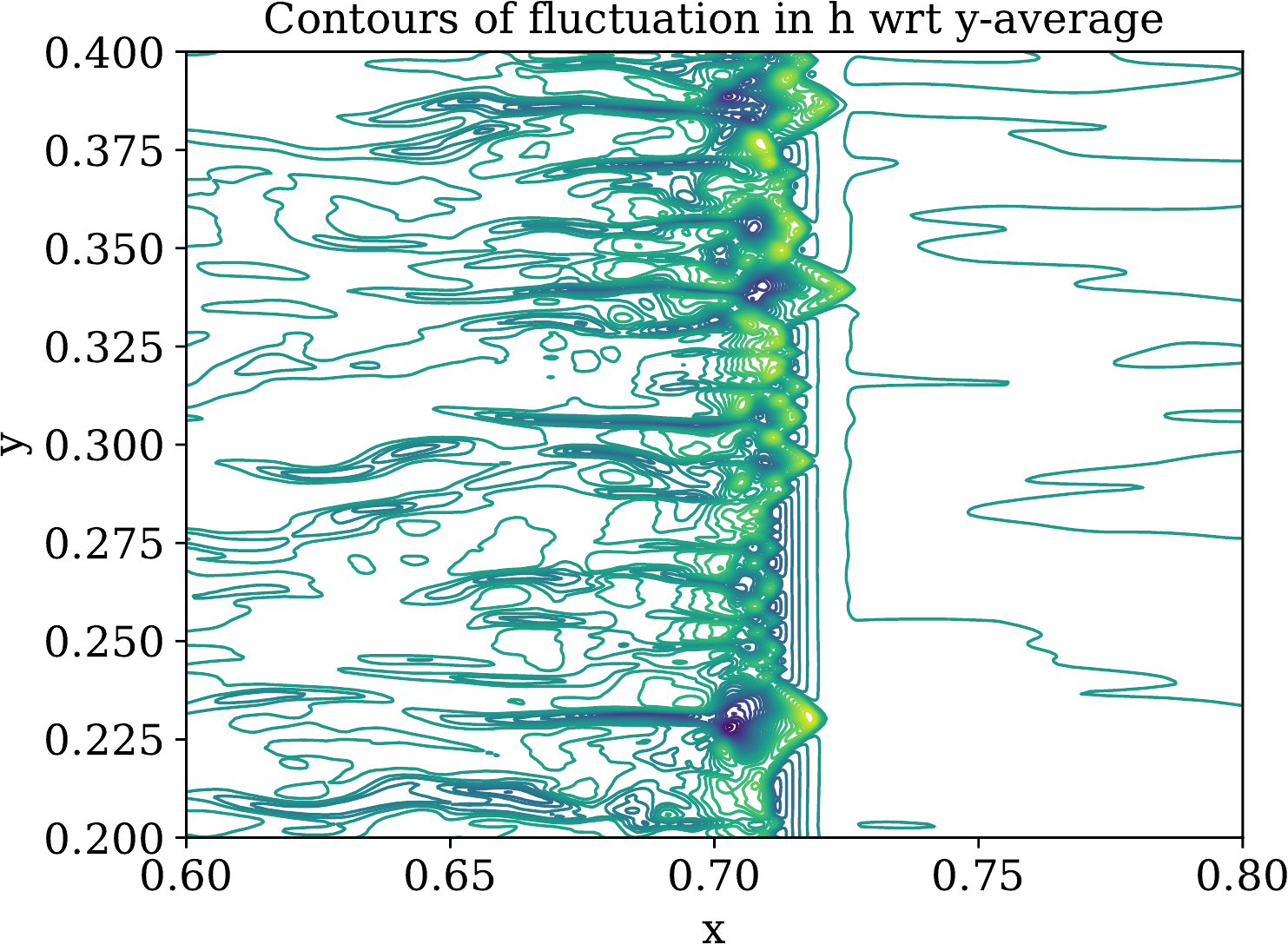}
\end{tabular}
\end{center}
\caption{2-D roll wave problem using second order scheme on $2080 \times 800$ mesh at $t=36$. Top: HLLC3, Bottom: HLLC5.  Left: Scatter plot of depth field and its $y$-average. Right: Zoomed view of fluctuations of depth field relative to $y$-average.}
\label{fig:roll2d4}
\end{figure}

\begin{figure}
\begin{center}
\begin{tabular}{cc}
\includegraphics[width=0.48\textwidth]{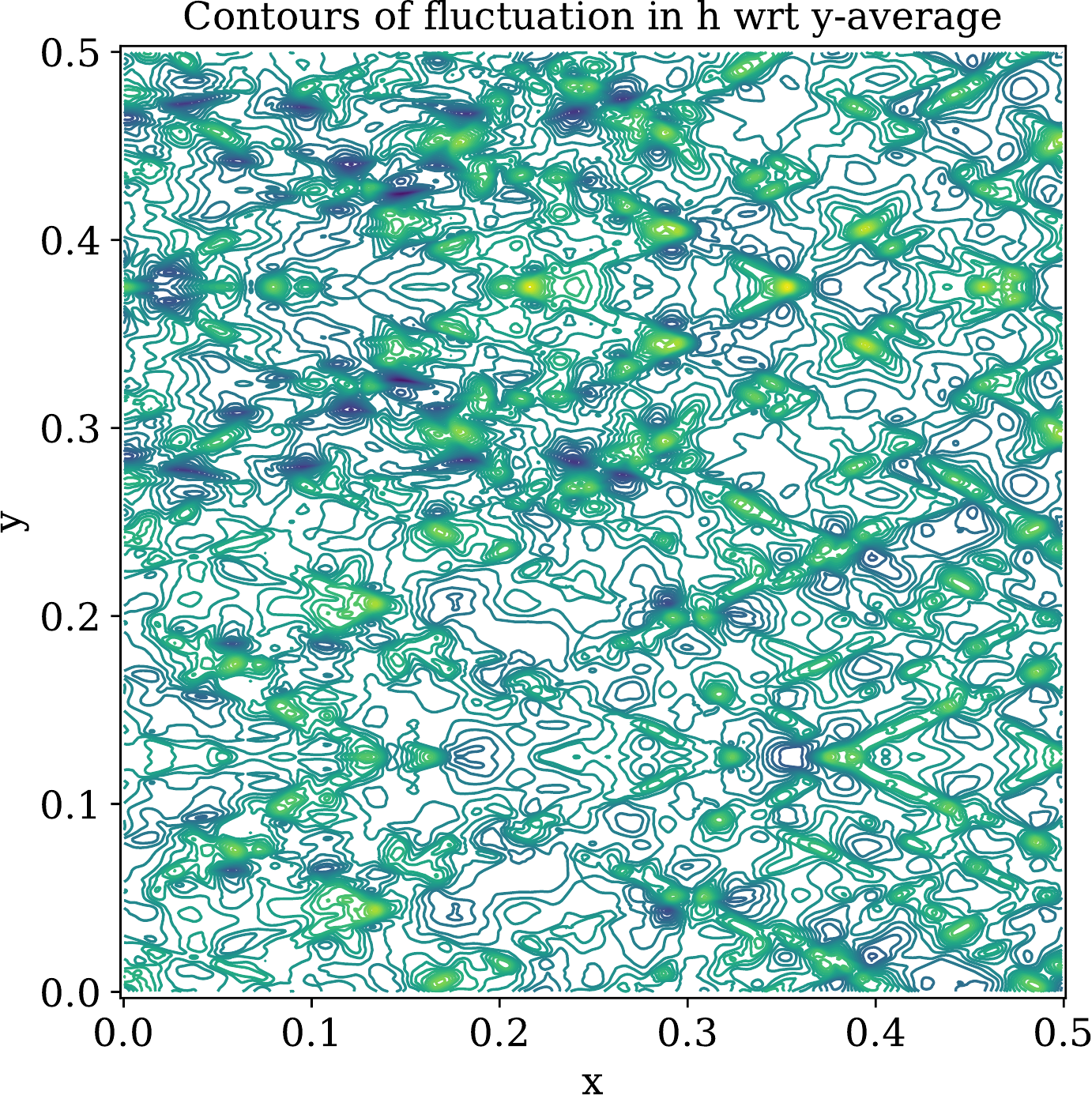} &
\includegraphics[width=0.48\textwidth]{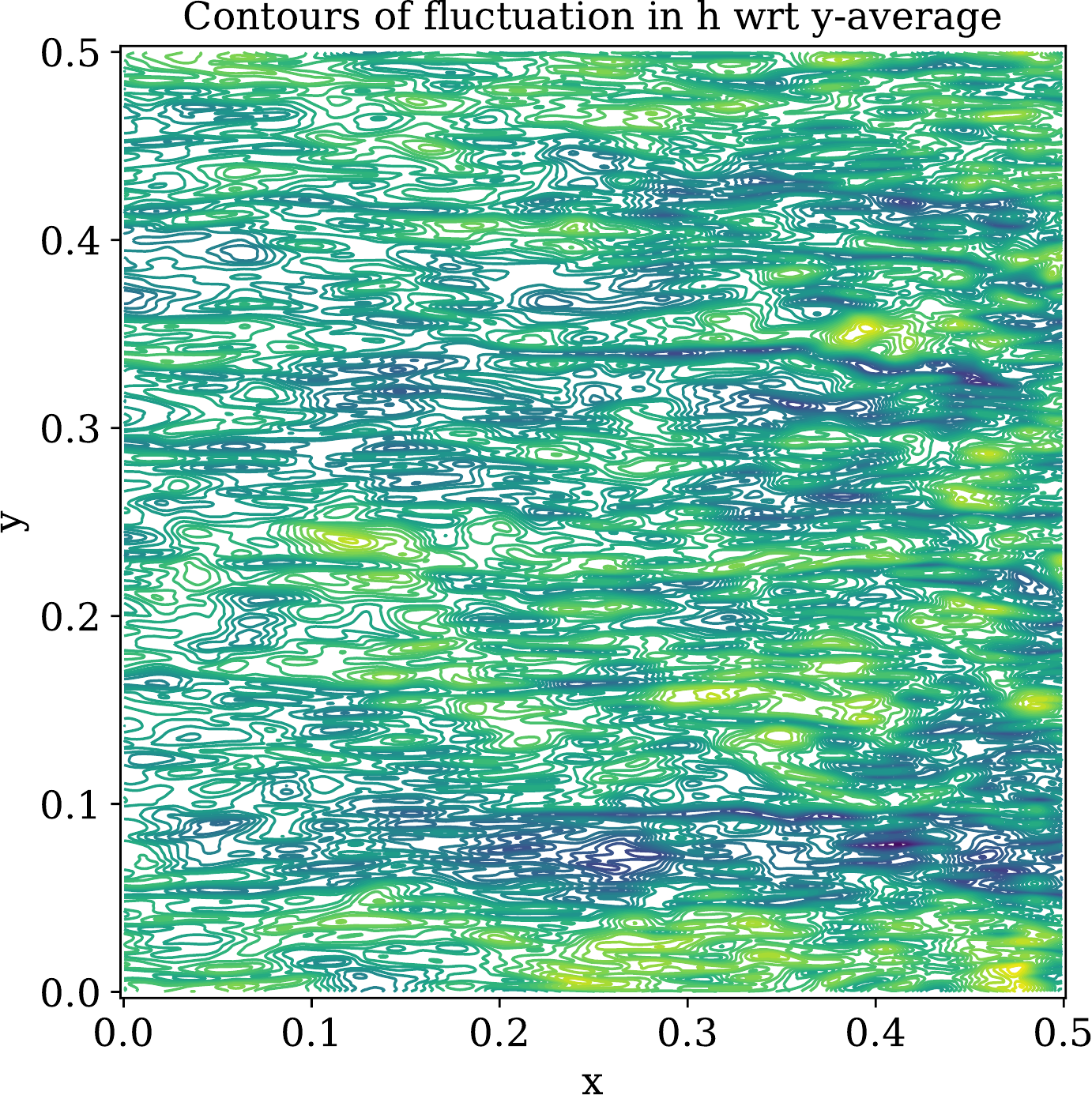}
\end{tabular}
\end{center}
\caption{2-D roll wave problem using second order scheme on $2080 \times 800$ mesh at $t=36$. Zoomed view of contour plots of fluctuations in $h$ wrt y-average. Left: HLLC3, Right: HLLC5.}
\label{fig:roll2d5}
\end{figure}

\begin{figure}
\begin{center}
\includegraphics[width=0.98\textwidth]{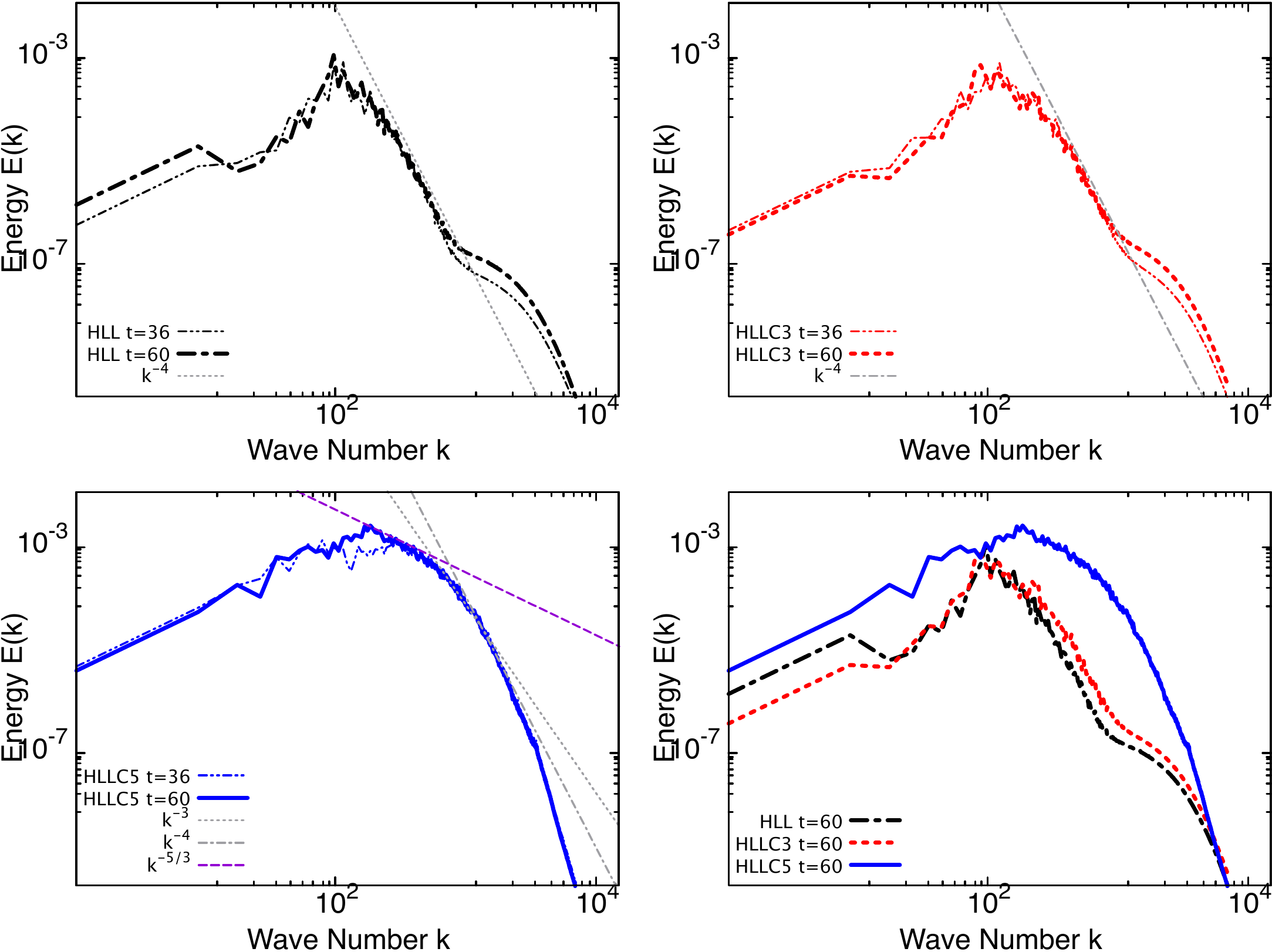}
\end{center}
\caption{Kinetic energy spectrum of fluctuations around the $y$-average for 2-D roll wave problem on $2080 \times 800$ mesh.}
\label{fig:tkespec}
\end{figure}

\paragraph{Remark}
\reva{
The shape observed with the HLLC5 solver on coarse grid have some similarities with the numerical instability named ``carbuncle phenomena" that pollute the bow-shock upstream in hypersonic flows. According to~\cite{Elling2009},  ``It is not known which numerical schemes are affected in which circumstances, what causes carbuncles to appear and whether carbuncles are purely numerical artefacts or rather features of a continuum equation or model". And then to conclude that "carbuncles can, in some circumstances, be valid vanishing viscosity limits. Trying to suppress them is making a physical assumption that may be false". Therefore, we cannot  get a definitive answer in this paper and further investigations are needed understand the origin of this behaviour in the present model. Moreover, in the context of shear shallow water model proposed here, the equations are genuinely non-conservative and the classical analysis of shock stability needs to be reconsidered. Nevertheless, results obtained here are consistent with previous papers. The turbulent pattern is shared by the different numerical schemes, including the HLL scheme, which is not expected to exhibit carbuncle-type solutions. 
}
\section{Summary and conclusions}
\label{sec:sum}
The present work deals with the shear shallow water model which is a higher order version of the classical shallow water model since it includes effects of vertical shear. The model has a non-conservative structure which poses difficulties in its numerical solution. Previous works have developed Riemann solvers by splitting the model into two parts, while in this work, we develop a unified Riemann solver using a more fundamental form of the equations which arises directly from the depth averaging process and is similar to the 10-moment equations in gas dynamics. By using a linear path in conservative variables, we develop path conservative Riemann solvers for this model. In particular, we develop HLL, 3-wave and 5-wave HLLC solvers and apply them to a set of test problems in 1-D and 2-D. A second order version of the scheme is developed using a predictor-corrector approach and the order of accuracy is checked numerically, while existing works have only demonstrated first order results. Among the three Riemann solvers, the results demonstrate that the HLLC5 solver is the most accurate in terms of resolving all the waves that arise in the solution. The 1-D results compare well with results from previous 5-wave solver available in the literature except in the 1-D dam break problem, where we see some differences in the shock. Such differences can be expected due to different jump conditions in our model compared to previous approaches. The 1-D roll wave problem develops the roll wave profile and compares well with Brock's experimental results. The 2-D roll wave problem shows significant differences in the performance of the Riemann solvers. Only the HLLC5 solver develops a turbulent-like solution which is similar to previous results in the literature and such a solution seems to correspond to real flows, see e.g.~\cite{Aranda2016}. The solution has a well defined $y$-average with fluctuations super-imposed on top of this. Moreover, the fluctuations have a kinetic energy spectrum with well-defined scaling laws being visible, though more refined calculations are necessary to make definite statements on the precise structure of the solution. The higher order scheme using the HLL-type solvers developed in this work is thus very promising for further studies on shear shallow water flows.
\section*{Acknowledgments}
The work reported in this paper was started when the first author was visiting Universit\'e C\^ote d'Azur, Nice, France, under a visiting professor grant. This work has been supported by the French government, through the {\em UCAJEDI Investments in the Future} project managed by the National Research Agency (ANR) with the reference number ANR-15-IDEX-01. The work of the first author is supported by the Department of Atomic Energy,  Government of India, under project no.  12-R\&D-TFR-5.01-0520.
\appendix
\section{Derivation of the SSW model}
\label{sec:avg}
The gravitational force is assumed to act in the negative $x_3$ direction. Let $x_3 = \xi(x_1,x_2,t)$ denote the height of the free-surface and $x_3 = b(x_1,x_2,t)$ denote the location of the bottom surface, see Figure~(\ref{fig:shallow}). Then the depth of water is $h = \xi - b$. The flow is assumed to be governed by the incompressible Euler equations which are given by\footnote{In this section, the subscripts like $i,j,k$ take values in $\{1,2,3\}$ and repeated subscripts indicate summation.}
\begin{align*}
\df{u_k}{x_k} =& 0 \\
\df{u_i}{t} + u_k \df{u_i}{x_k} + \frac{1}{\rho}\df{p}{x_i} =& - g \delta_{i3}
\end{align*}
where $(u_1,u_2,u_3)$ is the 3-D velocity field, $p$ is the pressure and we assume that the density $\rho$  is constant. A material point on the free surface $F(x_1,x_2,x_3,t) = x_3 - \xi(x_1,\xi_2,t) = 0$ remains on the free surface which means that\footnote{Greek subscripts like $\alpha,\beta, \gamma$ will takes values in $\{1,2\}$ only and repeated subscripts indicate summation.}
\begin{equation}
\label{eq:freesurf}
\dd{F}{t}=0 \quad\implies\quad \df{\xi}{t} + u_\alpha \df{\xi}{x_\alpha} - u_3 = 0 \qquad \textrm{on} \quad x_3 = \xi
\end{equation}
Similarly, for the bottom surface, we get
\begin{equation}
\label{eq:botsurf}
\df{b}{t} + u_\alpha \df{b}{x_\alpha} - u_3 = 0 \qquad \textrm{on} \quad x_3 = b
\end{equation}
Define the depth average $\davg{\phi} = \davg{\phi}(x_1,x_2,t)$ of any quantity $\phi = \phi(x_1,x_2,x_3,t)$ by
\[
\davg{\phi}(x_1,x_2,t) = \frac{1}{h} \int_{b}^\xi \phi(x_1,x_2,x_3,t) \ud x_3
\]
The fluctuation with respect to the average value is
\[
\phi' = \phi - \davg{\phi}
\]
and clearly
\[
\davg{\phi'} = 0
\]
Moreover, we note the following identities
\begin{equation}
\label{eq:avgx}
\int_b^\xi \df{\phi}{x_\alpha} \ud x_3 = \df{(h\davg{\phi})}{x_\alpha} - \df{\xi}{x_\alpha} \phi|_{x_3 = \xi} + \df{\xi}{x_\alpha} \phi|_{x_3 = b}
\end{equation}
\begin{equation}
\label{eq:avgt}
\int_b^\xi \df{\phi}{t} \ud x_3 = \df{(h\davg{\phi})}{t} - \df{\xi}{t} \phi|_{x_3 = \xi} + \df{\xi}{t} \phi|_{x_3 = b}
\end{equation}
Integrating the continuity equation over the depth of water and using~\eqref{eq:avgx} yields
\[
\df{}{x_\alpha}(h \davg{u}_\alpha) - u_\alpha(x_1,x_2,\xi,t)\df{\xi}{x_\alpha} + u_\alpha(x_1,x_2,b,t) \df{b}{x_\alpha} + u_3(x_1,x_2,\xi,t) - u_3(x_1,x_2,b,t) = 0
\]
and using \eqref{eq:freesurf},~\eqref{eq:botsurf} we get
\begin{equation}
\label{eq:ssw1}
\df{h}{t} + \df{}{x_\alpha}(h \davg{u}_\alpha) = 0
\end{equation}
Let $H,L$ denote the vertical and horizontal length scales; in the shallow water approximation, the vertical scale is much smaller than the horizontal scale (long wave approximation), i.e.,
\[
\varepsilon = \frac{H}{L} \ll 1
\]
Let $U$ denote the horizontal velocity scale; then the continuity equation implies that $u_3 = O(\varepsilon U)$. Using time and pressure scales as $L/U$ and $\rho U^2$, the $x_3$-momentum equation can be non-dimensionalized as 
\[
\varepsilon^2 \frac{D u_3'}{D t'} + \df{p'}{x_3'} = - \frac{1}{\Fr^2} \rho'
\]
where $D/Dt$ denotes material derivative, the prime quantities are non-dimensional and $\Fr = U/\sqrt{g H}$ is Froude number. Ignoring terms of $O(\varepsilon^2)$, we obtain the hydrostatic approximation
\[
\df{p}{x_3} = - \rho g \quad\implies\quad p - p_a = - \rho g (x_3 - \xi)
\]
where $p_a$ is the atmospheric pressure at the free surface which may be taken to be constant. Hence
\begin{equation}
\label{eq:hydroxy}
\df{p}{x_\alpha} =  \rho g \df{\xi}{x_\alpha}
\end{equation}
Now the horizontal momentum equation takes the form
\[
\df{u_\alpha}{t} + 2\df{K_{\alpha\beta}}{x_\beta} + \df{}{x_3}(u_\alpha u_3) + \rho g \df{\xi}{x_\alpha} = 0 \qquad \textrm{where} \qquad K_{\alpha\beta} = \half u_\alpha u_\beta
\]
Averaging this equation and using~\eqref{eq:freesurf},~\eqref{eq:botsurf},~\eqref{eq:avgx},~\eqref{eq:avgt}, we obtain
\begin{equation}
\label{eq:ssw2}
\df{(h \davg{u}_\alpha)}{t} + 2\df{\davg{K}_{\alpha\beta}}{x_\beta} + \df{}{x_\alpha}\left(\half \rho g h^2\right) = -\rho g h \df{b}{x_\alpha}
\end{equation}
where $\davg{K}_{\alpha\beta} = \half \davg{u}_\alpha \davg{u}_\beta + \half \davg{u'_\alpha u'_\beta}$ is not completely known to us and hence the set of equations~\eqref{eq:ssw1},~\eqref{eq:ssw2} does not form a closed system. If we set the second order velocity fluctuation term $\davg{u'_\alpha u'_\beta}$ to zero, then the equations are closed and we obtain the classical shallow water model. This simplification can be justified if the vertical shear is small, i.e., $\df{u_\alpha}{x_3} = O(\varepsilon^m)$ with $m \ge 1$~\cite{Teshukov2007}.

If the vertical shear is not small, i.e. $m < 1$, then it is not justifiable to ignore the second order velocity fluctuation terms, in which case these terms must be retained and we must derive additional equations to model them. Starting from the momentum equation, we can derive the following set of equations,
\[
2\df{K_{ij}}{t} + 2\df{}{x_k}(K_{ij} u_k) + \frac{u_i}{\rho} \df{p}{x_j} + \frac{u_j}{\rho} \df{p}{x_i} = -g( u_j \delta_{i3} + u_i \delta_{j3})
\]
and hence, using the hydrostatic condition~\eqref{eq:hydroxy}, we obtain
\[
2\df{K_{\alpha\beta}}{t} + 2\df{}{x_\gamma}(K_{\alpha\beta} u_\gamma) + 2\df{}{x_3}(K_{\alpha\beta} u_3) + g u_\alpha \df{\xi}{x_\beta} + g u_\beta \df{\xi}{x_\alpha} = 0
\]
Averaging the above equation over the depth and using~\eqref{eq:freesurf},~\eqref{eq:botsurf},~\eqref{eq:avgx},~\eqref{eq:avgt} yields
\begin{equation}
\df{(h \davg{K}_{\alpha\beta})}{t} + \df{}{x_\gamma}(h \davg{K_{\alpha\beta} u_\gamma}) + \half g h \davg{u}_\alpha \df{h}{x_\beta} + \half g h \davg{u}_\beta \df{h}{x_\alpha} = -\half g h \davg{u}_\alpha \df{b}{x_\beta} - \half g h \davg{u}_\beta \df{b}{x_\alpha}
\label{eq:ssw3}
\end{equation}
Let us identify $\p_{\alpha\beta} = \davg{u'_\alpha u'_\beta}$ and $E_{\alpha\beta} = h \davg{K}_{\alpha\beta}$. The remaining average term in the previous equation can be written as
\[
\davg{K_{\alpha\beta} u_\gamma} = \half \davg{u_\alpha u_\beta u_\gamma} = \davg{K}_{\alpha\beta} \davg{u}_\gamma + \half \davg{u}_\alpha \p_{\beta\gamma} + \half \davg{u}_\beta \p_{\alpha\gamma} + \half \davg{u'_\alpha u'_\beta u'_\gamma}
\]
The third order fluctuations are of $O(\varepsilon^{3 m})$ which are smaller than the second order fluctuations which are of $O(\varepsilon^{2 m})$~\cite{Teshukov2007}. This allows us to ignore the third order fluctuations in equation~\eqref{eq:ssw3} and we obtain a closed set of equations. The set of equations~\eqref{eq:ssw1},~\eqref{eq:ssw2}~,\eqref{eq:ssw3} is exactly the model given in~\eqref{eq:tenmom}.
\section{Solution of semi-implicit scheme}
\label{sec:imp}
The semi-implicit schemes take the form
\[
\con^{n+1} - \theta \dt \s(\con^{n+1}) = \tilde{\con}^{n+1}
\]
where $\tilde{\con}^{n+1}$ is the explicit update without the source term, $\theta=1$ for the first order scheme and $\theta=\half$ for the second order scheme. There is no source term in the continuity equation and we obtain $h^{n+1} = \tilde{\con}^{n+1}_1$. The momentum equations have the form\footnote{For simplicity of notation, we drop the superscript $n+1$ on some of the quantities in the following equations.}
\begin{align*}
m_1 + \frac{\theta C_f \dt}{h^2} m_1 \sqrt{m_1^2 + m_2^2} =& \tilde{\con}_2^{n+1} - \theta g h^{n+1} \df{b}{x} \dt =: a_1 \\
m_2 + \frac{\theta C_f \dt}{h^2} m_2 \sqrt{m_1^2 + m_2^2} =& \tilde{\con}_3^{n+1} - \theta g h^{n+1} \df{b}{y} \dt =: a_2
\end{align*}
Squaring and adding the two equations we get
\[
(m_1^2 + m_2^2) \left(1 + c \sqrt{m_1^2 + m_2^2}\right)^2 = a_1^2 + a_2^2, \qquad c = \frac{\theta C_f \dt}{h^2}
\]
This leads to a quartic equation for $m = \sqrt{m_1^2 + m_2^2}$
\[
m^2 (1 + cm)^2 - (a_1^2 + a_2^2) = 0 \quad\implies\quad c m^2 + m - \sqrt{a_1^2 + a_2^2} = 0
\]
where we choose the positive square root. The positive solution of the quadratic equation is
\[
m = \frac{-1 + \sqrt{1 + 4c\sqrt{a_1^2 + a_2^2}}}{2c}
\]
Then the momentum components are given by
\[
m_1 = \frac{a_1}{1 + cm}, \qquad m_2 = \frac{a_2}{1 + cm}
\]
From the $E_{11}, E_{22}$ equations we obtain
\begin{align}
\label{eq:p11}
\half h \p_{11} + \alpha |\vel|^3 \theta\dt \p_{11} =& S_{11} := \tilde{\con}_4^{n+1} - \half h^{n+1} (v_1^{n+1})^2  - \theta \dt \left[ g h^{n+1} v_1^{n+1} \df{b}{x} + C_f |\vel^{n+1}| (v_1^{n+1})^2 \right] \\
\label{eq:p22}
\half h \p_{22} + \alpha |\vel|^3 \theta\dt \p_{22} =& S_{22} := \tilde{\con}_6^{n+1} - \half h^{n+1} (v_2^{n+1})^2 - \theta \dt \left[ g h^{n+1} v_2^{n+1} \df{b}{y} + C_f |\vel^{n+1}| (v_2^{n+1})^2 \right] 
\end{align}
Adding the two equations, we get
\[
f(T) := \half h T + \alpha |\vel|^3 \theta\dt T - ( S_{11} + S_{22} ) = 0
\]
which contains only the trace $T$ as unknown. This is a non-linear equation since $\alpha = \alpha(h,T)$ is non-linear and a necessary condition for existence of positive solution is that $S_{11} + S_{22} > 0$. For $0 < T \le \phi h^2$, $f(T) = \half h T  - ( S_{11} + S_{22} )$ is negative for $T$ sufficiently small; for $T \gg \phi h^2$ we have $f(T) > 0$. Moreover $f(T)$ is a continuous, monotonically increasing function in $(0,\infty)$; this is clear for $T \in (0, \phi h^2]$ and for $T > \phi h^2$
\[
f(T)  = \half h T + \frac{C_r |\vel|^3  \theta \dt}{T} (T - \phi h^2) - (S_{11} + S_{22}), \qquad f'(T) = \half h + \frac{C_r |\vel|^3  \theta \dt \phi h^2}{T^2} > 0
\]
hence $f(T)=0$ has  a unique positive solution. To solve this equation we first compute, $\hat{T}$ from
\[
\half h \hat{T} = S_{11} + S_{22}
\]
If $\hat{T} \le \phi h^2$ then the solution is $\hat{T}$ itself; otherwise we solve the quadratic equation $f(T)=0$ and choose the unique positive root in $(\phi h^2, \infty)$ as the solution. Once $T$ is obtained, we get $\p_{11}, \p_{22}$ from~\eqref{eq:p11},~\eqref{eq:p22}; and also $\p_{12}$ can be determined from the $E_{12}$ equation which is of the form
\begin{align*}
\half h \p_{12} + \alpha |\vel|^3 \theta\dt \p_{12} = S_{12} := & \ \tilde{\con}_5^{n+1} - \half h^{n+1} v_1^{n+1} v_2^{n+1}  \\
& - \theta \dt \left[ \half g h^{n+1} v_2^{n+1} \df{b}{x} + \half g h^{n+1} v_1^{n+1} \df{b}{y} + C_f |\vel^{n+1}| v_1^{n+1} v_2^{n+1} \right]
\end{align*}
Since we know the tensor $\p$, we can update the energy tensor $E$ which completes the update of all the quantities.
\bibliographystyle{plain}
\bibliography{bbibtex}

\begin{thebibliography}{10}

\bibitem{Abgrall2010}
R{\'e}mi Abgrall and Smadar Karni.
\newblock A comment on the computation of non-conservative products.
\newblock {\em Journal of Computational Physics}, 229(8):2759--2763, April
  2010.

\bibitem{Aranda2016}
Alfredo Aranda, Nicol{\'a}s Amigo, Christian Ihle, and Aldo Tamburrino.
\newblock Digital image correlation applied to the calculation of the
  out-of-plane deformation induced by the formation of roll waves in a
  non-{{Newtonian}} fluid.
\newblock {\em Optical Engineering}, 55(6):064101, June 2016.

\bibitem{Berthon2002}
C.~Berthon, F.~Coquel, J.M. H{\'e}rard, and M.~Uhlmann.
\newblock An approximate solution of the {{Riemann}} problem for a realisable
  second-moment turbulent closure.
\newblock {\em Shock Waves}, 11(4):245--269, January 2002.

\bibitem{Berthon2006}
Christophe Berthon.
\newblock Numerical approximations of the 10-moment {{Gaussian}} closure.
\newblock {\em Mathematics of Computation}, 75(256):1809--1831, June 2006.

\bibitem{Bhole2019}
Ashish Bhole, Boniface Nkonga, Sergey Gavrilyuk, and Kseniya Ivanova.
\newblock Fluctuation splitting {{Riemann}} solver for a non-conservative
  modeling of shear shallow water flow.
\newblock {\em Journal of Computational Physics}, 392:205--226, September 2019.

\bibitem{Boffetta2012}
Guido Boffetta and Robert~E. Ecke.
\newblock Two-{{Dimensional Turbulence}}.
\newblock {\em Annual Review of Fluid Mechanics}, 44(1):427--451, January 2012.

\bibitem{Brock1969}
Richard~R. Brock.
\newblock Development of {{Roll}}-{{Wave Trains}} in {{Open Channels}}.
\newblock {\em Journal of the Hydraulics Division}, 95(4):1401--1428, 1969.

\bibitem{Brock1970}
Richard~R. Brock.
\newblock Periodic {{Permanent Roll Waves}}.
\newblock {\em Journal of the Hydraulics Division}, 96(12):2565--2580, 1970.

\bibitem{DalMaso1995}
Gianni Dal~Maso, Philippe~G. Lefloch, and Francois Murat.
\newblock Definition and weak stability of nonconservative products.
\newblock {\em J. Math. Pures Appl.}, 74:483--548, 1995.

\bibitem{Einfeldt1988}
Bernd Einfeldt.
\newblock On {{Godunov}}-{{Type Methods}} for {{Gas Dynamics}}.
\newblock {\em SIAM Journal on Numerical Analysis}, 25(2):294--318, April 1988.

\bibitem{Elling2009}
Volker Elling.
\newblock The carbuncle phenomenon is incurable.
\newblock {\em Acta Mathematica Scientia}, 29(6):1647--1656, November 2009.

\bibitem{Gavrilyuk2018}
S.~Gavrilyuk, K.~Ivanova, and N.~Favrie.
\newblock Multi-dimensional shear shallow water flows: {{Problems}} and
  solutions.
\newblock {\em Journal of Computational Physics}, 366:252--280, August 2018.

\bibitem{Godlewski1996}
Edwige Godlewski and Pierre-Arnaud Raviart.
\newblock {\em Numerical {{Approximation}} of {{Hyperbolic Systems}} of
  {{Conservation Laws}}}, volume 118 of {\em Applied {{Mathematical
  Sciences}}}.
\newblock {Springer New York}, {New York, NY}, 1996.

\bibitem{Gosse2001}
Laurent Gosse.
\newblock A well-balanced scheme using non-conservative products designed for
  hyperbolic systems of conservation laws with source terms.
\newblock {\em Mathematical Models and Methods in Applied Sciences},
  11(02):339--365, March 2001.

\bibitem{Harten1983a}
Amiram Harten, Peter~D. Lax, and Bram van Leer.
\newblock On {{Upstream Differencing}} and {{Godunov}}-{{Type Schemes}} for
  {{Hyperbolic Conservation Laws}}.
\newblock {\em SIAM Review}, 25(1):35--61, January 1983.

\bibitem{Ivanova2017}
K.A. Ivanova, S.L. Gavrilyuk, B.~Nkonga, and G.L. Richard.
\newblock Formation and coarsening of roll-waves in shear shallow water flows
  down an inclined rectangular channel.
\newblock {\em Computers \& Fluids}, 159:189--203, December 2017.

\bibitem{Levermore1998}
C.~David Levermore and William~J. Morokoff.
\newblock The {{Gaussian Moment Closure}} for {{Gas Dynamics}}.
\newblock {\em SIAM Journal on Applied Mathematics}, 59(1):72--96, January
  1998.

\bibitem{Pares2006}
Carlos Par{\'e}s.
\newblock Numerical methods for nonconservative hyperbolic systems: A
  theoretical framework.
\newblock {\em SIAM Journal on Numerical Analysis}, 44(1):300--321, January
  2006.

\bibitem{Richard2012}
G.~L. Richard and S.~L. Gavrilyuk.
\newblock A new model of roll waves: Comparison with {{Brock}}'s experiments.
\newblock {\em Journal of Fluid Mechanics}, 698:374--405, May 2012.

\bibitem{Richard2013}
G.~L. Richard and S.~L. Gavrilyuk.
\newblock The classical hydraulic jump in a model of shear shallow-water flows.
\newblock {\em Journal of Fluid Mechanics}, 725:492--521, June 2013.

\bibitem{Teshukov2007}
V.~M. Teshukov.
\newblock Gas-dynamic analogy for vortex free-boundary flows.
\newblock {\em Journal of Applied Mechanics and Technical Physics},
  48(3):303--309, May 2007.

\bibitem{Toro1994}
E.~F. Toro, M.~Spruce, and W.~Speares.
\newblock Restoration of the contact surface in the {{HLL}}-{{Riemann}} solver.
\newblock {\em Shock Waves}, 4(1):25--34, July 1994.

\bibitem{VanLeer1997}
Bram Van~Leer.
\newblock On {{The Relation Between The Upwind}}-{{Differencing Schemes Of
  Godunov}}, {{Engquist}}\textemdash{{Osher}} and {{Roe}}.
\newblock In M.~Yousuff Hussaini, Bram {van Leer}, and John Van~Rosendale,
  editors, {\em Upwind and {{High}}-{{Resolution Schemes}}}, pages 33--52.
  {Springer Berlin Heidelberg}, {Berlin, Heidelberg}, 1997.

\bibitem{volpert1967}
A.~I. Volpert.
\newblock The spaces {{BV}} and quasilinear equations.
\newblock {\em Mathematics of the USSR-Sbornik}, 2(2):225--267, February 1967.

\bibitem{Whitham1999}
G.~B. Whitham.
\newblock {\em Linear and {{Nonlinear Waves}}}.
\newblock {John Wiley \& Sons, Inc.}, {Hoboken, NJ, USA}, June 1999.

\end{thebibliography}
\end{document}